\documentclass[10pt]{article}
\usepackage{amssymb}
\usepackage{amsmath}
\usepackage{amsthm}
\usepackage{graphicx}
\usepackage[utf8]{inputenc}
\usepackage{amsfonts}
\usepackage{graphicx}
\usepackage{xcolor}
\usepackage{breqn}
\usepackage[utf8]{inputenc}
\usepackage{float}
\usepackage{hyperref}
\usepackage{cite}
\def\bc{\begin{center}}
	\def\ec{\end{center}}

\def\s2c{\vskip 2cm}
\def\bt{\begin{Theorem}}
	\def\et{\end{Theorem}}
\def\bd{\begin{Definition}}
	\def\ed{\end{Definition}}
\def\bl{\begin{Lemma}}
	\def\el{\end{Lemma}}
\def\be{\begin{Example}}
	\def\ee{\end{Example}}
\def\bcor{\begin{Corollary}}
	\def\ecor{\end{Corollary}}
\def\br{\begin{Remark}}
	\def\er{\end{Remark}}

\newtheorem{Lemma}{Lemma}[section]
\newtheorem{Theorem}[Lemma]{Theorem}
\newtheorem{Definition}[Lemma]{Definition}

\newtheorem{Corollary}[Lemma]{Corollary}
\newtheorem{Remark}[Lemma]{Remark}
\usepackage{graphicx, epsfig, graphics, graphpap}

\usepackage{breqn,epstopdf}
\newtheorem{theorem}{Theorem}[section]

\newtheorem{lemma}[theorem]{Lemma}
\newtheorem{remark}{Remark}

\date{}
\title{An efficient explicit jump HOC immersed interface approach for transient incompressible viscous flows}
\author {Raghav Singhal$^{1}$, Jiten C. Kalita$^{1}$
\\
Deparment of Mathematics, Indian Institute of Technology Guwahati, Assam 781039, India
	\\author's email :  raghav2016@iitg.ac.in \quad and jiten@iitg.ac.in}
	
\begin{document}
	\maketitle
	\author
	\noindent {\bf Abstract} :
In the present work, we propose a novel hybrid explicit jump immersed interface approach in conjunction with a higher order compact (HOC) scheme for simulating transient complex flows governed by the streamfunction-vorticity ($\psi$-$\zeta$) formulation of the Navier-Stokes (N-S) equations for incompressible viscous flows. A new strategy has been adopted for the jump conditions at the irregular points across the interface using Lagrangian interpolation on a Cartesian grid. This approach, which starts with the discretization of parabolic equations with discontinuities in the solutions, source terms and the coefficients across the interface, can easily be accommodated into simulating flow past bluff bodies immersed in the flow. The superiority of the approach is reflected by the reduced magnitude and faster decay of the errors in comparison to other existing methods. It is seen to handle several fluid flow problems having practical implications in the real world very efficiently, which involves flows involving multiple and moving bodies. This includes the flow past a stationary circular and a twenty-four edge cactus cylinder, flows past two tandem cylinders, where in one situation both are fixed and in another, one of them is oscillating transversely with variable amplitude in time. To the best of our knowledge, the last two examples have been tackled for the first time by such an approach employing the $\psi$-$\zeta$ formulation in finite difference set-up. The extreme closeness of our computed solutions with the existing numerical and experimental results exemplifies the accuracy and the robustness of the proposed approach.
	
\section{Introduction}Parabolic partial differential equations (pde) with discontinuous coefficients play significant roles in the fields of electrostatics, porous media, multiphase flows, material science, underwater acoustics, biology (blood flow models) and several other fields. For example, one may consider the case of wave propagation in a composite or irregular medium with different material properties, which has numerous engineering applications. Besides, the parabolic equations can be suitably reconstructed in the shape of the unsteady Navier-Stokes (N-S) equations that models incompressible viscous flows. As such, devising computationally efficient numerical algorithms and obtaining highly accurate solutions of the N-S equations is one of the primary goals of Computational Fluid Dynamics (CFD). In order to deal with bodies of complicated geometries or moving bodies immersed in fluid flows in Cartesian grids, immersed interface methods have been of great efficacy in the past three decades. Recently, we proposed a new higher-order accurate finite difference explicit jump Immersed Interface Method (HEJIIM) for solving two-dimensional elliptic problems \cite{singhalkalita}. In the current work, we extend this idea to transient problems, viz., to parabolic problems with singular source and discontinuous coefficients irregular regions on a compact Cartesian mesh.    

In most of the parabolic problems with discontinuities across some interface, they are mostly expressed in terms of natural jump conditions in the dependent variable given by $[u]=0$ and its normal derivative $[\beta u_{\textbf{n}}]=0$ across the interface. For example, one may cite the process of conductive heat transfer over mixed media, which has been extensively investigated in the existing literature. Although many of the parabolic equations possess analytical solutions, even for problems having homogeneous jump conditions, the discontinuity present in the diffusion coefficient doesn't allow one to have the analytical solution easily. To overcome these issues, numerical approaches become handy techniques for such types of problems. However, as the solution is not smooth over the whole physical and computational domains, standard numerical algorithms cannot be applied to attain accurate solutions. One must conceive specific numerical procedures in the neighbourhood of these discontinuities to achieve reasonable approximation of the solutions thereat.

Immersed Boundary Method (IBM) was first introduced in  1972 by Peskin \cite{peskin1972flow} to simulate blood and cardiac mechanics. The main feature of Peskin's approach was that the entire simulation was accomplished on a Cartesian grid which did not conform to the shape of the heart. The method approximated the prescribed boundary conditions of the immersed objects  by incorporating the forcing term in the form of the Dirac delta function to the right side of the N-S equations. The approach was restricted to problems having continuous solutions only, and was first-order accurate. Superior to Peskin's formulation, Li et al. \cite{leveque1994immersed} devised an accurate second-order method that deals with problems with singular source terms and discontinuous coefficients on the irregular domain, namely Immersed Interface Method (IIM). They incorporated the interface jump conditions in the solution and the flux at the point of discontinuity. IBMs, IIMs, and their numerous variants have become increasingly popular and highly relevant to numerically solving the initial boundary-value problem on irregular domains. 

Several remarkably designed numerical approaches have already been developed in the literature to solve the parabolic interface problems by incorporating the jump conditions into the discretization process \cite{li2018multiscale, li2017matched, papac2010efficient, wei2018spatially, zhao2015matched}. Notwithstanding, most of these involve finite element and finite volume methods where the interface is captured by body-fitted approach \cite{attanayake2011convergence, chen1998finite, song2018symmetric}, leading to computationally expensive grid generation.  In order to avoid this, Li et al.  introduced second-order accurate immersed interface method for moving interface problem on Cartesian mesh \cite{li1997immersed}, by using generalized Taylor series expansions to modify the standard finite difference (FD) discretization on the irregular points to recover the loss of accuracy at the interface. Adams et al. \cite{adams2002immersed}  introduced the second-order maximum principle immersed interface method to solve the linear parabolic equation. First-order derivatives are approximated by an explicit scheme and the diffusion part by Cranck Nicolson. \cite{bouchon2010immersed} Bouchon and Peichl presented a method to parabolic equations with mixed boundary conditions, where they applied immersed interface algorithm to discretize Neumann condition and Shortley-Weller approximation for the Dirichlet condition. The most significant and inherent aspect of these approaches lie in their clarity in achieving the solution on Cartesian mesh, which can be generated very speedily and allows users to simulate flows containing moving objects with complex geometries with ease. On the other hand, body-conformal mesh requires the generation of a new mesh at each time-step which may impact the robustness, accuracy and computational cost for similar problems detrimentally. 

In the current work, we propose a new higher-order compact finite difference Immersed Interface Method  for solving two-dimensional parabolic problems, more specifically for transient problems involving bluff bodies immersed in incompressible viscous flows on Cartesian mesh. Such problems are governed by the unsteady N-S equations which are parabolic in nature with singular source and discontinuous coefficients in irregular domains. $\psi$-$\zeta$ form of the N-S equations has been utilized for this purpose as in  \cite{calhoun2002cartesian, linnick2005high, russell2003cartesian}. Note that Calhoun \cite{calhoun2002cartesian} presented a second-order finite volume approach in an unstable region and imposed no-slip flow condition to find the vorticity sources, while Linnack and Fasel \cite{linnick2005high} introduced a fourth-order compact difference scheme based on Weigmanm \cite{wiegmann2000explicit} approach. Russel and Wang \cite{russell2003cartesian} satisfied the no-penetration condition in the streamfucntion by superimposing a homogenous solution to Poisson's equation for moving boundaries and a no-slip condition for the surface vorticity of the objects. Similar to others IIMs approaches, to ensure the accuracy of the numerical solutions, they adopted special strategies close to the embedded boundary, resulting in the loss of compactness of stencil. In contrast, the proposed scheme maintains its compactness on a nine point stencil at both the regular and irregular points. In order to treat the jump across the interface, we modified the HEJIIM \cite{singhalkalita} in such a way that at each time step, the scheme maintains fourth order spatial accuracy throughout the whole computational domain. 

Using the proposed scheme, firstly we solve one problem with circular interface in a rectangular region having analytical solutions. Then we simulate flow past stationary as well as moving bluff bodies immersed in fluids governed by the N-S equations.  Our simulations include flow situations involving multiple and moving bodies as well. For the problem having analytical solution, our results are excellent match with the analytical ones and for the fluid flow problems, our simulations are extremely close to the experimental and available numerical results.

The paper is organized in the following way. In section 2, we detail the development of the proposed scheme, section 3 discusses issues involving the steamfunction-vorticity formulations along with brief descriptions of the associated fluid dynamic forces and solution of the algebraic systems, section 4 deals with the numerical examples and finally in conclusion, we summarize our achievements.
 
\section{Mathematical Formulation} 
\begin{figure}[!t]
\minipage{0.33\textwidth}
  \includegraphics[width=\linewidth]{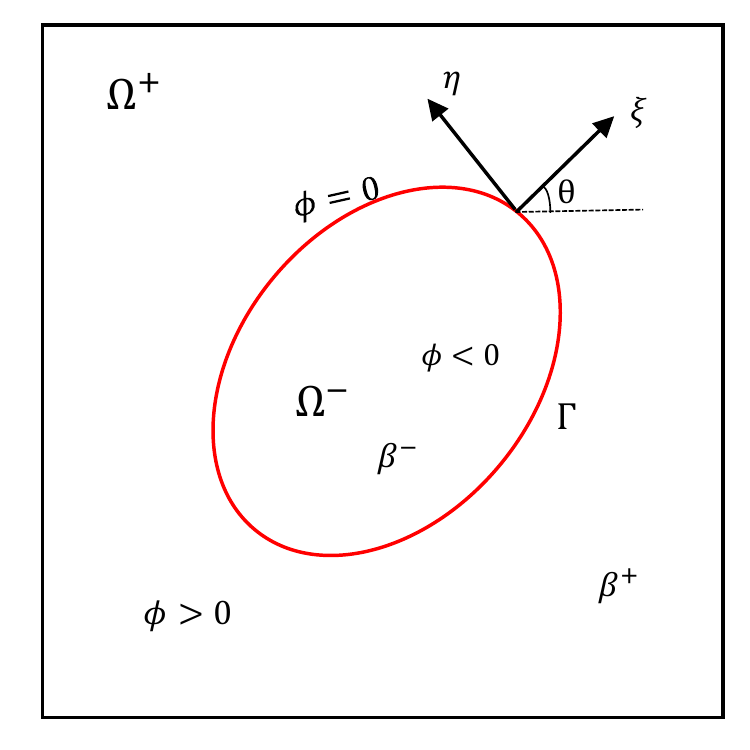}
  \begin{center}
   (a)
  \end{center}
\endminipage\hfill
\minipage{0.33\textwidth}
  \includegraphics[width=\linewidth]{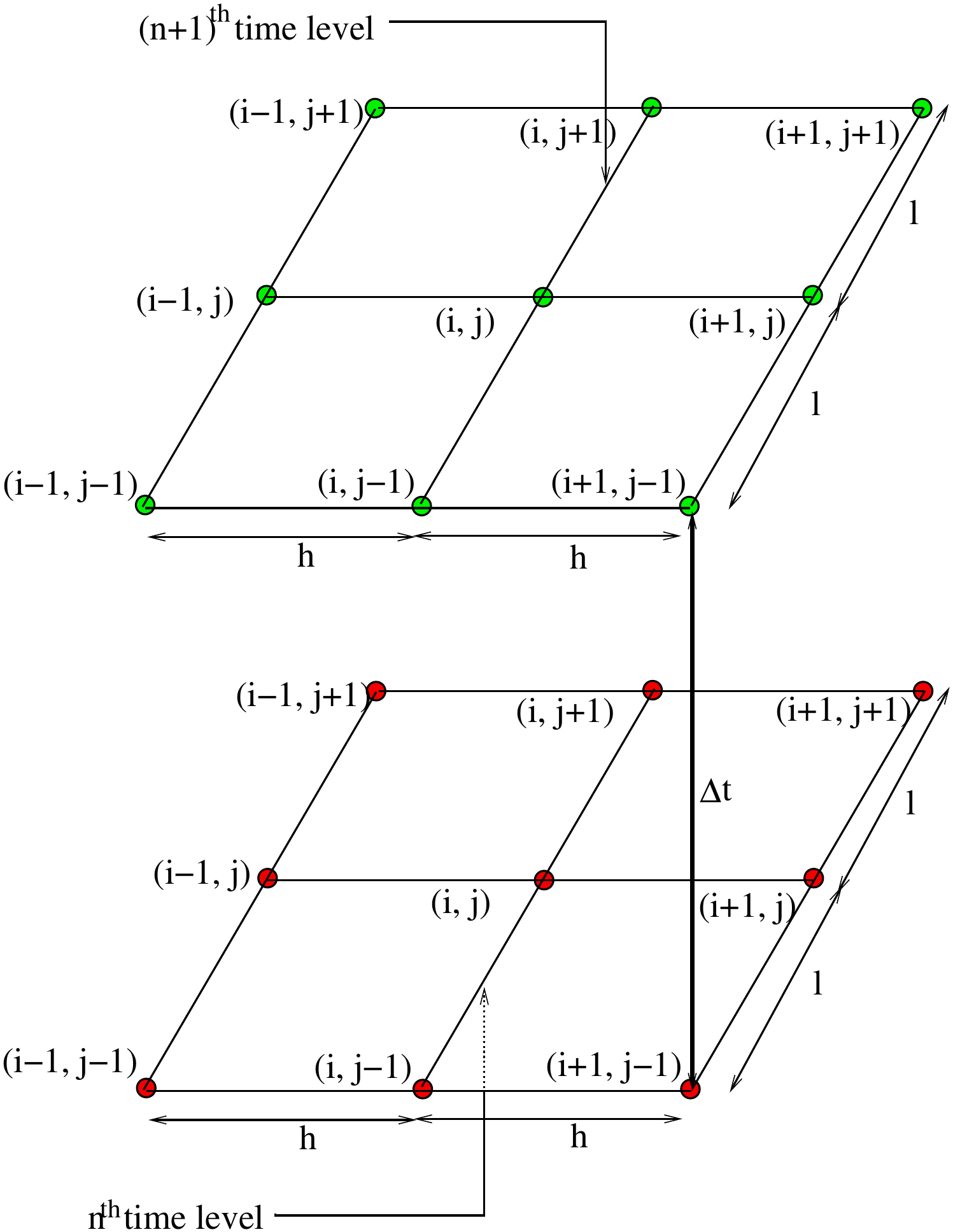}
  \begin{center}
   (b)
  \end{center}
\endminipage\hfill
\minipage{0.33\textwidth}%
  \includegraphics[width=\linewidth]{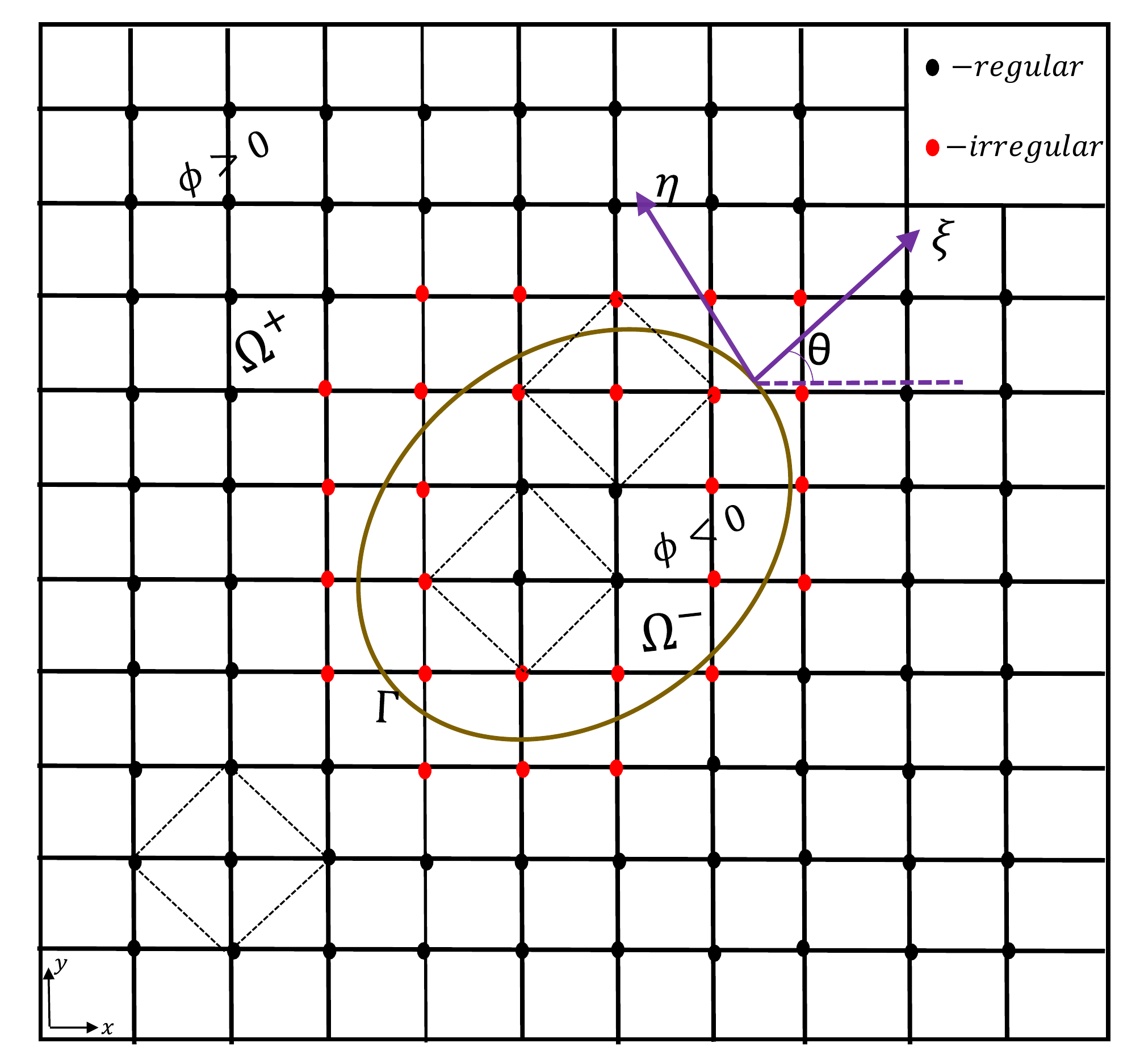}
  \begin{center}
   (c)
  \end{center}
\endminipage
\caption{{\sl (a) Schematic of the problem domain, (b) the unsteady HOC nine point stencil and, (c) regular and irregular points along with the local coordinates on an interfacial point.} }
\label{sch_irr}
\end{figure}

A two dimensional Parabolic interface problem may be modelled as: 
   \begin{equation}
 \lambda u_{t}= \nabla. (\beta \nabla u) + \kappa u -f+b \delta\{(x-x^{*})(y-y^{*})\} \\ 
\quad \textnormal{in} \quad \Omega \times (0, \infty)\;,\;\;(x^{*},y^{*})\in \Gamma \label{s1}
  \end{equation}
with specified  initial and boundary conditions. Here $\Omega$ is an open bounded subset in $\mathbb{R}^{2}$ and $\mathbf{x}=(x,y)$ is an interior point in the domain having an interface  $\Gamma $ immersed in it (see figure \ref{sch_irr}(a)), and  $(x^{*},y^{*}) \in \Gamma$ is an interfacial point. It is assumed that  $\beta(\mathbf{x},t) \in C^{1}(\Omega^{\pm}\backslash \Gamma)$ and $\kappa(\mathbf{x},t)$, $f(\mathbf{x},t) \in C(\Omega^{\pm}\backslash \Gamma)$ may have finite jump across the interface $\Gamma$. Moreover  $b(\mathbf{x},t) \in C(\Gamma)$ and all the parameter $\beta_{x}$ and $\beta_{y}$ are considered to be bounded, hence the solution $u(\mathbf{x},t) \in C^{2}(\Omega^{\pm}\backslash \Gamma) $ . If $\textbf{n}=(n_{1}, n_{2})^{T}$ is the unit outward normal vector to the interface at a point $(x^{*},y^{*}) \in \Gamma$ inside  $\Omega^{-}$ (figure \ref{sch_irr}(a)) and $a(\mathbf{x},t) \in C^{2}(\Gamma)$, the jump conditions in the solution and the flux across the interface may be expressed as 
  \begin{equation}
[u]_{\Gamma}= u^{+} - u^{-}= a(\mathbf{x},t),\label{s2}
\end{equation}  
\begin{equation}
\left[\beta \frac{\partial u}{ \partial \textbf{n}}\right]_{\Gamma}= \beta^{+} \nabla u^{+} - \beta^{-} \nabla u^{-} = b(\mathbf{x},t) \label{s3}
\end{equation}
In other words, the jump conditions in the solution and flux are incorporated in the numerical formulation across the interface to solve such problems, and are defined in (\ref{s2}) and (\ref{s3}) respectively, where subscripts $+$ and $-$ represents the subdomains $\Omega^{+}$ and $\Omega^{-}$ respectively.

In most of the practical cases, the interface has a complex shape. In order to represent the boundary of a bluff body immersed in fluid or the interface, we use the concept of \textbf{level set function} devised by Osher and Sethian in 1996 \cite{osher1988fronts}. They split the domain into sub-domains i.e $\Omega = \Omega^- \cup \Gamma \cup \Omega^+$ and defined the zero level set function $\phi (x,y)$ in two dimensions as a smooth function, 
$$\left\{\begin{array}{ll}
&\phi (x,y) < 0, \;\;\;\;\;\;\; \textnormal{if}\;\; (x,y) \in \Omega^- \\
&\phi (x,y) = 0, \;\;\;\;\;\;\;\textnormal{if}\;\; (x,y) \in \Gamma \\
&\phi (x,y) > 0 \;\;\;\;\;\;\;\;\; \textnormal{if}\;\; (x,y) \in \Omega^+.
\end{array}\right.$$  

In order to discretize equation \eqref{s1} in the finite difference framework in Cartesian grid, we assume the problem to be of rectangular shape given by $\Omega= [x_0,x_f] \times [y_0,y_f]$. The grid is generated by vertical and horizontal lines intersecting across the points $(x_i,y_j)$ given by
$$x_i=x_0+ih, \;\;\;\ y_j=y_0+jl, \;\;\;\; i=0,1,2,...,M-1 \;\;\;, \textnormal{and} \quad j=0,1,2,...,N-1.$$
The step length along $x$- and $y$-directions are defined as $h = \displaystyle{\frac{x_f-x_0}{(M-1)}}$ and $l =\displaystyle{\frac{y_f-y_0}{(N-1)}}$ respectively. The grid points generated this way throughout the whole domain are further sub-categorized into $:$ regular and irregular points. A grid point $\mathbf{x}_{ij}$ is defined as a regular point if all the five points corresponding to a standard central finite difference stencil lies only on one side of the interface $\Gamma$, i.e., either in $\Omega^{-}$ or $\Omega^{+}$. A grid point which is not regular is an irregular point, i.e., $\phi_{ij}^{max} \phi_{ij}^{min} \leq 0$ (see figure \ref{sch_irr}(c)), where 
\begin{equation}
\phi_{ij}^{min}=min \lbrace\phi_{i+1},\phi_{i-1},\phi_{j+1},\phi_{j-1}\rbrace,
\end{equation}  
\begin{equation}
\phi_{ij}^{max}=max \lbrace\phi_{i+1},\phi_{i-1},\phi_{j+1},\phi_{j-1}\rbrace.
\end{equation}  

\subsection{Discretization on the regular points}
We have restructured the HOC finite difference scheme developed by Kalita {\it et al.} \cite{kalita2002class} for the two-dimensional transient convection-diffusion equation to discretize equation \eqref{s1} at the  regular points. Note that the last term on the right hand of equation \eqref{s1} vanishes at regular points and as such \eqref{s1} may be recast in convection-diffusion-reaction form as
 \begin{equation}
 \lambda u_{t} + \beta_xu_x+\beta_yu_y+\beta\nabla^2u + \kappa u =f  
 \label{cd}
  \end{equation}
We use the methodology prescribed in \cite{kalita2002class} to obtain a high order compact finite difference approximation of \eqref{cd} by using a uniform time step $\Delta t$. This is accomplished by first substituting the high order spatial derivatives  appearing in the truncation error terms of the central difference approximation of \eqref{cd} by lower order derivatives obtained from successive differentiation of the original differential equation \eqref{cd}. Next, the time derivative is approximated numerically by a  Crank-Nicolson type of discretization. These result in a spatially fourth and temporally second order accurate finite difference scheme on compact uniform grid requiring nine point stencils both at the $n$ and $(n+1)^{\rm th}$ time levels as shown in figure \ref{sch_irr}(b). As such, it is termed as a $(9,9)$ scheme \cite{kalita2002class}. Eventually, the HOC finite difference approximation of \eqref{cd} can be written as   
\begin{eqnarray}
 \lambda_{ij} \left[ 1 + \frac{h^2}{12} \left(\delta_{xx} +  \frac{(c-2 \beta_x)}{\beta} \delta_{x}  \right) + \frac{l^{2}}{12} \left( \delta_{yy} +  \frac{(d-2 \beta_y)}{\beta} \delta_{y} \right) \right] \left(  u_{ij}^{n+1} - u_{ij}^{n} \right)= \frac{\Delta t}{2}(F_{ij}^{n+1}-F_{ij}^{n}) + \nonumber\\ 
  \frac{\triangle t}{2} \left[ A_{ij} \delta^{2}_{x} + B_{ij} \delta^{2}_{y} + C_{ij} \delta_{x} + D_{ij} \delta_{y} + E_{ij} \delta^{2}_{x}\delta^{2}_{y} + H_{ij} \delta_{x}\delta^{2}_{y} + K_{ij} \delta^{2}_{x}\delta_{y} +  L_{ij} \delta_{x} \delta_{y}  + M_{ij}\right] (u_{ij}^{n+1} + u_{ij}^{n})\nonumber \\
   + O((\triangle t)^{2},(h^{4}, l^{4})) \label{s4}
 \end{eqnarray} 
 where $\delta^{2}_{x}$, $\delta^{2}_{y}$ , $ \delta_{x}$ , $\delta_{y}$, $\delta_{x} \delta_{y}$, $\delta_{x}\delta^{2}_{y}$, $\delta^{2}_{x}\delta_{y}$ and $\delta^{2}_{x}\delta^{2}_{y}$ are second order accurate central difference operators along $x$- and $y$- directions and,
 \begin{equation*}
A_{ij}=\beta_{ij} + \frac{h^{2}}{12 } \left(\beta_{xx}+ 2 c_{x}+\kappa + \frac{(c-2\beta_x)}{\beta} (\beta_x+c)\right)_{ij} +\frac{l^{2}}{12 } \left( \beta_{yy} + \beta_{y} \frac{(d-2\beta_y)}{\beta}\right)_{ij},
 \end{equation*}
  \begin{equation*}
B_{ij}=\beta_{ij} + \frac{h^{2}}{12 } \left( \beta_{xx} + \beta_{x} \frac{(c-2\beta_x)}{\beta}\right)_{ij}  +\frac{l^{2}}{12 }  \left(\beta_{yy}+ 2 d_{y}+\kappa + \frac{(d-2\beta_y)}{\beta} (\beta_y+d)\right)_{ij},
 \end{equation*}
\begin{equation*}
C_{ij}=c_{ij} + \frac{h^{2}}{12 } \left(c_{xx}+ 2\kappa_{x} + \frac{(c-2\beta_x)}{\beta} (c_x+ \kappa)\right)_{ij} +\frac{l^{2}}{12 } \left( c_{yy} + c_{y} \frac{(d-2\beta_y)}{\beta}\right)_{ij},
 \end{equation*}
  \begin{equation*}
D_{ij}=d{ij} + \frac{h^{2}}{12 } \left( d_{xx} + d_{x} \frac{(c-2\beta_x)}{\beta}\right)_{ij}  +\frac{l^{2}}{12 }  \left(d_{yy}+ 2 \kappa_{y} + \frac{(d-2\beta_y)}{\beta} (d_y+ \kappa)\right)_{ij},
 \end{equation*}
\begin{equation*}
E_{ij}=\beta_{ij}\bigg(\frac{h^{2} }{12 } + \frac{l^{2}}{12 }\bigg) \textnormal{, } \quad  H_{ij}=c_{ij} \bigg(\frac{h^{2} }{12 } + \frac{l^{2}}{12 }\bigg) \textnormal{, } \quad  K_{ij}=d_{ij}\bigg(\frac{h^{2} }{12 } + \frac{l^{2}}{12 }\bigg), 
\end{equation*}
\begin{equation}
L_{ij}=\frac{h^{2}}{12}\left(2 d_x + d \frac{(c- 2 \beta_x)}{\beta} )\right)_{ij} + \frac{l^{2}}{12}\left(2 c_y + c \frac{(d- 2 \beta_y)}{\beta} )\right)_{ij} ,\nonumber 
\end{equation}
\begin{equation}
M_{ij}=\kappa_{ij} + \frac{h^{2}}{12 } \left( \kappa_{xx} + \kappa_{x} \frac{(c-2 \beta_x)}{\beta} \right)_{ij} +\frac{l^{2}}{12} (\beta_{ij}\left( \kappa_{yy} + \kappa_{y} \frac{(d-2 \beta_y)}{\beta} \right)_{ij},\nonumber 
\end{equation}
\begin{equation}
F_{ij}=f_{ij} + \frac{h^{2}}{12 } \left(f_{xx} + f_{x} \frac{(c-2 \beta_x)}{\beta}  \right)_{ij}+ \frac{l^{2}}{12} \left( f_{yy} + f_y  \frac{(d-2 \beta_y)}{\beta} \right)_{ij}.\nonumber
\end{equation}
Rewrite the equation (\ref{s4}) into 
\begin{equation}
\sum_{i=1}^9 ~c_{i}u_{i}^{n+1}=\sum_{i=1}^9 ~c_{i}^{'}u_{i}^{n}+ \frac{\Delta t}{2}\left(F_{ij}^{n+1}-F_{ij}^{n} \right), \label{reg}
\end{equation}
where  ${\displaystyle c_{i}= \lambda m_{i}-\frac{\Delta t}{2} n_{i}}$, $\quad$ ${\displaystyle c^{'}_{i}
=\lambda m_{i}-\frac{\Delta t}{2} n_{i}}$, \quad and \quad $m_{1}=m_{3}=m_{7}=m_{9}=0$, $m_{5}={\displaystyle\frac{8}{12}}$,\\ ${\displaystyle m_{2}=\frac{1}{12}\left(1- l\frac{(d-2 \beta_y)}{\beta} \right)}$,
\begin{eqnarray}
m_{8}=\frac{1}{12}\left(1+ l\frac{(d-2 \beta_y)}{\beta} \right), \quad m_{4}=\frac{1}{12}\left(1- h\frac{(c-2 \beta_x)}{\beta} \right), \quad m_{6}=\frac{1}{12}\left(1+ h\frac{(c-2 \beta_x)}{\beta} \right),
\end{eqnarray}
 
\begin{equation}
n_{1}=\frac{E_{ij}}{h^{2}l^{2}}-\frac{H_{ij}}{2hl^{2}}-\frac{K_{ij}}{2h^{2}l} + \frac{L_{ij}}{4hl} \textnormal{, } \quad  n_{2}=\frac{B_{ij}}{l^{2}}-\frac{D_{ij}}{2l}-\frac{2E_{ij}}{h^{2}l^{2}}+\frac{K_{ij}}{h^{2}l} \textnormal{, } \quad n_{3}=\frac{E_{ij}}{h^{2}l^{2}}+\frac{H_{ij}}{2hl^{2}}-\frac{K_{ij}}{2h^{2}l} - \frac{L_{ij}}{4hl} \textnormal{, }
\end{equation}
\begin{equation}
n_{4}=\frac{A_{ij}}{h^{2}}-\frac{C_{ij}}{2h}-\frac{2E_{ij}}{h^{2}l^{2}}+\frac{H_{ij}}{hl^{2}} \textnormal{, }\quad  n_{5}=-\frac{2A_{ij}}{h^{2}}-\frac{2B_{ij}}{l^{2}}+\frac{4E_{ij}}{h^{2}l^{2}}+M_{ij} \textnormal{, } \quad  n_{6}=\frac{A_{ij}}{h^{2}}+\frac{C_{ij}}{2h}-\frac{2E_{ij}}{h^{2}l^{2}}-\frac{H_{ij}}{ hl^{2}} \textnormal{, }
\end{equation}
\begin{equation}
n_{7}=\frac{E_{ij}}{h^{2}l^{2}}+\frac{H_{ij}}{2hl^{2}}+\frac{K_{ij}}{2h^{2}l} + \frac{L_{ij}}{4hl} \textnormal{, }\quad  n_{8}=\frac{B_{ij}}{l^{2}}+\frac{D_{ij}}{2l}+\frac{2E_{ij}}{h^{2}l^{2}}-\frac{K_{ij}}{h^{2}l} \textnormal{, } \quad n_{9}=\frac{E_{ij}}{h^{2}l^{2}}-\frac{H_{ij}}{2hl^{2}}+\frac{K_{ij}}{2h^{2}l} - \frac{L_{ij}}{4hl}  \textnormal{. }
\end{equation}

It is worth mentioning that the scheme developed in \cite{kalita2002class} was devoid of the reaction term present in \eqref{cd} and of any variable diffusion coefficients.

\subsection{Discretization on the irregular points}
The schematic of the irregular points across the interface along with the regular ones in the computational plane can be seen in figure \ref{sch_irr}(c). Let ($\eta, \xi$) represent the local coordinate system at an interfacial point $(x^\star,y^\star)$ with $\eta$ and $\xi$ representing the tangent and normal direction respectively at the point along the interface. Then for approximating the jump conditions \eqref{s3} on a Cartesian mesh at the point $(x^{\star}, y^{\star})$, we have
\begin{eqnarray*}
\xi=&(x-x^{\star})  cos(\theta )+ (y-y^{\star}) sin(\theta),\\
\eta=&-(x-x^{\star})  sin(\theta )+ (y-y^{\star}) cos(\theta).
\end{eqnarray*}
where $\theta$ is the angle between $x$-axis and $\xi$-direction. The jump conditions for the derivatives up to third order can be calculated by the following formulas.
\begin{equation}
[u_{x}]=cos(\theta)[u_{\xi}]-sin(\theta)[ u_{\eta}], \label{aj1}
\end{equation}
\begin{equation}
[u_{y}]=sin(\theta)[u_{\xi}]+cos(\theta)[ u_{\eta}] \label{aj2}
\end{equation}
\begin{equation}
\lbrack u_{xx}\rbrack=cos^{2}(\theta)[u_{\xi \xi}]-2cos(\theta)sin(\theta)[u_{\xi \eta}] + sin^{2}(\theta)[u_{\eta \eta}] \label{aj3}
\end{equation}
  \begin{equation}
\lbrack u_{yy}\rbrack=sin^{2}(\theta)[u_{\xi \xi}]+2cos(\theta)sin(\theta)[u_{\xi \eta}] + cos^{2}(\theta) [u_{\eta \eta}]    \label{aj4}
  \end{equation}
  \begin{equation}
\lbrack u_{xxx}\rbrack=cos^{3}(\theta) [u_{\xi \xi \xi}] - 3 cos^{2}(\theta)sin(\theta) [u_{\xi \xi \eta}]  + 3 cos(\theta)sin^{2}(\theta) [u_{\xi \eta \eta}] -sin^{3}(\theta)[u_{\eta \eta \eta}]   \label{aj5}
  \end{equation}
  \begin{equation}
\lbrack u_{yyy}\rbrack=sin^{3}(\theta) [u_{\xi \xi \xi}] + 3 cos(\theta)sin^{2}(\theta) [u_{\xi \xi \eta}]  + 3 cos(\theta)sin^{2}(\theta) [u_{\xi \eta \eta}] +cos^{3}(\theta)[u_{\eta \eta \eta}]   \label{aj6}
  \end{equation}

In the subsequent sections, we explore several feasible scenarios for the irregular points and the discretization of the equation thereat.

\subsubsection{Irregular points lying on grid lines parallel to $x$-axis only}\label{sectx}
\begin{figure}[!h]
\begin{minipage}[b]{.5\linewidth}  
\includegraphics[scale=1.2]{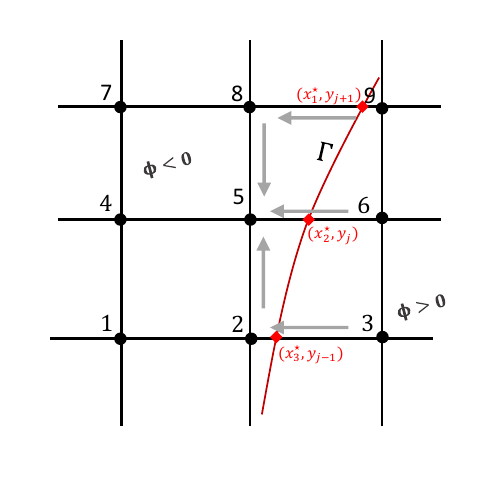} 
 \centering (a)
\end{minipage}            \hspace{-2.mm}
\begin{minipage}[b]{.5\linewidth}
\includegraphics[scale=1.2]{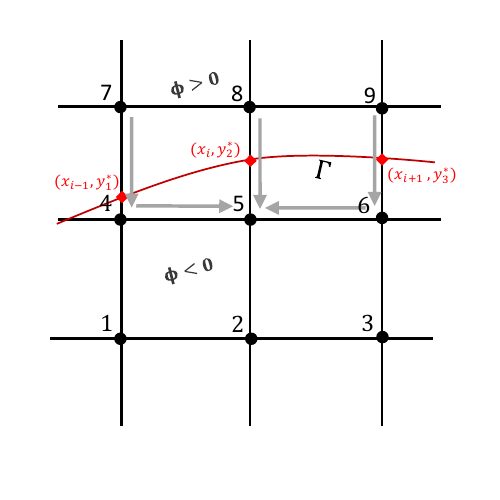} 
 \centering (b)
\end{minipage} 
\caption{{\sl Stencils around the irregular points lying on grid lines parallel to (a) $x$-axis only and (b) $y$-axis only.} }
\label{x_irr}
\end{figure}

Here, we describe the case when the irregular point lies only on grid lines parallel to $x-$axis. Let us assume that such a grid line meets the interface between $(x_{i},y_{j})$ and $(x_{i+1},y_{j})$ (i.e $\phi_{i+1,j} \times \phi_{i,j} <0$) at the point $(x_{2}^{\star},y_{j})$, with the possibility of the interface cutting the grid lines above and below  $y_j$ level at the points $(x_{1}^{\star},y_{j+1})$ and $(x_{3}^{\star},y_{j-1})$ respectively on a nine point compact stencil as shown in figure \ref{x_irr}(a). From this figure, it is clear that  $\phi>0$ on these three points $(x_{i+1},y_{j+1})$, $(x_{i+1},y_{j})$ and $(x_{i+1},y_{j-1})$, while the remaining six points lie on the other side of the interface where $\phi<0$. 

For approximating the mixed derivatives appearing in \eqref{s4}, we apply Taylor series expansion to approximate $ u(x_{i+1},y_{j+1})$, $ u(x_{i+1},y_{j})$ and $ u(x_{i+1},y_{j-1})$ about the grid point $(x_{i},y_{j})$ by including the jumps in the solution and the derivatives along the $x$-direction at the interfacial points. This can be accomplished by firstly expanding the Taylor series in the direction of irregularity i.e. along $x$- axis as indicated by the arrowheads in figure \ref{x_irr}(a) and then moving in the other direction along $y$-axis. The following lemma ensures the high order accuracy of these approximations, detailed proof of which can be found in the authors' paper \cite{singhalkalita} along with other possible scenarios.

\begin{lemma} \label{l1}
 Let $u^{-} \in C^{k+1} [x_0,x_{1}^{\star}] \times [y_0,y_f]$, $u^{+} \in C^{k+1} [x_{1}^{\star},x_f] \times [y_0,y_f]$, $h=x_{i+1}-x_{i}$, $h_{1}^{+}=x_{i+1}-x_{1}^{\star}$\\ and $h_{1}^{-}$=$x_{i}-x_{1}^{\star}$ then we have the following inequality
\begin{eqnarray} \label{p1}
\bigg\Vert u(x_{i+1},y_{j+1}) -\sum_{p=0}^{k}\sum_{q=0}^{k-p} \frac{h^{p} l^{q}}{p! q!} \frac{\partial^{p+q} u}{\partial x^{p} \partial y^{q}}(x_{i},y_{j}) - \sum_{r=0}^{k} \frac{(h_{1}^{+})^{r} }{r!} \left[\frac{\partial^r u}{\partial x^{r}}(x_{1}^{\star},y_{j+1})\right]  \bigg \Vert \leq \nonumber
 K \frac{h^{k+1}}{(k+1)!}+\frac{M}{(k+1)!}(|h|+|l|)^{k+1}
\end{eqnarray}
where $K$=$\max(\max_{x \in [x_{i},x_{1}^{\star}) } \mid u^{k+1}(x_{1}^{\star},y_{j+1})\mid$ , $\max_{x \in (x_{1}^{\star},x_{i+1}] } \mid u^{k+1}(x_{1}^{\star},y_{j+1})\mid)$
\end{lemma}
\begin{remark}\label{r2}
 Let  $h_{3}^{+}=x_{i+1}-x_{3}^{\star}$ and $h_{3}^{-}$=$x_{i}-x_{3}^{\star}$ then we have the following inequality
\begin{eqnarray}
\bigg\Vert u(x_{i+1},y_{j-1}) -\sum_{p=0}^{k}\sum_{q=0}^{k-p} \frac{h^{p} (-l)^{q}}{p! q!} \frac{\partial^{p+q} u}{\partial x^{p} \partial y^{q}}(x_{i},y_{j}) - \sum_{r=0}^{k} \frac{(h_{3}^{+})^{r} }{r!} \left[\frac{\partial^r u}{\partial x^{r}}(x_{3}^{\star},y_{j-1})\right]  \bigg \Vert \leq O(h^{k+1},l^{k+1}).
\end{eqnarray}
\end{remark}  
For $k=3$, the above Lemma and Remarks guarantee that the approximation of \eqref{cd} is fourth order accurate in spatial direction (see table \ref{ex1_table1} for test case $1$), yielding 
\begin{equation}
\delta_{xy}u (x_{i},y_{j}) =\frac{\partial^{2} u}{\partial x \partial y} -\frac{1}{4hl}\left(\sum_{r=0}^{k}\frac{(h_{1}^{+})^{r}}{r!} \left[\frac{\partial^r u}{\partial x^{r}}(x_{1}^{\star},y_{j+1})\right]-\sum_{r=0}^{k}\frac{(h_{3}^{+})^{r} }{r!}  \left[\frac{\partial^r u}{\partial x^{r}}(x_{3}^{\star},y_{j-1})\right] \right)+O(h^{2},l^{2}),
\end{equation}
\begin{equation}
\delta_{xx} \delta_{y}u_{ij} =\frac{\partial^{3} u}{\partial x^{2} \partial y}+\frac{1}{2h^{2}l}\left(\sum_{r=0}^{k}\frac{(h_{3}^{+})^{r} }{r!}  \left[\frac{\partial^r u}{\partial x^{r}}(x_{3}^{\star},y_{j-1})\right] -\sum_{r=0}^{k}\frac{(h_{1}^{+})^{r}}{r!} \left[\frac{\partial^r u}{\partial x^{r}}(x_{1}^{\star},y_{j+1})\right]\right)+O(h^{2},l^{2}), 
\end{equation}
\begin{eqnarray}\label{a17}
\delta_{x} \delta_{yy}u_{ij} =\frac{\partial^{3} u}{\partial x \partial y^{2}} +\frac{1}{2hl^{2}}\left( 2 \sum_{r=0}^{k}\frac{(h_{2}^{+})^{r}  }{r!} \left[\frac{\partial^r u}{\partial x^{r}}(x_{2}^{\star},y_{j})\right]- \sum_{r=0}^{k}\frac{(h_{1}^{+})^{r}}{r!} \left[\frac{\partial^r u}{\partial x^{r}}(x_{1}^{\star},y_{j+1})\right]-\sum_{r=0}^{k}\frac{(h_{3}^{+})^{r} }{r!}  \left[\frac{\partial^r u}{\partial x^{r}}(x_{3}^{\star},y_{j-1})\right] \right) \nonumber \\ +O(h^{2},l^{2}),\\
\delta_{xx}\delta_{yy}u_{ij} =\frac{\partial^{4} u}{\partial x^{2} \partial y^{2}}+\frac{1}{h^{2}l^{2}}\left(2  \sum_{r=0}^{k}\frac{(h_{2}^{+})^{r}} {r!} \left[\frac{\partial^r u}{\partial x^{r}}(x_{2}^{\star},y_{j})\right]-\sum_{r=0}^{k}\frac{(h_{1}^{+})^{r}}{r!} \left[\frac{\partial^r u}{\partial x^{r}}(x_{1}^{\star},y_{j+1})\right]-\sum_{r=0}^{k}\frac{(h_{3}^{+})^{r} }{r!}  \left[\frac{\partial^r u}{\partial x^{r}}(x_{3}^{\star},y_{j-1})\right] \right)\nonumber \\ +O(h^{2},l^{2}).
\end{eqnarray}
With these, equation \eqref{reg} at the irregular point $(x_i,y_j)$ reduces to
\begin{eqnarray}
\sum_{i=1}^9 ~c_{i}u_{i}^{n+1} -c_{3} \sum_{r=0}^{k}\frac{(h_{1}^{+})^{r}}{r!} \left[\frac{\partial^r u}{\partial x^{r}}(x_{1}^{\star},y_{j+1})\right]^{(n+1)}- c_{6} \sum_{r=0}^{k}\frac{(h_{2}^{+})^{r}} {r!} \left[\frac{\partial^r u}{\partial x^{r}}(x_{2}^{\star},y_{j})\right]^{(n+1)}+ c_{9} \sum_{r=0}^{k}\frac{(h_{3}^{+})^{r} }{r!}  \left[\frac{\partial^r u}{\partial x^{r}}(x_{3}^{\star},y_{j-1})\right]^{(n+1)} \nonumber \\=
\sum_{i=1}^9 ~c_{i}^{'}u_{i}^{n}-c^{'}_{3} \sum_{r=0}^{k}\frac{(h_{1}^{+})^{r}}{r!} \left[\frac{\partial^r u}{\partial x^{r}}(x_{1}^{\star},y_{j+1})\right]^{(n)}-c^{'}_{6}  \sum_{r=0}^{k}\frac{(h_{2}^{+})^{r}} {r!} \left[\frac{\partial^r u}{\partial x^{r}}(x_{2}^{\star},y_{j})\right]^{(n)}+ c^{'}_{9}\sum_{r=0}^{k}\frac{(h_{3}^{+})^{r} }{r!}  \left[\frac{\partial^r u}{\partial x^{r}}(x_{3}^{\star},y_{j-1})\right]^{(n)} \nonumber \\
+ \frac{\Delta t}{2}\left(F_{ij}^{n+1}-F_{ij}^{n} \right). \label{irr_reg}
\end{eqnarray}

\subsubsection{Irregular points lying on grid lines parallel to $y$-axis only}\label{secty}
The treatment for irregular points lying only on grid lines parallel to $y-$axis is similar to the cases described in the above section. The details of the spatial discretizations of interfacial points lying between $(x_{i},y_{j})$ and $(x_{i},y_{j+1})$ as shown in figure \ref{x_irr}(b) and other possible cases including treatment of irregular points lying simultaneously on grid lines parallel to both x-axis and y-axis can be found in \cite{singhalkalita}.

\section{Streamfunction-Vorticity Formulation}
The Navier-Stokes equations provide the mathematical framework for incompressible viscous flows and as such, are the backbones in both theoretical and computational fluid dynamics studies. In particular, the numerical solutions of the N-S equations has played an important role in the recent advances in the flow past bluff bodies in the field of aerospace engineering. The non-dimensional form of the N-S equations in the primitive variable formulation in two dimensions can be written as:
\begin{equation}
\nabla . \mathbf{u} = 0 \label{cont}
\end{equation}
\begin{equation}
\frac{\partial \mathbf{u}}{\partial t} + \mathbf{u} . \nabla \mathbf{u} = - \nabla p + \frac{1}{Re} \triangle \mathbf{u}. \label{NS}
\end{equation} 
where $\mathbf{u}=(u, v)$ is the velocity field, $t$ the time, $p$ the pressure and $\displaystyle Re=\frac{UL}{\nu}$ is the Reynolds Number with $U$ and $L$  being some characteristic velocity and length, and $\nu$ the kinematic viscosity of the fluid.

The main interest of this study is to simulate flow fields for the different fixed and moving interfaces in two-dimensional laminar flows. Due to the presence of the pressure term, the direct solution of \eqref{cont}-\eqref{NS} has been a costly affair despite representing the fluid phenomena accurately. In order to overcome this, the streamfunction $\psi$ is introduced as:
\begin{equation}
u=\frac{\partial \psi}{\partial y} \quad \textnormal{and}  \quad v=- \frac{\partial \psi}{\partial x}
\label{vel}
\end{equation}
which allows the preservation of the incompressibility condition of the continuity equation \eqref{cont}. Taking curl of \eqref{NS} reduces it into the vorticity transport equation
\begin{equation}
\frac{\partial\zeta}{\partial t}=\frac{1}{Re}\nabla^{2} \zeta- \mathbf{u}. \nabla \zeta \label{omega}
\end{equation} 
where in the 2D flow field the vorticity vector $\mathbf{\omega}$ is given by $\displaystyle \nabla \times \mathbf{u} =\mathbf{\omega}=\zeta\hat{k}=\left( \frac{\partial v}{\partial x}-  \frac{\partial u}{\partial y}\right)\hat{k}$, $\hat{k}$ being the unit vector normal to the $xy$-plane. From the definition of streamfunction $\psi$ and the scalar vorticity $\zeta$ above, one can get the following Poisson equation for the streamfunction
\begin{equation}
\nabla^{2} \psi=-\zeta. \label{psi}
\end{equation}

Let the surface of the bluff body be represented by $S_b$ with curvilinear coordinates $s$ along it and $\hat{\mathbf{n}}$, $\hat{\mathbf{\tau}}$ be the outward unit normal and tangent vectors respectively. If $\mathbf{u}_S(s,t)$ is the velocity of the fluid on the surface $S_b$, the corresponding boundary conditions for $\psi$ is
\begin{equation}
\left .\frac{\partial \psi}{\partial n}\right |_{S_b}=\mathbf{\hat{n}} \cdot \mathbf{u}_S(s,t), \quad \left .\frac{\partial \psi}{\partial \tau}\right |_{S_b}=-\mathbf{\hat{\tau}} \cdot \mathbf{u}_S(s,t).
\label{bc_psi}
\end{equation}

The system of equations \eqref{vel}-\eqref{bc_psi} completely describes what is known as the Streamfunction-Vorticity ($\psi$-$\zeta$) formulation of the N-S equations. Note that, being parabolic in nature, equation \eqref{omega} can be discretized using \eqref{reg} and \eqref{irr_reg}, while \eqref{psi} being elliptic, one may use the approach developed by the authors in \cite{singhalkalita}. Over the past few decades, this formulation has been employed for 2D fluid flow computations with great success because of its ease of implementation. In particular, large number of researchers have utilized it in order to check the efficiency of newly developed methods by numerical solving a variety of challenging fluid flow problems. However, owing to the non-specification of vorticity values at the no-slip boundaries, one needs to devise specific approach to approximate the vorticity values thereat. This task becomes more trickier in the immersed interface framework over Cartesian grids for curved boundaries, which is detailed in the next section.

\subsection{Treatment of irregular points for $\psi$ and $\zeta$}
\begin{figure}[!h]
\begin{center} 
\includegraphics[scale=.6]{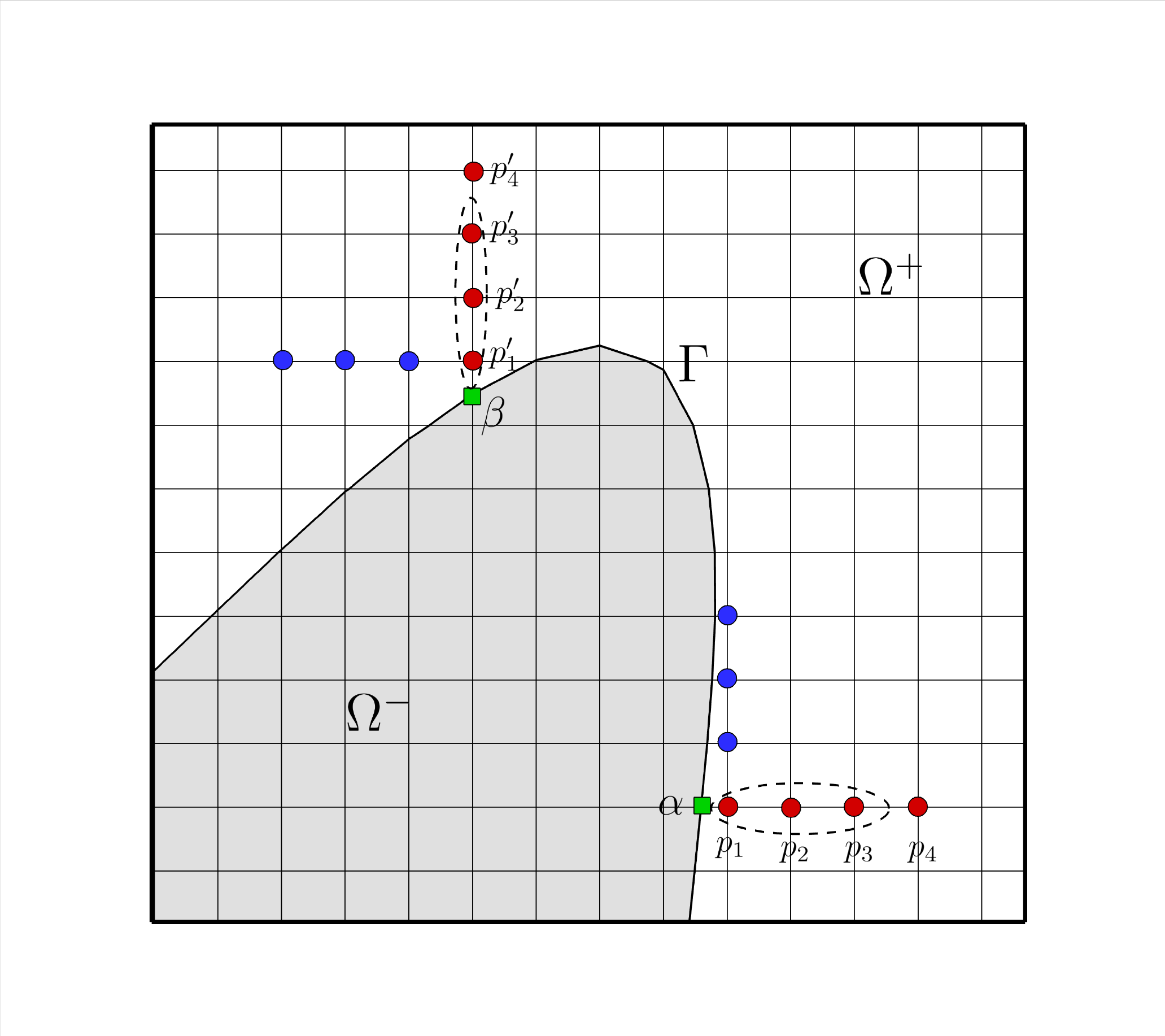}  
\end{center}
\caption{{\sl Schematic of the neighbourhood of the interface for the evaluation of correction terms.} }
\label{sch_jmp}
\end{figure} 
A quick look at equation \eqref{irr_reg} reveals that evaluation of the solution at the $(n+1)^{\rm th}$ level requires the correction terms on the interfacial points across the interface both at $(n+1)^{\rm th}$ and $(n)^{\rm th}$ levels. However, such evaluation for the streamfunction and vorticity across the irregular points is not that straightforward. Therefore, a specific interpolation strategy is adopted to calculate the jumps by mapping the values of $u$ and its derivatives at the regular and irregular points across the interface.  Consider two interfacial points $\alpha$ and $\beta$ (denoted by the green solid squares) as shown in figure \ref{sch_jmp}, which corresponds to the classification of irregular points represented in the section \ref{sectx} and \ref{secty} respectively. If  $u^{(n)}(\alpha)$ denotes the $n^{\rm th}$ order partial derivative of $u$ at the point $\alpha$ lying along $x$-axis, the jump invariably involves $u$ and $u^{(n)}(\alpha)$ thereat from either sides of the interface. However, for the test cases involving fluid flows under consideration here, the physical constraints allow the solution to be zero inside the immersed boundary. 

Note that for such cases, all these derivatives must be calculated via evaluation of the variable $u(\mathbf{x})$ through the one sided finite difference formula at an irregular point. For example, in order to evaluate $u$ at the interfacial point $\alpha$, we employ a Lagrangian interpolation polynomial by making use of the first irregular node $p_{1}$ on its right hand side (One can choose either left or right side depending upon the location of the interfacial point.) and the subsequent regular nodes $p_{2}$ and $p_{3}$  as shown in figure \ref{sch_jmp}. A unique polynomial of degree two is
\begin{equation}
P(x)= \sum_{j=1}^3 u(p_{j}) \textit{l}_{j} (x) + \frac{u^{(4)} (\xi(\mathbf{x})) }{ 4!} \prod_{i=1}^3 (\alpha - p_{i}), 
\label{lg}
\end{equation}
where $\xi(\mathbf{x})$ is some number lying in the interval $\displaystyle \left( min\lbrace p_{i} \rbrace, max\lbrace p_{i} \rbrace \right)_{1 \leq i \leq 3}$ and 
\begin{equation}
 \textit{l}_{j} (x)= \prod_{i \neq j} \frac{(x - p_{i})}{(p_{j}-p_{i})}.
 \label{lg_b}
\end{equation}
In the correction of jumps, the successive one-directional derivatives of the variables are also of utmost importance, which can be approximated by the successive differentiation of \eqref{lg}. As can be seen from \eqref{lg}, it involves the differentiation of the basis  functions $\textit{l}_{j}$ at each node $p_j$. This can be accomplished by taking its logarithm
\begin{equation}
ln \left(\textit{l}_{j} (x) \right)= ln \left(\prod_{i \neq j} \frac{(x - p_{i})}{(p_{j}-p_{i})} \right)= \sum_{i \neq j} ln \left( \frac{(x - p_{i})}{(p_{j}-p_{i})} \right). \nonumber
\end{equation}
Differentiating the above, we have:
\begin{equation}
\frac{\textit{l}^{\prime}_{j} (x)}{\textit{l}_{j} (x)}= \sum_{i \neq j} \frac{1/(p_{j}-p_{i})}{(x - p_{i})/(p_{j}-p_{i})}=\sum_{i \neq j} \frac{1}{(x - p_{i})}, \nonumber
\end{equation}
which yields 
\begin{equation}
\textit{l}^{\prime}_{j} (x)=\textit{l}_{j} (x) \left( \sum_{i \neq j} \frac{1}{(x - p_{i})} \right).
\label{lg_d1}
\end{equation}
Applying product rule for the derivatives in \eqref{lg_d1}
\begin{eqnarray}
\textit{l}^{\prime \prime}_{j} (x)&=& \textit{l}^ \prime_{j} (x) \left( \sum_{i \neq j} \frac{1}{(x - p_{i})} \right) + \textit{l}_{j} (x) \left( \sum_{i \neq j} \frac{1}{(x - p_{i})} \right)^{\prime}, \nonumber \\
&=& \textit{l}^ \prime_{j} (x) \left( \sum_{i \neq j} \frac{1}{(x - p_{i})} \right) + \textit{l}_{j} (x) \left( \sum_{i \neq j} \frac{-1}{(x - p_{i})^2} \right). 
\label{lg_d2}
\end{eqnarray}
Substituting \eqref{lg_d1}, \eqref{lg_d2} reduces to
\begin{equation}
\textit{l}^{\prime \prime}_{j} (x)=\textit{l}_{j} (x) \left[ \left( \sum_{i \neq j} \frac{1}{(x - p_{i})} \right)^{2} - \left( \sum_{i \neq j} \frac{1}{(x - p_{i})^{2}} \right) \right]. 
\end{equation}

On the surface of the bluff bodies, which accounts for the interface, the jump condition for $\psi$ can be computed using 
\begin{equation}
[\psi]= 0,\quad \left[ \frac{\partial \psi}{\partial x} \right]=-\left[v \right],\quad \left[ \frac{\partial \psi}{\partial y} \right]=\left[u \right], \label{jpsi1}
\end{equation}
\begin{equation}
\left[ \frac{\partial^2 \psi}{\partial x^2} \right]=-\left[ \frac{\partial v}{\partial x} \right] \rm{and} \label{jpsi2}
\end{equation}
\begin{equation}
\left[ \frac{\partial^2 \psi}{\partial y^2} \right]=\left[ \frac{\partial u}{\partial y} \right].\label{jpsi3}
\end{equation}

In equation \eqref{jpsi2} which is utilized for interfacial points lying only on $x$-axis, $\displaystyle \frac{\partial v}{\partial x}$ is computed using the one-sided second order approximation
\begin{equation}
\left[ \frac{\partial v}{\partial x} \right]=  \left . \frac{\partial v}{\partial x} \right|_{\alpha^+}=\frac{1}{h \delta h (\delta h+h)}\left( -h(2\delta h+h)v(\alpha)+(\delta h+h)^2v(p_1)-(\delta h)^2v(p_2)\right).\label{jpsi4}
\end{equation}
where $\delta h=|p_1-\alpha|$ and $v(\alpha)$ is computed using equation \eqref{lg}. Likewise  in \eqref{jpsi3}, which is typical of irregular points lying on $y$-axis, the roles of $\alpha$, $p_1$ and $p_2$ in equation \eqref{jpsi2} are carried out by $\beta$, $p_1'$ and $p_2'$ respectively for computing $\displaystyle \frac{\partial u}{\partial y}$. While exact jump conditions  are not difficult to find for streamfunction because of the availability of exact boundary conditions on solid surfaces, no such conditions are available for vorticity. 

For vorticity, at the point $\alpha$, $\displaystyle \left [\zeta\right]=\zeta(\alpha^+)$ is evaluated at the current time level $(n+1)$ by making use of equation \eqref{lg}, which requires the value of $\displaystyle \zeta(p_1)^{(n+1)}$. However, $p_1$ being an irregular point,  $\displaystyle \zeta(p_1)^{(n+1)}$ is not readily available thereat. In order to circumvent this, a one-sided $O(h^3)$ approximation is utilized to compute $\zeta$ by a one-sided discretization of $-\nabla^2 \psi (p_1)$ in $x$ and $y$-directions, viz., making use of the nodes next right and above $p_1$, denoted by blue and red dots respectively as shown in figure \ref{sch_jmp}. We use
\begin{equation}
\frac{\partial^2 \psi}{\partial x^2}(p_1)= \frac{1}{h^2} \left(2 \psi(p_1) - 5 \psi (p_2) + 4 \psi(p_3) - \psi (p_4) \right) + O(h^3)\label{one_sided_xx}
\end{equation} 
and likewise for $\displaystyle \frac{\partial^2 \psi}{\partial y^2}(p_1')$. Again  $\displaystyle \left [\frac{\partial \zeta}{\partial x}\right]=\left .\frac{\partial \zeta}{\partial x}\right|_{\alpha^+}$ at $\alpha$ and $\displaystyle \left [\frac{\partial \zeta}{\partial y}\right]=\left .\frac{\partial \zeta}{\partial y}\right|_{\beta^+}$ at $\beta$, the procedure for approximating which is similar to the ones for finding out the jumps of first order derivatives described in equation \eqref{jpsi4}. The approximations for the second order jump condition is 

\begin{equation}
\displaystyle \left [\frac{\partial^2 \zeta}{\partial x^2}\right]=\left .\frac{\partial^2 \zeta}{\partial x^2}\right|_{\alpha^+}=\frac{2}{h \delta h(\delta h+h)}\left( h v(\alpha)-(\delta h+h)v(p_1)+\delta h v(p_2)\right).
\end{equation}
Likewise, jump conditions for higher order derivatives can also be estimated.
\subsection{Fluid dynamic forces on the body}
\begin{figure}[!h]
 \begin{minipage}[b]{.45\linewidth}  
\includegraphics[scale=0.48]{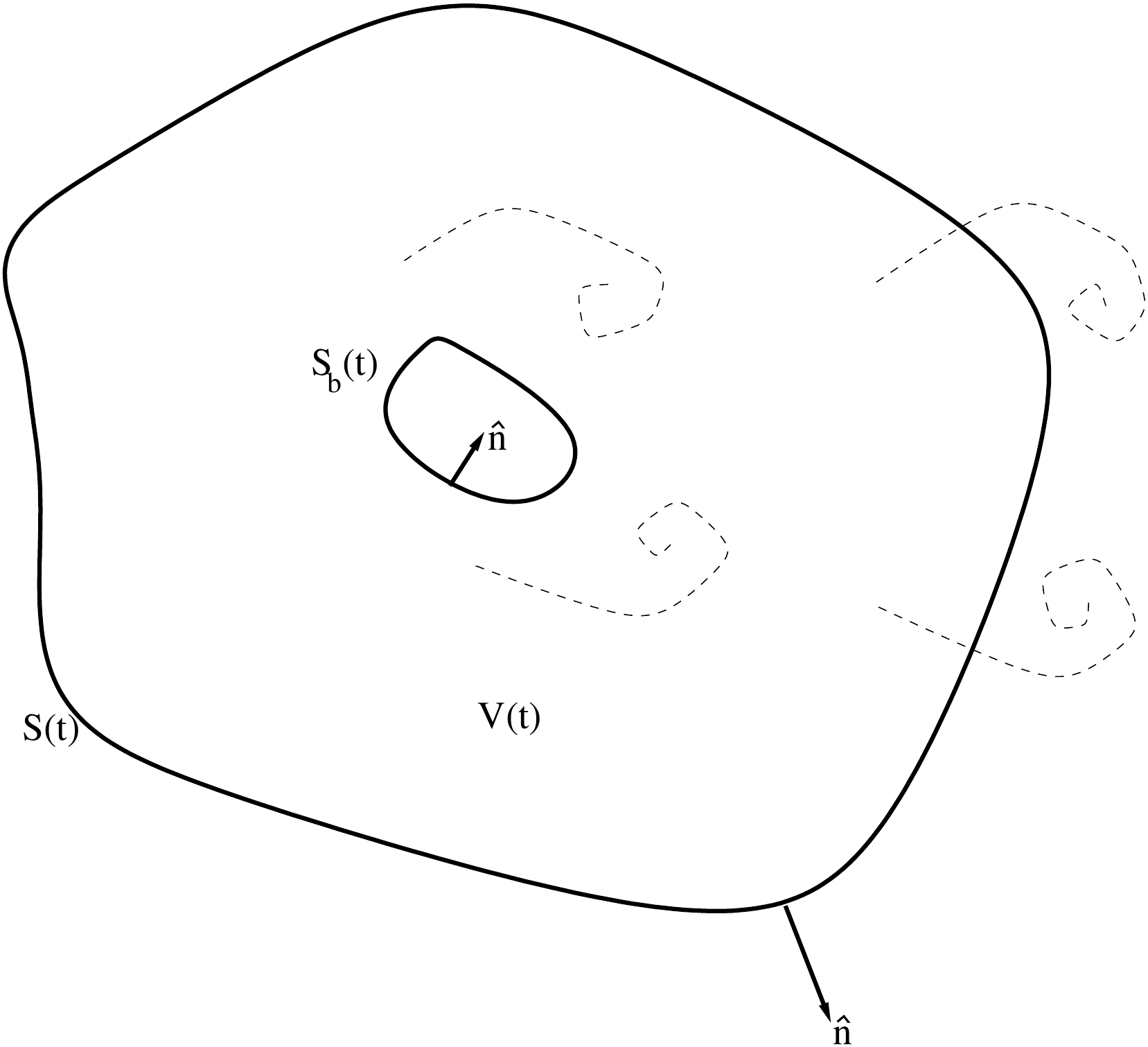} 
\centering (a) 
\end{minipage}           
\begin{minipage}[b]{.45\linewidth}
\includegraphics[scale=0.37]{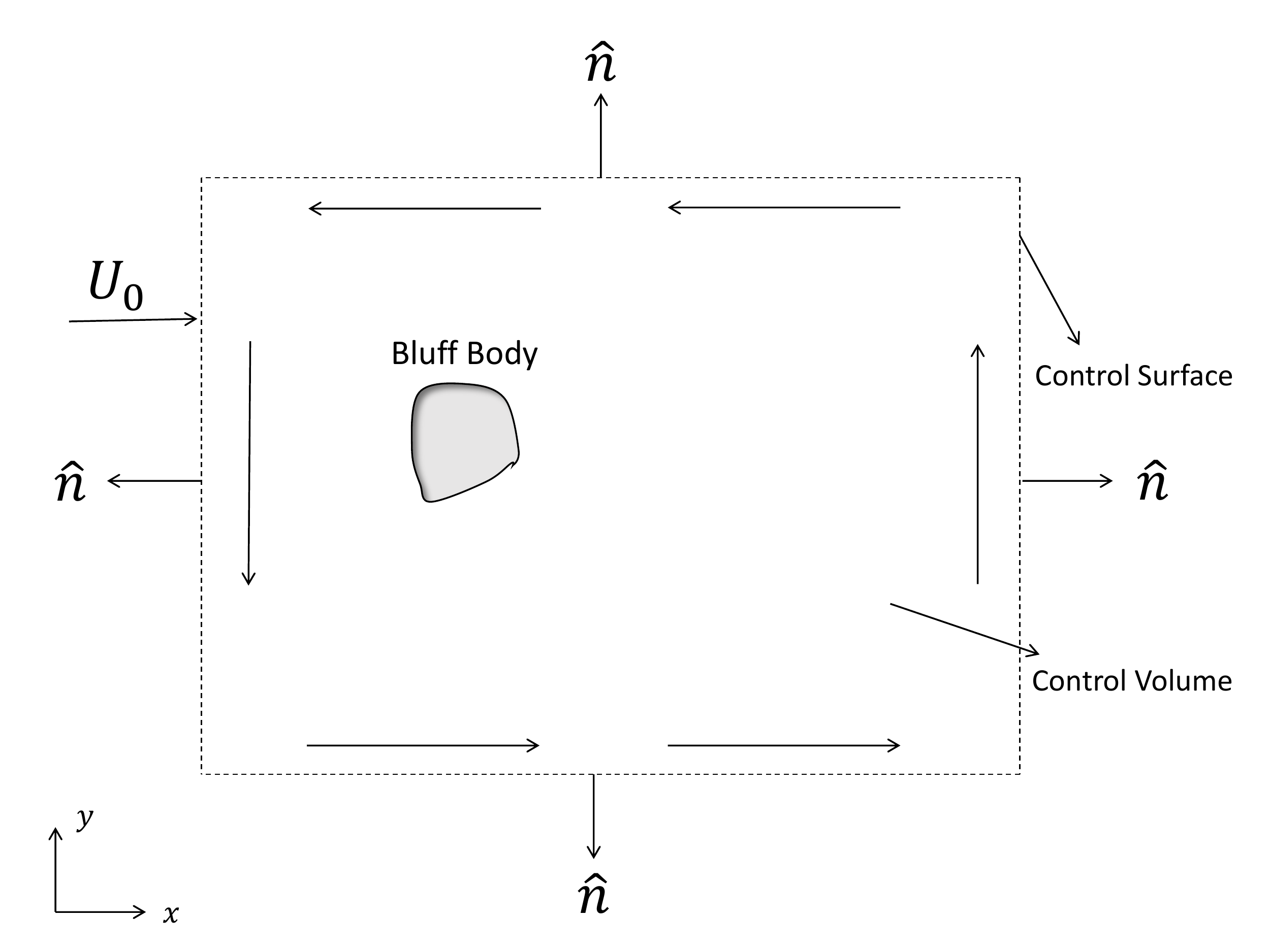}
\centering (b) 
\end{minipage} 
\caption{{\sl  Schematic of (a) the domain of integration for body force evaluation on a bluff body and (b) rectangular region representing the Control volume.} }
\label{sch_dl}
\end{figure}

When a body is immersed into a fluid in relative motion, the fluid exerts a force on the bluff body which can be derived from the equations of motion \eqref{cont}-\eqref{NS}. We have utilized the momentum approach adopted by Noca {\it et al.} \cite{noca1999comparison}, who devised a formula that does not require explicit knowledge of the pressure term. Note that Equation \eqref{NS} is nothing but confirmation of Newton's second law, which states that the time rate of change within the control volume  is equal to the net force. The momentum balance is written in integral form by considering an arbitrary time-dependent control volume $V(t)$ bounded externally by a control surface $S(t)$ and internally by the body surface $S_b(t)$ as shown in figure \ref{sch_dl}(a). Thus $V(t)$ is a simply connected region. The fluid dynamic force $\vec{F}$ in dimensional form acting on body enclosed by a fixed control volume can be written as
\begin{equation}
\vec{F}=-\frac{d}{dt} \int_{V(t)} \rho\vec{u} dV + \oint_{S(t)} \hat{n}. \gamma_{mom} \,ds - \oint_{S_b(t)} \rho\hat{n}. (\vec{u} - \vec{u_{s}}) \,ds \label{m1}
\end{equation} 
where $\rho$ is the density of the fluid, $\hat{n}$ is a unit normal vector, $\vec{u}$ is the flow velocity, $u_s$ is the velocity of the surface of the body. The term $\displaystyle \gamma_{mom}$ is  a tensor accumulating several terms evaluated on fixed control volume given by,
\begin{equation}
\gamma_{mom}= \frac{\rho}{2} |\vec{u}|^{2} \mathbf{I} + \rho\left [(\vec{u_{s}} - \vec{u}) \vec{u}- \vec{u} (\vec{x} \times \vec{\zeta}) + \zeta (\vec{x} \times \vec{u})\right] - \rho\left[ \left( \vec{x}. \frac{\partial \vec{u}}{\partial t} \mathbf{I} - \vec{x} \frac{\partial \vec{u}}{\partial t} \right) \right] + \left[ \vec{x}. (\nabla . \mathbf{T}) \mathbf{I} - \vec{x} (\nabla. \mathbf{T}) \right] + \mathbf{T} \label{m2}
\end{equation} 
where $\mathbf{I}$ is the unit tensor and $\mathbf{T}$ is the viscous stress tensor $\displaystyle \mathbf{T}=\mu(\nabla \vec{u}+\nabla \vec{u}^T)$, $\mu$ being the dynamic viscosity of the fluid.

We choose $V$ to be a rectangular box of unit depth and net fluxes are calculated across its boundary. For a 2D flow, $V$ reduces to a rectangular region and $S$ reduces to a counter-clockwise rectangular curve (figure \ref{sch_dl}(b)).  Under these assumptions, if $F_D$ and $F_L$ are the drag and lift forces and making use of the fact that $\hat{n}ds=dy\hat{i}-dx\hat{j}$, making use of \eqref{m2} in \eqref{m1} and after some complicated vector algebra (see Appendix), we arrive at
\begin{eqnarray}
\label{dr_1}
\left [ \begin{array}{c} F_D \\
                         F_L
\end{array} \right ]
&=& -\rho\iint_V \left [ \begin{array}{c} \left(\frac{\partial u }{\partial t} + u \frac{\partial u}{\partial x} + v \frac{\partial u}{\partial y}\right)   \\
                            \left(\frac{\partial v }{\partial t} + u \frac{\partial v}{\partial x} + v \frac{\partial v}{\partial y}\right) 
\end{array} \right ]dxdy \nonumber \\
&+& \oint_S\left [ \begin{array}{c}  \rho \left(-v (u_{s} -u)+yv \zeta -y \frac{\partial v}{\partial t}\right)+\nu y \nabla^{2} u   \\
                         \rho \left( -\frac{1}{2} (u^2+v^2) - v (v_{s} - v) -xv \zeta - x \frac{\partial u}{\partial t} \right)  -\nu  \left( x \nabla^2 u +  \frac{\partial u}{\partial y} + \frac{\partial v}{\partial x} + 2 \frac{\partial v }{\partial y}   \right)
\end{array} \right ]dx \nonumber \\
&+& \oint_S\left [ \begin{array}{c}    \rho \left( \frac{1}{2} (u^2+v^2) + u (u_{s} - u) -yu \zeta - y \frac{\partial v}{\partial t} \right) +\nu \left( y \nabla^2 v + 2 \frac{\partial u}{\partial x} + \frac{\partial u}{\partial y} + \frac{\partial v }{\partial x}   \right)   \\
                          \rho  \left(u (v_{s} -v)+xu \zeta +x \frac{\partial u}{\partial t}\right)-\nu  x \nabla^{2} v 
\end{array} \right ]dy
\end{eqnarray}

Normalizing the drag and lift forces by the characteristic velocity $U$ and characteristic dimension $L$ of the bluff body under consideration, the drag and lift coefficients, viz., $C_D$ and $C_L$ respectively, reduces to
\begin{equation}
\label{dr_2}
\left [ \begin{array}{c} C_D \\
                         C_L
\end{array} \right ]
=\left [ \begin{array}{c} F_D/\left (\frac{1}{2}\rho U^2 L\right ) \\
                         F_L/\left (\frac{1}{2}\rho U^2 L\right )
\end{array} \right ]
\end{equation}
Subsequently, under the assumption of a solid bluff body and fixed control volume, and making use of \eqref{dr_1}, equation \eqref{dr_2} in terms of the non-dimensionalized variables can be written as
\begin{eqnarray}
\label{cdcl}
\left [ \begin{array}{c} C_D \\
                         C_L
\end{array} \right ]
= -2\iint_V \left [ \begin{array}{c} \frac{\partial u }{\partial t}    \\
                           \frac{\partial v }{\partial t}
\end{array} \right ]dxdy 
&+& 2\oint_S\left [ \begin{array}{c}   \left(uv+yv \zeta -y \frac{\partial v}{\partial t}\right)+\frac{1}{Re} y \nabla^{2} u   \\
                         \left( \frac{1}{2} (v^2-u^2) -xv \zeta - x \frac{\partial u}{\partial t} \right)  -\frac{1}{Re}  \left( x \nabla^2 u +  \frac{\partial u}{\partial y} + \frac{\partial v}{\partial x} + 2 \frac{\partial v }{\partial y}   \right)
\end{array} \right ]dx \nonumber \\
&+& 2\oint_S\left [ \begin{array}{c}     \left( \frac{1}{2} (v^2-u^2)  -yu \zeta - y \frac{\partial v}{\partial t} \right) +\frac{1}{Re} \left( y \nabla^2 v + 2 \frac{\partial u}{\partial x} + \frac{\partial u}{\partial y} + \frac{\partial v }{\partial x}   \right)   \\
                           \left(-uv+xu \zeta +x \frac{\partial u}{\partial t}\right)-\frac{1}{Re}  x \nabla^{2} v 
\end{array} \right ]dy
\end{eqnarray}
Both the line and surface integrals in the above formula are numerically computed by Trapezoidal rule.
\subsection{Solution of the System of Algebraic Equation}  
  The N-S equations \eqref{omega}-\eqref{psi}  in $\psi$-$\zeta$ formulation can easily be recast into the parabolic equation \eqref{s1}. For example, \eqref{omega} can be obtained from \eqref{s1} by setting $\lambda=1$, $u=\zeta$, $\beta=-1/Re$, $f=0$, $\beta_x=u$ and $\beta_y=v$. Likewise, \eqref{psi} is nothing but the steady-state version of \eqref{s1} with $u=\psi$, $f=-\zeta$ and $\beta_x=\beta_y=0$. Equation \eqref{omega} is discretized at regular points by employing \eqref{s4} and at irregular points by the procedure described in sections \ref{sectx} and \ref{secty}. \eqref{psi} is discretized by the procedure developed by the authors in their recent work \cite{singhalkalita}.

Once vorticity $\zeta$ and streamfunction $\psi$ have been computed, making use of equation \eqref{vel}, HOC approximation the velocities $u$ and $v$ are given by \cite{kalita2001fully}
\begin{equation}
u_{ij}=\delta_y \psi +   \frac{l^{2}}{6} \left(\delta_y \zeta + \delta_x^2 \delta_y \psi \right) + O(h^4,l^4)\label{u_vel}
\end{equation}
\begin{equation}
v_{ij}=-\delta_x \psi - \frac{h^{2}}{6} \left(\delta_x \zeta + \delta_x \delta_y^2 \psi \right) + O(h^4,l^4)\label{v_vel}
\end{equation}
In matrix form, the discretized HOC form of the vorticity transport equation (\ref{omega}) at the interior nodes is given by
\begin{equation}
\tilde{S} \zeta^{(n+1)} = \tilde{S^{'}} \zeta^{(n)} + \tilde{C}_{\zeta}^{(n)}+ \tilde{C}_{\zeta}^{(n+1)}\label{mzeta}
\end{equation} 
where $\tilde{S}$, $\tilde{S^{'}}$ is HOC matrix to the equation (\ref{omega}) and $\tilde{C}_{\zeta}^{(n)}, \;\tilde{C}_{\zeta}^{(n+1)}$ are the vorticity correction vectors at the irregular points corresponding to the $n^{\rm th}$ and $(n+1)^{\rm th}$ time levels respectively.  Likewise, the HOC matrix representation of equation (\ref{psi}) at the interior nodes is 
\begin{equation}
\tilde{T} \psi^{(n+1)}+ \left[ I+ \frac{h^2}{12} T \right] \zeta^{(n+1)} =0 \label{mpsi}
\end{equation} 
where $\tilde{T}$ and T are the HOC and CDS matrices corresponding to the interior nodes.
On the other hand, the boundary conditions \eqref{bc_psi} may be expressed
\begin{equation}
N \psi^{(n+1)} + \tilde{B} \zeta^{(n+1)} = \tilde{U}^{(n+1)}
\end{equation}
\begin{equation}
\psi_{B}=0
\end{equation}
where $N$ is the matrix resulting from normal derivative boundary conditions, $\tilde{B}$ is the vorticity boundary matrix, $\tilde{U}^{(n+1)}$ is the current velocity vector owing to \eqref{u_vel}-\eqref{v_vel}, and the subscript $I$, and $B$ denotes the interior and boundary respectively. The following matrix equations provide a complete picture of the discretized equations at the regular, irregular and the boundary points simultaneously in concise and compact form
\begin{equation}
\begin{bmatrix}
\tilde{T}_{R} & \tilde{T}_{IR} & T_{B} \\
O & O & I 
\end{bmatrix}
\begin{bmatrix}
\psi_{R}\\
\psi_{IR}\\
\psi_{B}
\end{bmatrix}^{(n+1)}
 = -
\begin{bmatrix}
I+\frac{h^{2}}{12} & \frac{h^{2}}{12} T_{B}  \\
O & O 
\end{bmatrix}
\begin{bmatrix}
\zeta_{I}\\
\zeta_{B}
\end{bmatrix}^{(n)}+
\begin{bmatrix}
\tilde{C}_{\psi}\\
0
\end{bmatrix}^{(n)}\label{block1}
\end{equation}
\begin{equation}
\begin{bmatrix}
\tilde{S}_{R} & \tilde{S}_{IR} & \tilde{S}_{B} \\
\tilde{B}_{B} & O & \tilde{B}_{B} 
\end{bmatrix}
\begin{bmatrix}
\zeta_{R}\\
\zeta_{IR}\\
\zeta_{B}
\end{bmatrix}^{(n+1)}
 = -
\begin{bmatrix}
O & O  \\
N_{I} & N_{B} 
\end{bmatrix}
\begin{bmatrix}
\psi_{I}\\
\psi_{B}
\end{bmatrix}^{(n+1)}+
\begin{bmatrix}
F\\
\tilde{U}
\end{bmatrix}^{(n)}\label{block2}
\end{equation}
where,
\begin{equation}
F = \begin{bmatrix}
\tilde{S}'_{R} & \tilde{S}'_{IR} & \tilde{S}'_{B} 
\end{bmatrix}
\begin{bmatrix}
\zeta_{R}\\
\zeta_{IR}\\
\zeta_{B}
\end{bmatrix}^{(n)}+
(\tilde{C}_{\zeta}^{(n+1)}+\tilde{C}_{\zeta}^{(n)}) \label{block3}
\end{equation}
and the subscripts $R,\; {\rm and}\; IR$ represent the regular and irregular interior points respectively.

For a grid of size $M \times N$, the matrices $\tilde{S}$, $\tilde{S}'$, $\tilde{T}$ and $T$ are of order $MN$ and $\zeta^{(n)}$, $\zeta^{(n+1)}$, $\psi^{(n)}$, $\psi^{(n+1)}$, $\tilde{C}^{(n)}$, $\tilde{C}^{(n+1)}$ are vectors of length $MN$ in equations \eqref{mzeta} and \eqref{mpsi}. Apart from the grid size of the computational domain, the size of the block matrices in equations \eqref{block1}-\eqref{block3} depends on the geometry of the immersed body which determines the number of regular and irregular points.

An inner-outer iteration procedure is essential for the time marching solutions of transient fluid flow problems governed by \eqref{omega}-\eqref{psi}. Once $u$, $v$, $\zeta$ and $\psi$ are presented with appropriate initial and boundary conditions, firstly  \eqref{omega} and then \eqref{psi} is solved. Once $\psi$ is available, $u$ and $v$ are computed by utilizing \eqref{u_vel} and \eqref{v_vel}. This completes one outer time iteration.

The inner iterations are composed of solving \eqref{mzeta} and \eqref{mpsi} by efficient iterative solvers at each time step. We have accomplished this in our computations by employing the BiCGStab Stabilized \cite{kelley1995iterative} iterative solver along with Incomplete LU decomposition as preconditioner with the help of the Lis Library \cite{lis}. The inner iterations were stopped when the residual vectors arising out of equations \eqref{mzeta} and \eqref{mpsi} fell below $10^{-13}$. All our computations were performed on a Intel Xeon processor-based PC with a 32 GB RAM.  
\section{Numerical test cases}
In order to demonstrate the effectiveness of the proposed approach, it is applied to several problems. The first of these has analytical solution and the remaining are flow past bluff bodies immersed in fluids, for both the stationary and moving cases. Also considered are flows involving multiple bodies. In all the fluid problems under consideration, the flow is governed by the unsteady N-S equations for incompressible viscous flows.    
 \subsection{Test Case 1: Interface problem having analytical solution}  
As our first test case, we validate our algorithm to solve the parabolic equation given by
\begin{equation}
u_t= \nu \nabla^2 u\quad (x,y,t)\in \Omega \times (0,T]    
\end{equation}
with initial and boundary conditions 
$$\displaystyle u(x,y,0)=u_0(x,y), \quad \quad (x,y)\in \Omega,$$
$$\displaystyle u(x,y,t)=u_b(x,y,t), \quad (x,y,t)\in \partial\Omega \times (0,T], $$
The computational domain $\Omega$ is the square $[0,1] \times [0,1]$ and the solution has a discontinuity across the the circular interface $\Gamma$ of radius $0.25$ centered at the point $(0.5,0.5)$. The level set function and analytical solution are respectively defined by $\phi = (x-0.5)^2 + (y-0.5)^2 - (0.25)^2$ and 
\begin{equation}\label{ex1}
u(x,y,t)=\left\{\begin{array}{rr}
&e^{-t \nu \pi^2 (k_x^2+k_y^2)} cos(k_x \pi x) cos(k_x \pi y), \quad \phi \geq 0\\
&0, \quad
 \phi < 0.
\end{array}\right.
\end{equation}
respectively, where $(k_x,k_y)$ are wave number and set the values of both are 2. The initial and boundary conditions are obtained from the equation \eqref{ex1}. We apply a no flux condition on the boundary of the circular interface i.e
\begin{equation}
\frac{\partial u}{\partial \textit{n}}\equiv \nabla u . \hat{n}= -4\pi e^{-t \nu \pi^2 (k_x^2+k_y^2)} (k_x sin(k_x \pi x) (x-0.5)+ k_y sin(k_x \pi y) (y-0.5))
\end{equation}
where $\hat{n}=(x-0.5,y-0.5)$ is the normal vector to the circle and the jump conditions  approximated using equations \eqref{aj1}-\eqref{aj6}.
\begin{table}[!h]
\caption{ Grid refinement analysis of maximum error for Test Case 1 at $t=2.5$ for $\nu=1/200$ with $\triangle t= 10^{-3}$.}
\begin{center}
\begin{tabular}{|c|c|c|c|c|}  \hline
N & Present (k=2) & ROC & Present (k=3) & ROC  \\ \hline
20  & $ 2.57 \times 10^{-4}$ & $-$  & $ 3.96 \times 10^{-5}$ & $-$ \\
40  & $ 7.11 \times 10^{-5}$ & $1.85$ & $ 5.61 \times 10^{-6}$ & $2.81$ \\
80  & $ 9.82 \times 10^{-6}$ & $2.84$ & $ 4.59 \times 10^{-7}$ & $3.61$ \\
160 & $ 1.61 \times 10^{-6}$ & $2.60$ & $ 3.21 \times 10^{-8}$ & $3.83$ \\
320 & $ 2.25 \times 10^{-7}$ & $2.83$ & $ 2.38 \times 10^{-9}$ & $3.75$ \\
  \hline 
  \end{tabular}
\end{center}\label{ex1_table1}
\end{table}

In table \ref{ex1_table1}, we present the maximum error $\parallel E_{N}\parallel_{\infty}$ resulting from our computation on gradually increasing grid sizes $N \times N$ and show the effect of the parameter $k$ described in equation \eqref{irr_reg}. Expectedly, a higher value of $k$ yields a better convergence rate (ROC), which is defined as
 $$\textnormal{Order}=\frac{log(\parallel E_{N}\parallel_{\infty} / \parallel E_{(N/2)}\parallel_{\infty} )}{log(2)}, $$
where $\parallel E_{N/2}\parallel_{\infty}$ is the maximum error of the previous coarser grid having half the points in either direction than the current grid size. While the errors resulting from our computation corresponding to $k=3$ decay at a rate close to four, the ones from the simulation of Calhoun \cite{calhoun1999cartesian} could obtain an ROC close to two only.

\begin{figure}[!h]
\begin{minipage}[b]{.49\linewidth}  
\includegraphics[scale=0.42]{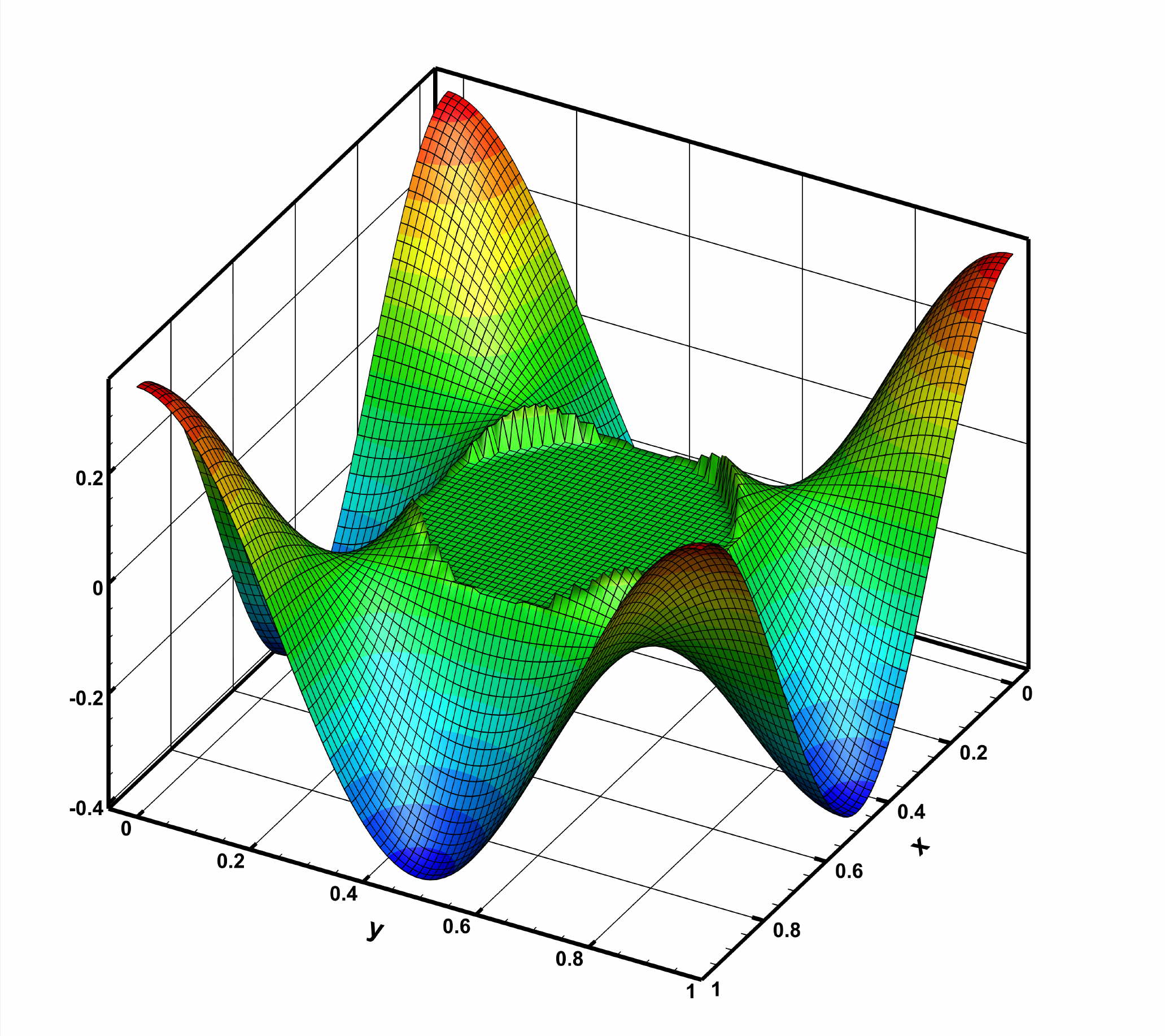} 
\end{minipage}           
\begin{minipage}[b]{.49\linewidth}
\includegraphics[scale=0.42]{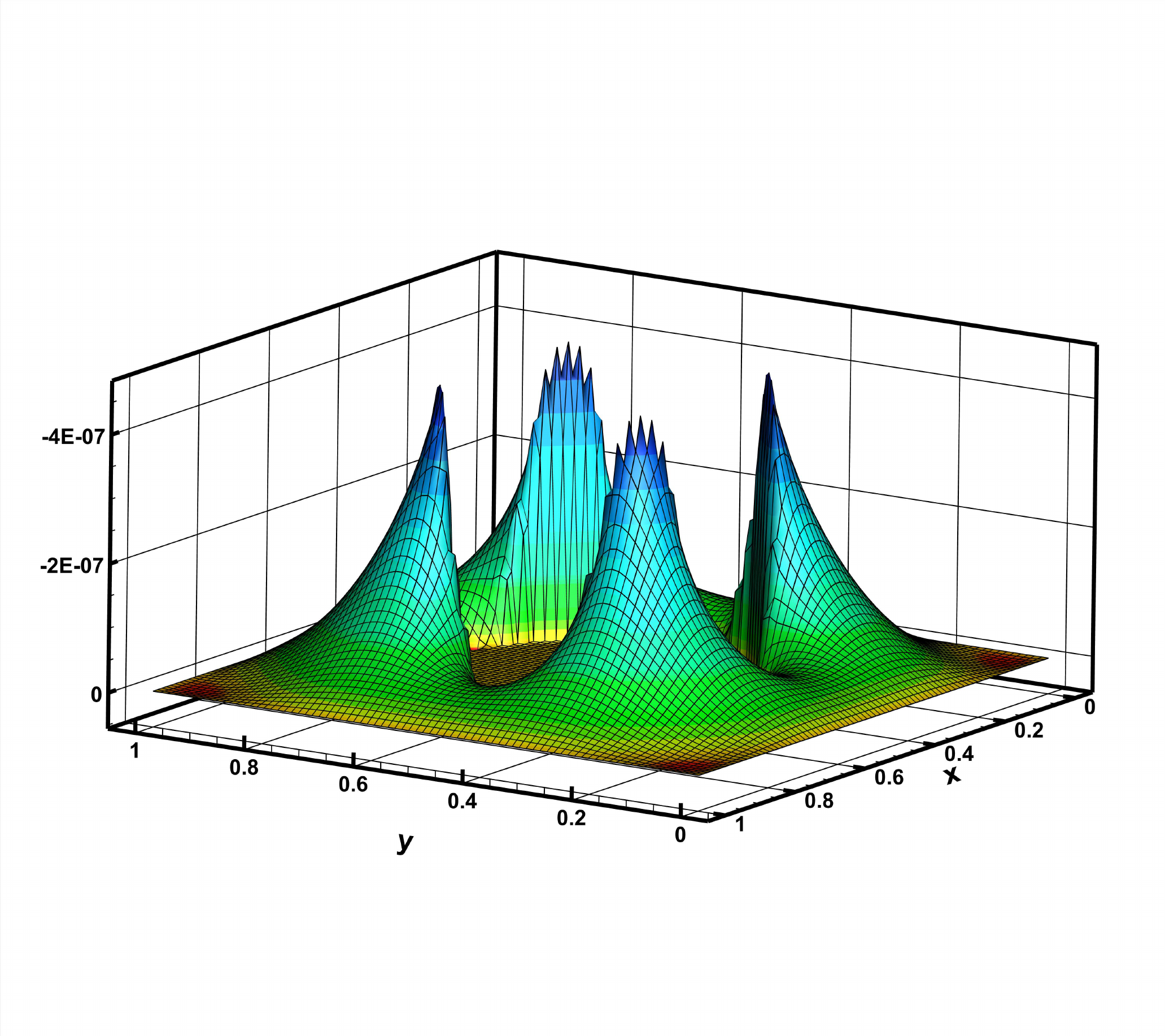} 
\end{minipage} 
\caption{{\sl  Surface plots of (a) the numerical solution and (b) error on a grid of size $80 \times 80$ for Test Case $1$ for $k=3$.} }
\label{case1}
\end{figure}

We also present the surface plots of our numerical solution on a grid of size $80 \times 80$ side by side with the surface plots of errors in figures \ref{case1}(a)-(b). Figure \ref{case1}(a) clearly demonstrates that the sharp interface has been resolved very efficiently by our approach. Note that the errors from our computation (see figure \ref{case1}(b)) are much lower in magnitude the ones in  \cite{calhoun1999cartesian}. 
\subsection{Flow Past Stationary Bluff Bodies}
\begin{figure}[!h]
 \begin{center} 
\includegraphics[scale=0.65]{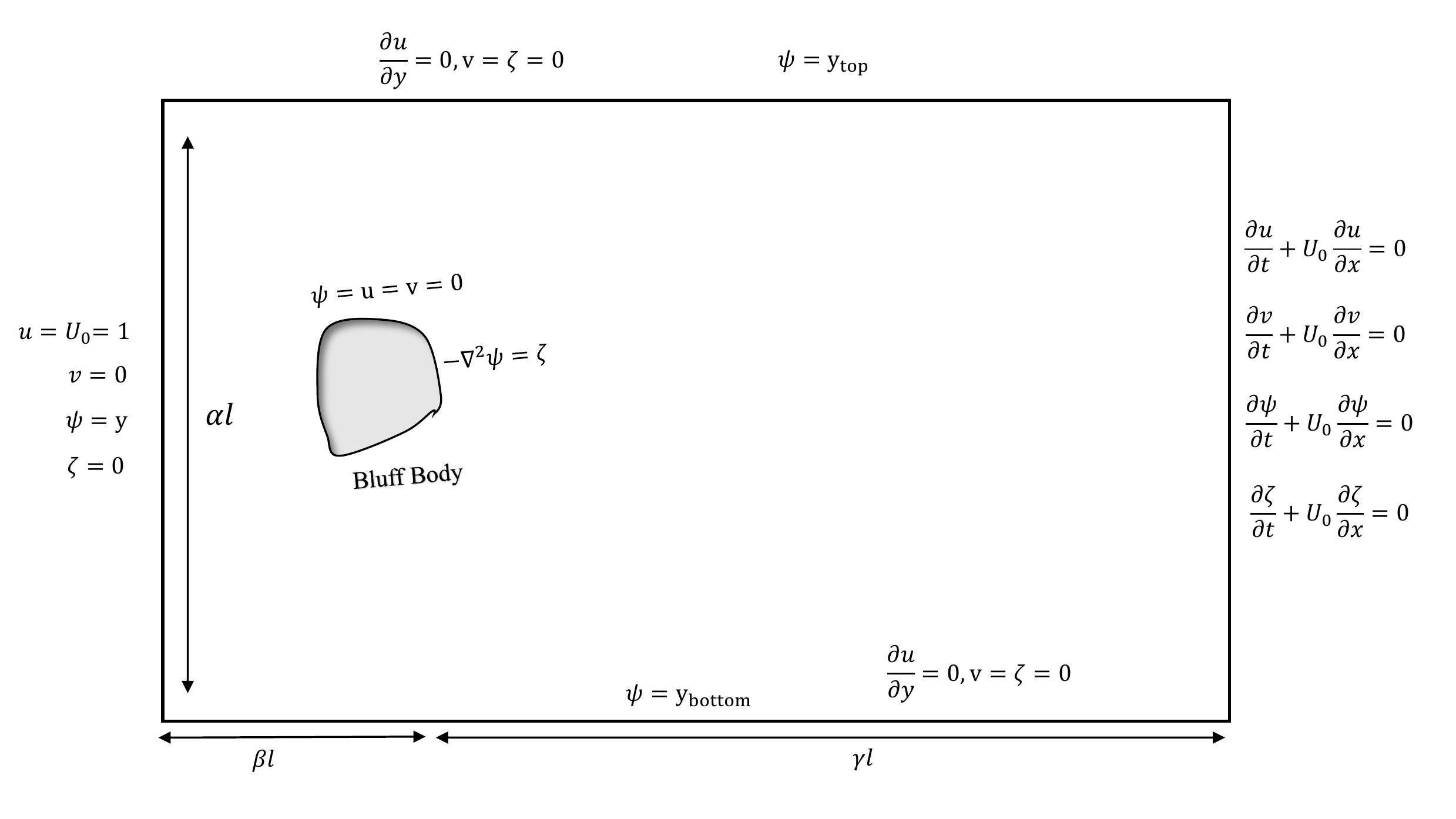}  
 \end{center}
 \caption{{\sl Schematic of the flow past bluff bodies in uniform flow.} }
\label{setup_bluff}
\end{figure} 
The study of the flow past bluff bodies holds an important place in many engineering applications, particularly in the field of naval architecture. Such flows are very complex and are highly characterised by the generation and shedding of vortical structures \cite{ chi2020directional, kalita2009transformation, kumar2021comprehensive, linnick2005high, park2016pre, xu2006immersed}. This section is concerned not only with the simulation of flow past stationary bluff bodies, but also with moving ones immersed in fluids.  Moreover flow situations involving multiple bodies are also considered. The problems have been chosen in such a way that the ability of the current approach in handling complicated geometry and varied flow situations can be established. As would be seen later on, while most of the previous studies involved computations either by finite volume or finite element approach in extremely finer grids  \cite{babu2008aerodynamic, behara2019flow, nepali2020two,   papaioannou2008effect, prasanth2009vortex,shaaban2018flow, zhao2014two} for this kind of flows, our approach accomplishes the same in relatively coarse grids, that too, in FD set-up.

In figure \ref{setup_bluff}, we show a schematic of the computational domain along with the boundary conditions used for the simulation. While choosing the dimensions, ample care was taken so that the simulation is free from any entrance effects and there is no hindrance in the smooth shedding of the vortices once the vortices formed on the surface of the bodies start detaching from them. Besides, in all the computations time-steps are chosen in the range $10^{-2} \leq \Delta t \leq 10^{-3}$ according to the flow situation.

\subsubsection{Test case 2: Flow Past a Stationary Circular Cylinder}\label{st_cyl}
In our first test case for flow past stationary bluff bodies, we consider the problem of flow around an impulsively started circular cylinder in a free-stream with uniform velocity. There exists an enormous number of numerical and experimental results for this problem and as such, is a perfect test case for examining the efficiency of the proposed approach by comparing the results obtained from our computations with the benchmark results available in literature. This problem also act as a prelude to tackling problems with moving immersed interfaces. 

The schematic for this problem has already been shown in figure \ref{setup_bluff} where the bluff body is now the circular cylinder. Here, Reynolds number is described as $\displaystyle Re=\frac{U_{0}l}{\nu}$, where $l$ is the cylinder diameter, $U_{0}$ is the free stream velocity at the inlet,  and $\nu$ is the kinematic viscosity of the fluid. For our simulations, we assume $l$ is to be $1.0$ with the center of cylinder fixed at $(0, 0)$. We have chosen $\beta=5.0$, $\alpha=10.0$, and $\gamma=25.0$ in figure \ref{setup_bluff} such that the corresponding dimensions of the computational domain are $-5.0 \leq x \leq 25.0$ and $5.0 \leq y \leq 5.0$; as such $y_{top} =5.0$ and $y_{bottom} =-5.0$. The boundary conditions at the far-field and on the surface of the cylinder are as follows:

\begin{itemize}
\item $u=1$, $v=0$, $\psi=y$ and $\zeta=0$ on the left wall $x=-5.0$,
\item $\displaystyle \frac{\partial u}{\partial y}=0$, $v=0$, $\psi=y_{bottom}$ and $\zeta=0$ on the bottom wall $y=-5.0$,
\item $\displaystyle \frac{\partial u}{\partial y}=0$, $v=0$, $\psi=y_{top}$ and $\zeta=0$ on the top wall $y=5.0$,
\item $\displaystyle \frac{\partial u}{\partial t}+U_0\frac{\partial u}{\partial x}=0$, $\displaystyle \frac{\partial v}{\partial t}+U_0\frac{\partial v}{\partial x}=0$, $\displaystyle \frac{\partial \psi}{\partial t}+U_0\frac{\partial \psi}{\partial x}=0$ and $\displaystyle \frac{\partial \zeta}{\partial t}+U_0\frac{\partial \zeta}{\partial x}=0$ on the right wall $x=25.0$.
\item $\displaystyle u=v=\psi=0$, $\zeta=-\nabla^2\psi$ on the surface of the cylinder.
\end{itemize}

\begin{figure}[!h]
\begin{minipage}[b]{.5\linewidth}  
\includegraphics[scale=0.425]{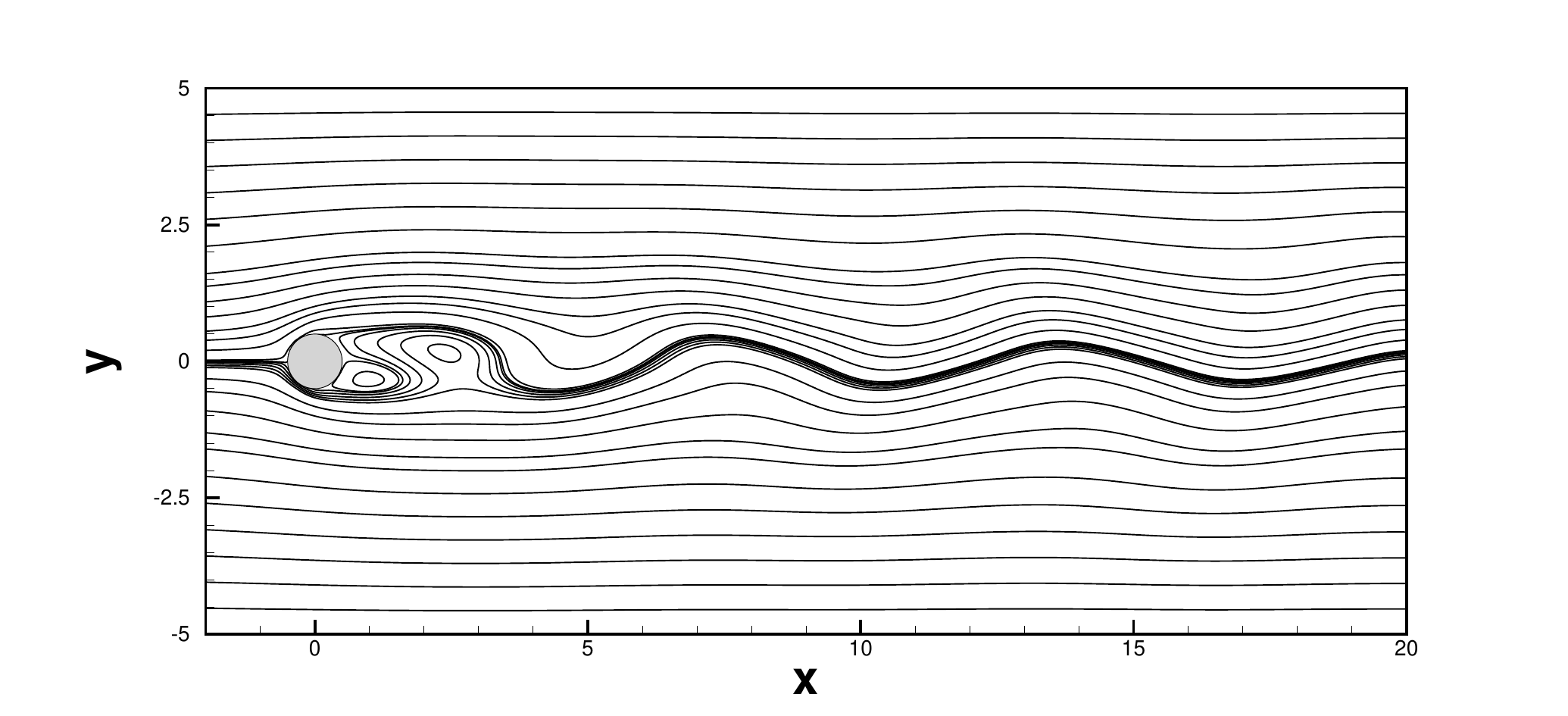} 
 \centering (a)
\end{minipage}            \hspace{-2.mm}
\begin{minipage}[b]{.5\linewidth}
\includegraphics[scale=0.425]{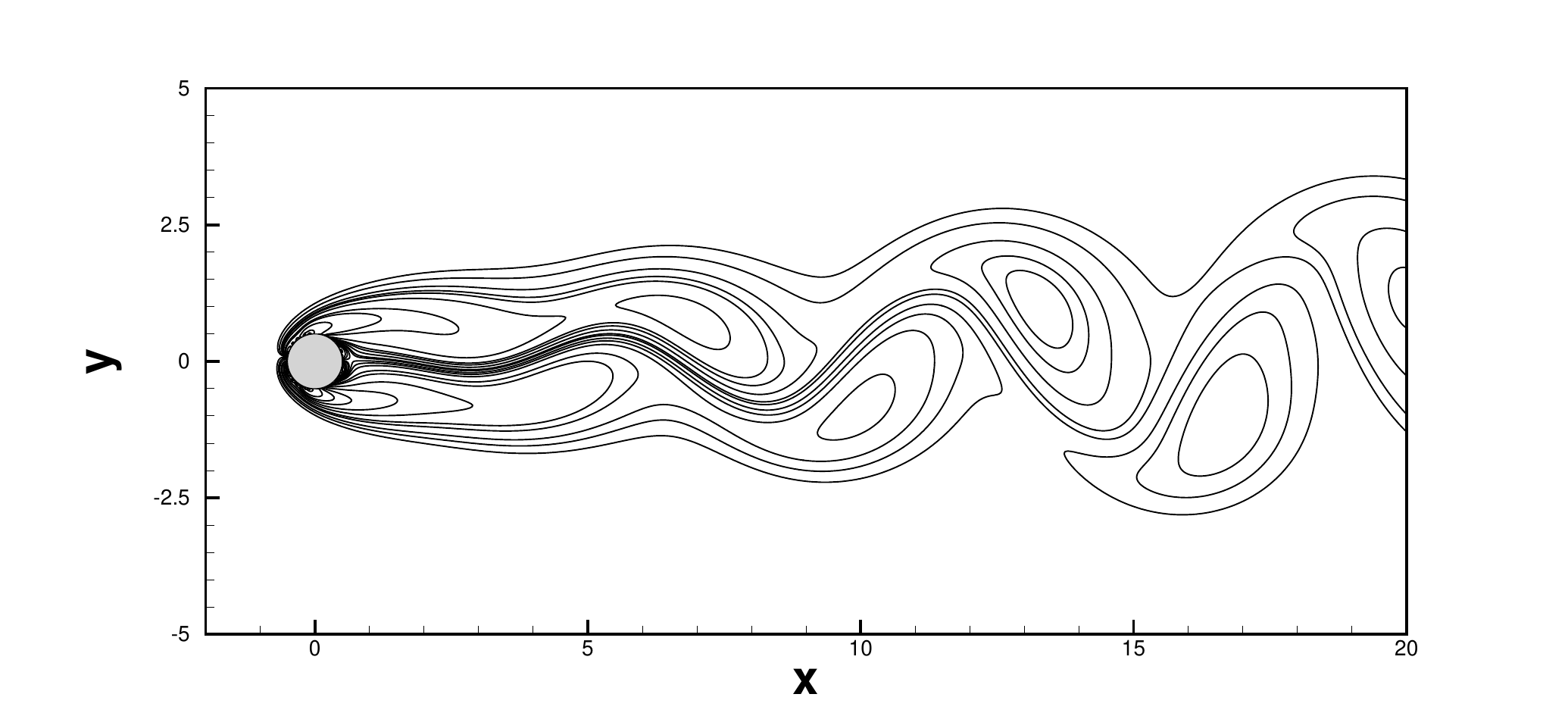} 
 \centering (a)
\end{minipage} 
\begin{minipage}[b]{.5\linewidth}  
\includegraphics[scale=0.425]{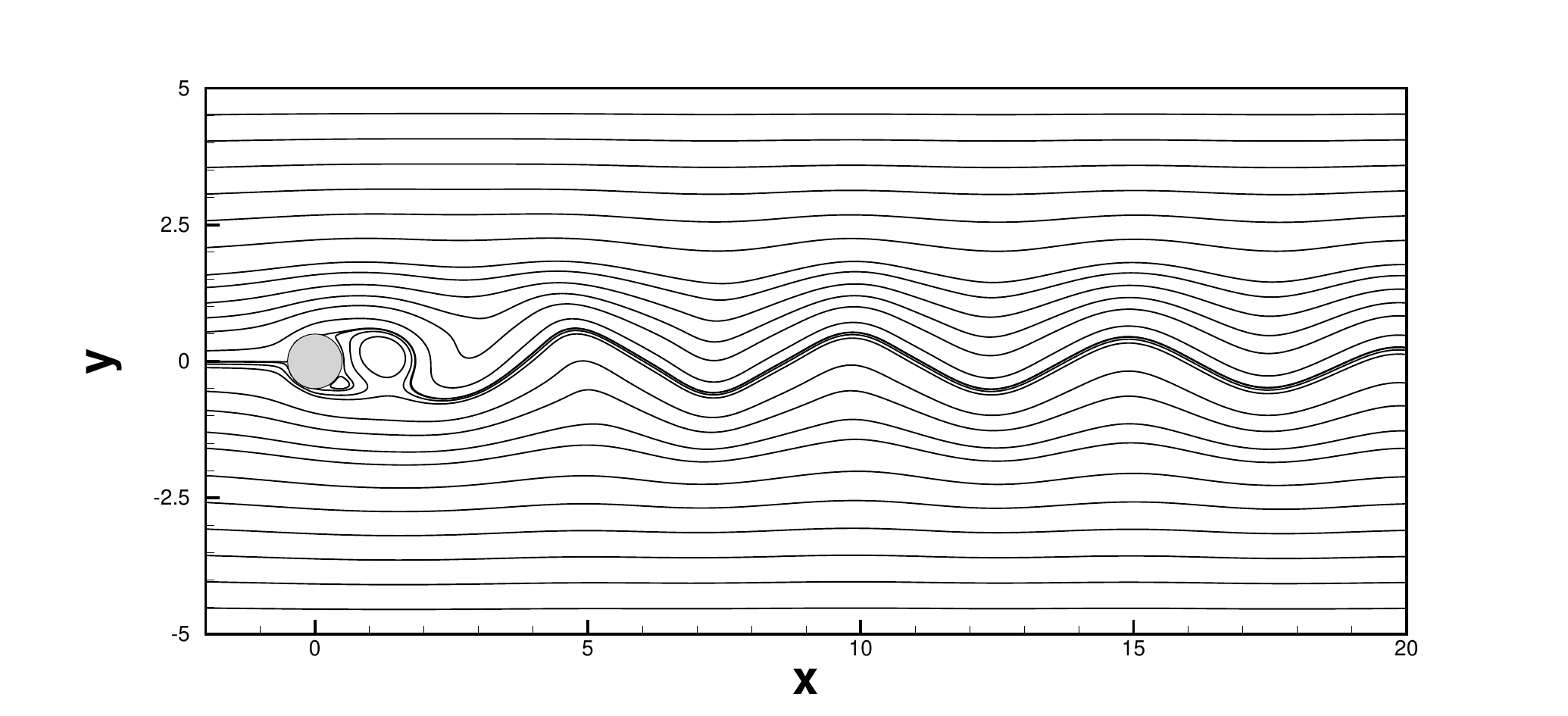} 
 \centering (b)
\end{minipage}            \hspace{-2.mm}
\begin{minipage}[b]{.5\linewidth}
\includegraphics[scale=0.425]{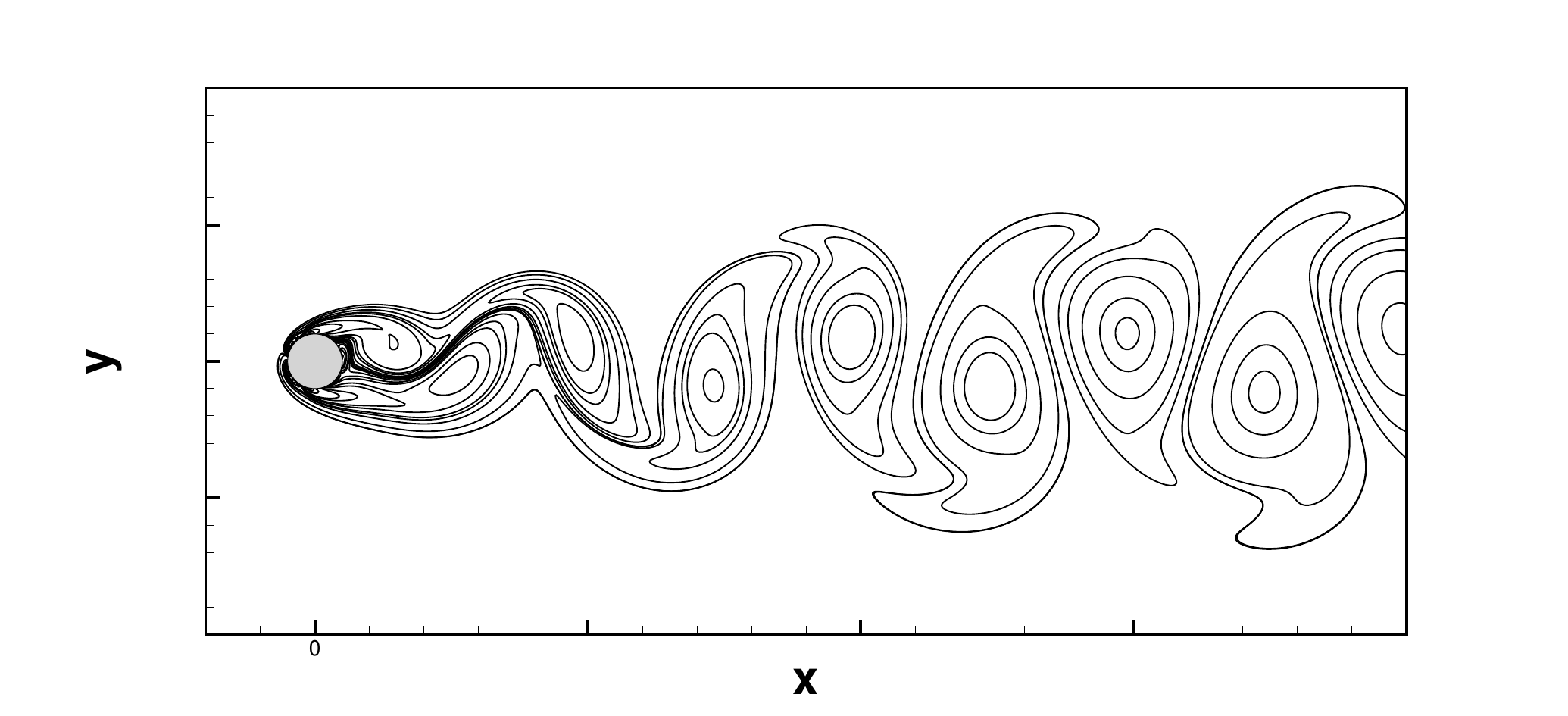} 
 \centering (b)
\end{minipage} 
\begin{minipage}[b]{.5\linewidth}   
\includegraphics[scale=0.425]{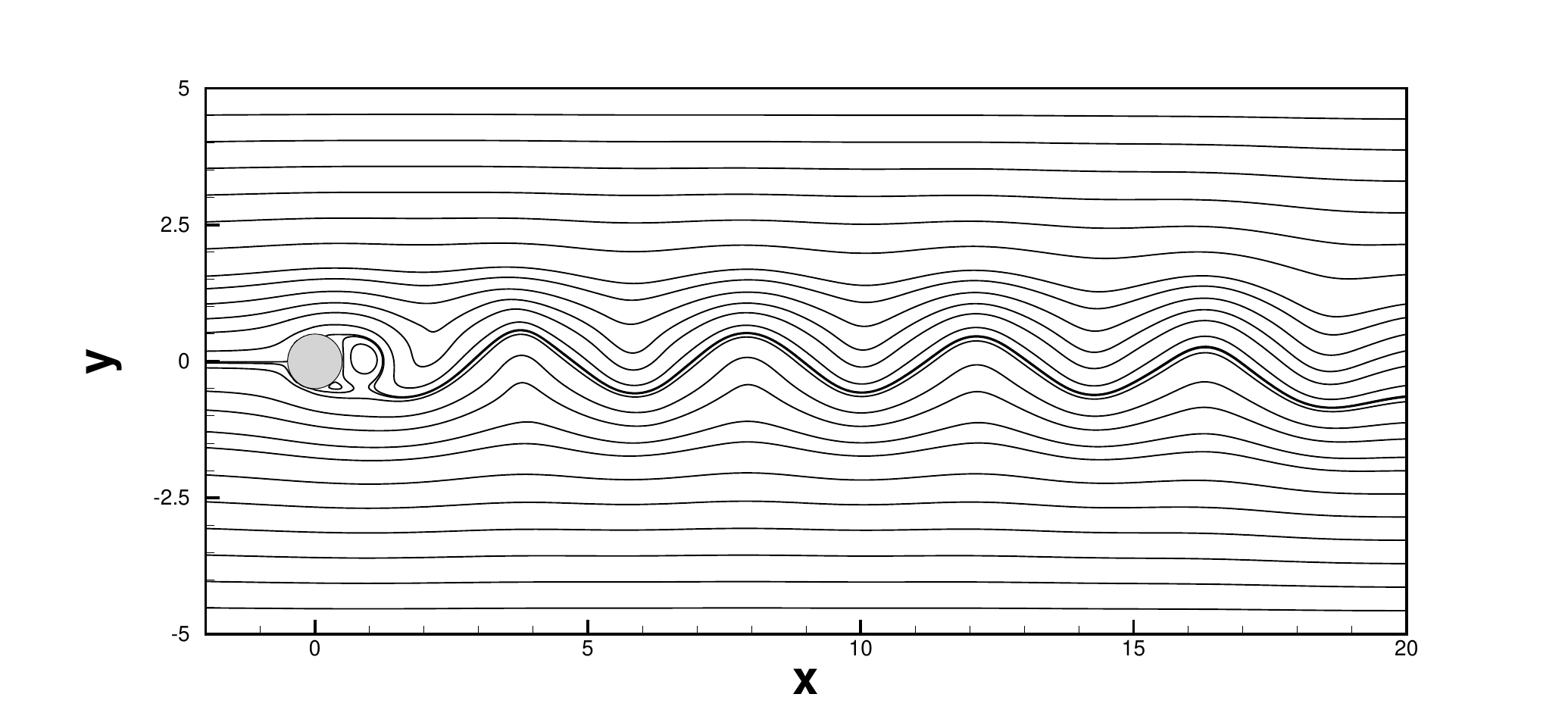} 
 \centering (c)
\end{minipage}          \hspace{-2.mm}
\begin{minipage}[b]{.5\linewidth}
\includegraphics[scale=0.425]{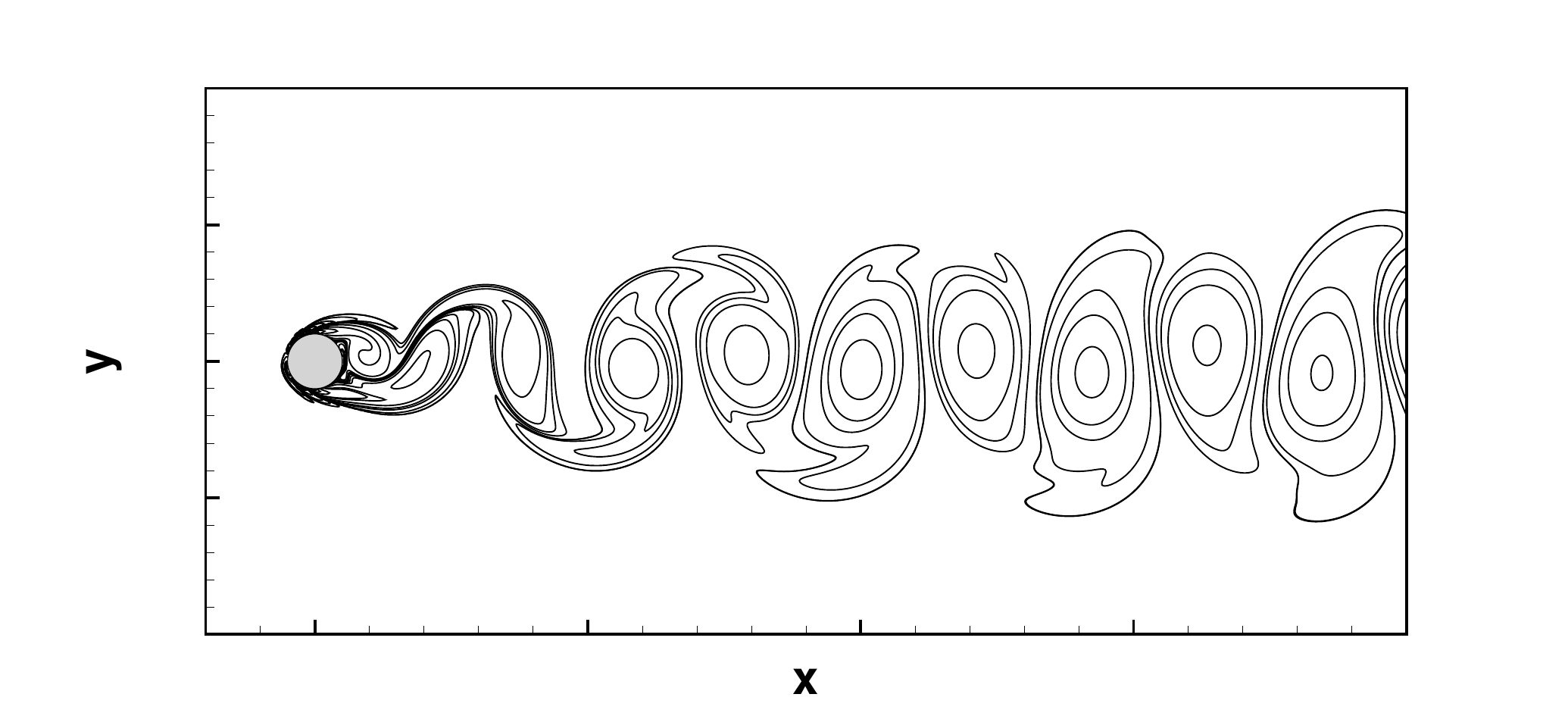} 
 \centering (c)
\end{minipage} 
\caption{{\sl  Simulation of flow past circular cylinder problem by present method: Streamlines (left) and Vorticity contours (right) for (a) $Re=50$, (b) $Re=100$ and (c) $Re=200$.} }
\label{sf_vt}
\end{figure}

Literature suggests that the flow for an impulsively started stationary cylinder becomes unsteady beyond a critical Reynolds number $45 \leq Re_C \leq 48$. In the following, we present our computational results for Reynolds numbers $Re=50$, $100$ and $200$. Note that in many studies, when the Reynolds number under consideration is slightly above $R_C$ as in the case of $Re=50$, the flow is artificially perturbed  \cite{dipankar2007suppression, russell2003cartesian} in order to break the symmetry of the flow. However, in the computation through our approach, asymmetry sets in naturally without the need of such perturbation. The flow for the range of $Re$ chosen, eventually becomes periodic and is fraught with the vortex shedding phenomenon characterised by the existence of von K$\acute{\rm a}$rm$\acute{\rm a}$n vortex street. Once shedding process starts some times after the flow symmetry is broken about $y=0$ line, vortices are shed alternatively from the either side of the $y=0$ line in a regular fashion. We depict this process in figure  \ref{sf_vt} where the left panel shows the instantaneous streamlines and the right, the vorticity contours  for   $Re=50$ (top), $100$ (middle) and $200$ (bottom) respectively. As one can see from these figures, shedding becomes more prominent with increase in $Re$ value.  

\begin{figure}[!h]
 \begin{minipage}[b]{.45\linewidth}  
\includegraphics[scale=0.4]{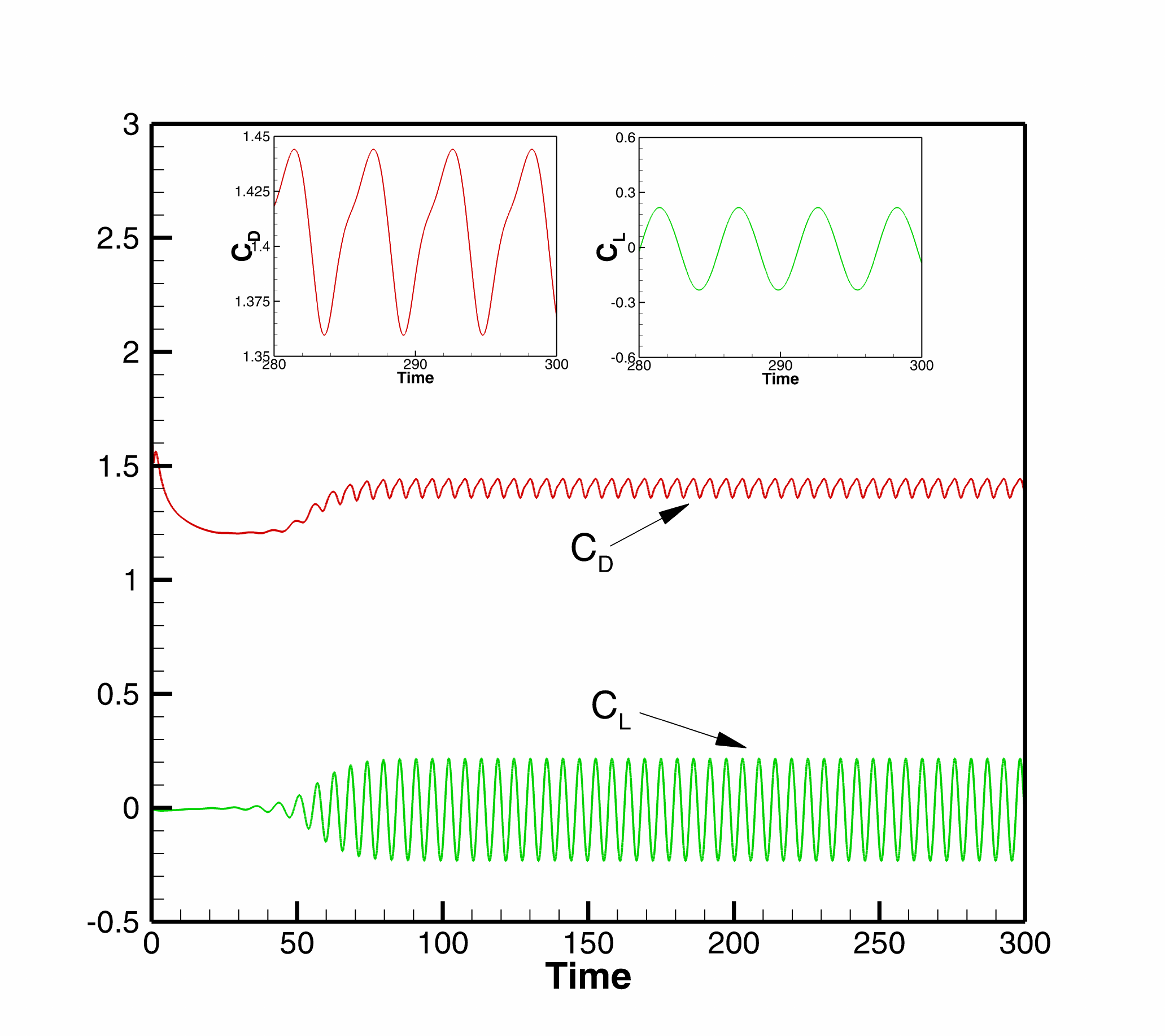} 
\centering (a) 
\end{minipage}           
\begin{minipage}[b]{.45\linewidth}
\includegraphics[scale=0.4]{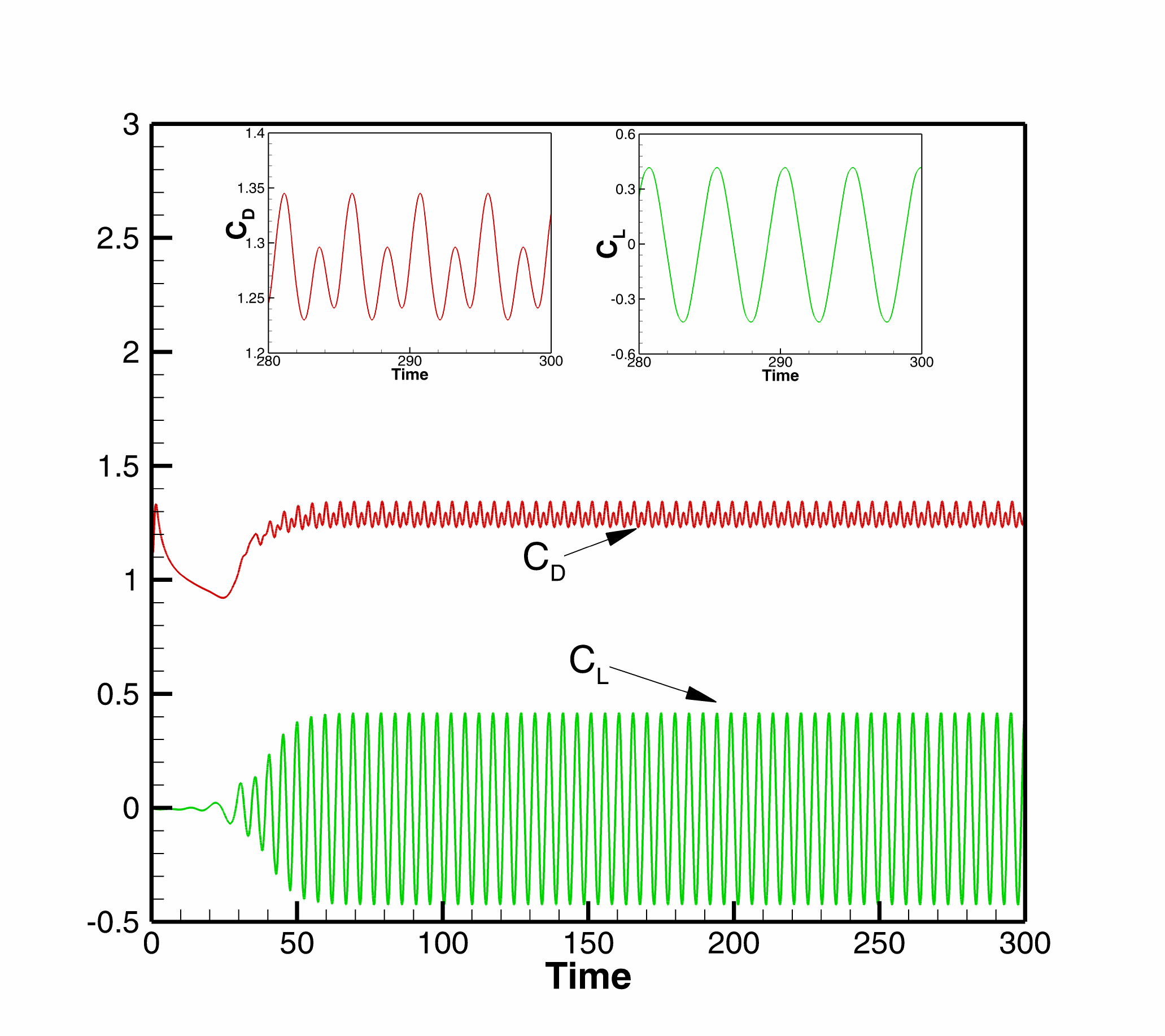}
\centering (b) 
\end{minipage} 
 \caption{{\sl History of drag and lift coefficients for the flow past an impulsively started stationary cylinder for (a) $Re=100$ and (b) $Re=200$.} }
\label{drag_lift}
\end{figure}

We further compute the drag and lift coefficients $C_D$ and $C_L$ respectively by the formulas \eqref{cdcl} and the Strouhal number $St$, which describes the mechanism of the oscillatory flow during the shedding process. It is defined as $\displaystyle St=\frac{fD}{U_0}$, where $f$ is the dominant frequency of the periodic lift variations, extracted from a time sample of $C_L$s. In figures \ref{drag_lift}(a)-(b), we plot the time histories of the drag and lift coefficients for $Re=100$ and $200$ respectively. We also compare our computed Strouhal numbers, drag and lift coefficients for the same Reynolds numbers with established experimental and numerical results in table \ref{dl} and obtain excellent comparison.     
\begin{table}[]
\centering
\caption{ Comparison of Strouhal number, drag and lift coefficients of the periodic flow for $Re =100$ and $200$.}
\begin{tabular}{|c|ccc|ccc|}
\hline
$Re$                               & \multicolumn{3}{c|}{$100$}                & \multicolumn{3}{c|}{$200$}                               \\ \hline
Reference &
$St$ &
\multicolumn{1}{c|}{$C_D$} &
$C_L$ &
\multicolumn{1}{c|}{$St$} &
\multicolumn{1}{c|}{$C_D$} &
$C_L$ \\ \hline
Frank \it{et al.} \cite{franke1990numerical}& -       & -                 & -           & $0.194$ & $1.31$            & $\pm 0.65$    \\
Williamson \cite{williamson1996vortex}                         & $0.163$ & -                 & -           & $0.185$ & -                 & -                  \\
Calhoun \cite{calhoun2002cartesian}                           & $0.175$ & $1.330 \pm 0.014$            &$\pm 0.298$           & $0.202$ & $1.172 \pm 0.058$                 &$\pm 0.668$                  \\
Le {\it{et al.}} \cite{le2006immersed}    & $0.160$ & $1.37 \pm 0.009$  & $\pm 0.323$ & $0.187$ & $1.34 \pm 0.030$  & $\pm 0.430$ \\
Berthelsen and Faltinsen \cite{berthelsen2008local}          & $0.169$ & $1.38 \pm 0.010$  & $\pm 0.340$ & $0.200$ & $1.37 \pm 0.046$ & -   \\
Russel \& Wang   \cite{russell2003cartesian}  & $0.169$ & $1.380 \pm 0.007$           & $\pm 0.300$ & $0.195$ & $1.290 \pm 022$           & $\pm 0.708$  \\
S.Sen  \cite{sen20154oec}                          & $0.165$ & $1.394 \pm 0.007$ & $\pm 0.191$ & $0.197$ & $1.375 \pm 0.038$ & $\pm 0.500$ \\
Present Study & $0.180$ &$1.402 \pm 0.042$ &$\pm 0.232$ &$0.210$ &$1.288 \pm 0.058$ &$\pm 0.425$ \\ \hline
\end{tabular}
\label{dl}
\end{table}
\subsubsection{Test case 3: Flow Past a cactus shaped cylinder}
\begin{figure}[!ht]
\begin{minipage}[b]{.45\linewidth}  
\includegraphics[scale=0.425]{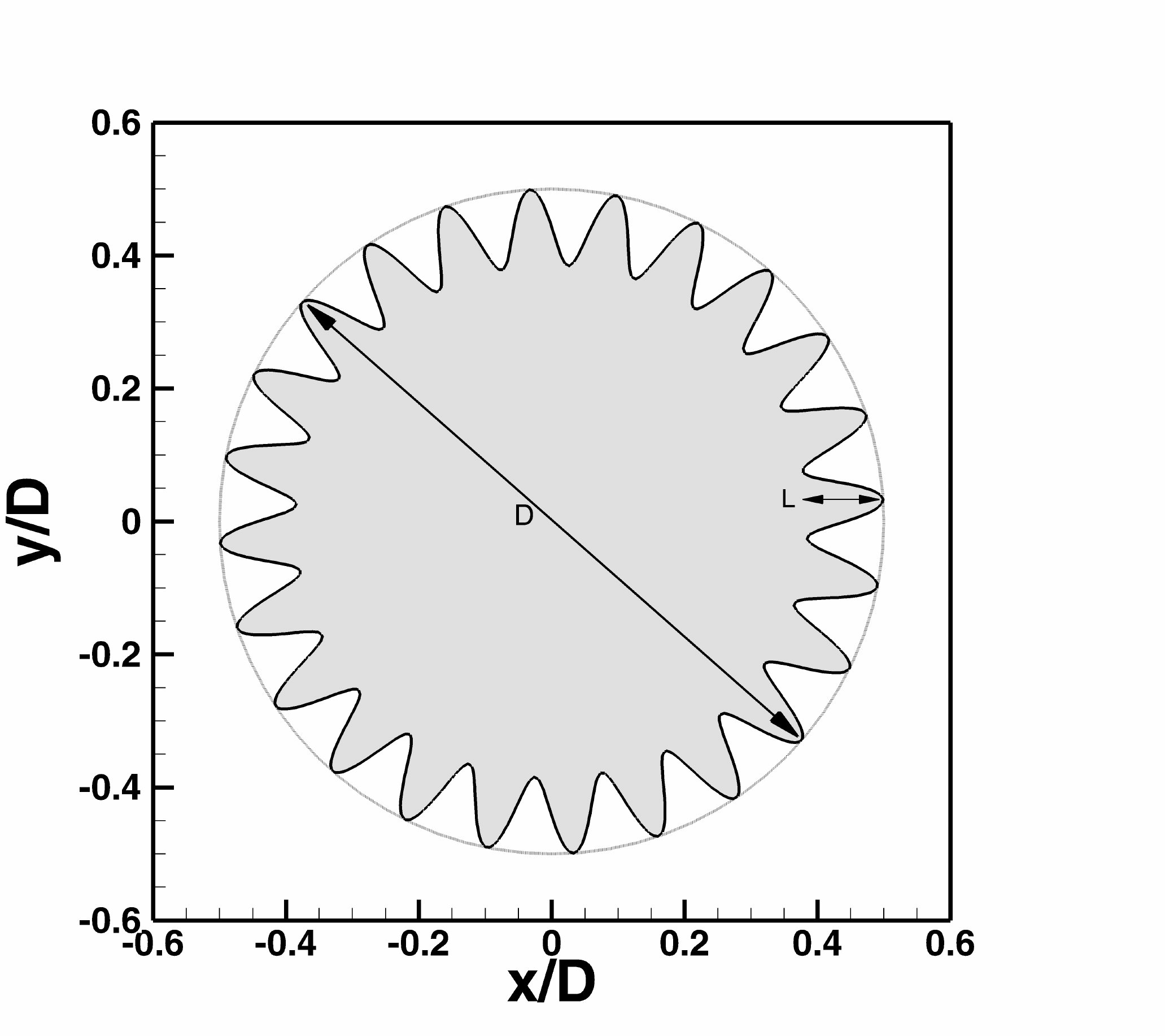} 
\centering (a) 
\end{minipage}           
\begin{minipage}[b]{.45\linewidth}
\includegraphics[scale=0.425]{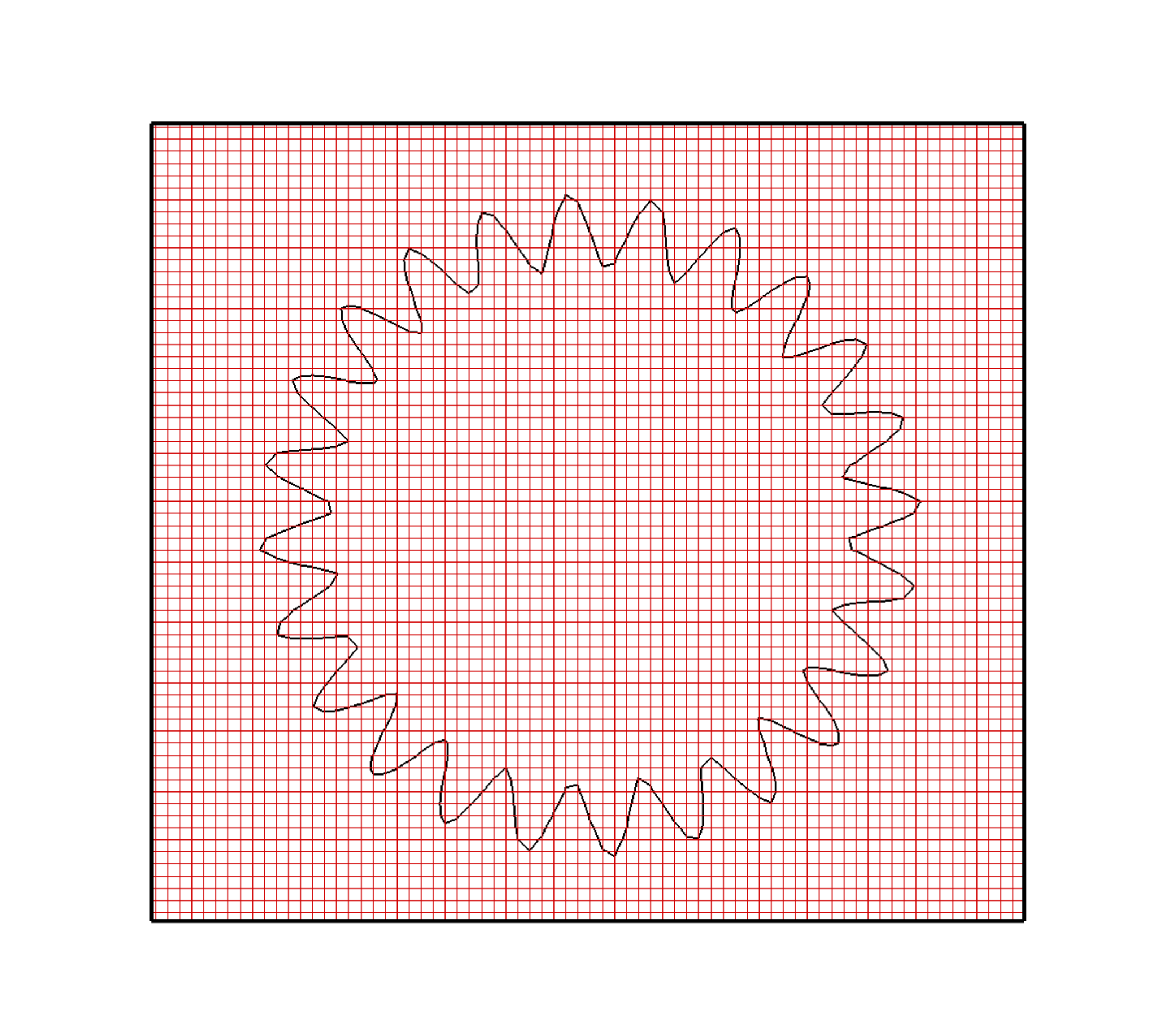}
\centering (b) 
\end{minipage} 
\caption{{\sl (a) Schematic of the $24$ spike cactus geometry and (b) the mesh around the surface of the immersed cactus on a grid corresponding to step length $h=l=0.018D$.} }
\label{cac_sch}
\end{figure} 

In this part of the flow simulation, we have considered the flow over a more complex geometry, i.e a closed curve in the shape of the cross-section of a cactus shaped cylinder. The flow configuration is similar to the flow past an impulsively started circular cylinder except the fact that the circle  is now replaced by a closed curve in the shape of the cross section of a cactus plant. Such simulations are capable of producing realistic results in some areas in the field of biology. One may cite the example of the flow past the Saguaro species of cacti, found in desert regions and which can withstand very high velocities despite its shallow root systems \cite{babu2008aerodynamic, talley2001experimental, zhdanov2019angle}. Recent studies \cite{talley2001experimental, talleymungal} have established that the cavities between two successive spikes of a cactus plant produce a dampening effect on the fluctuating drag and lift forces. To the best of our knowledge, all the earlier simulations for flows involving cactus shapes were carried out in the finite element framework.

We have used the following level set function for generating the cactus shape centered at $(x_c, y_c)$ is given by
\begin{equation}
\phi\left( r, \theta \right)  = r- r_{0}-L\sin(w \theta)
\end{equation}
where $r$=  $\sqrt{(x-x_{c})^{2}+(y-y_{c})^{2}}$, $\theta$= $arctan((y-y_{c})/(x-x_{c}))$, $r_0$, $w$ are parameters determining the base and the number of spikes and L is the maximum height of the spikes. In all our simulations, the center of the cactus shaped region is assumed to be at the origin. The geometry of the cactus along with the mesh around the surface is depicted in figures \ref{cac_sch}(a)-(b) respectively. We have chosen a  spike ratio of value $0.105$ for a $24$ spike cactus which is nothing but the ratio between the maximum height $L$ of the spike and the total diameter $D$ (set as $1$ here) of the the cylinder. Note that a spike ratio zero corresponds to a smooth cylinder described in section \ref{st_cyl}.  

\begin{figure}[!h]
 \begin{minipage}[b]{.45\linewidth}  
\includegraphics[scale=0.4]{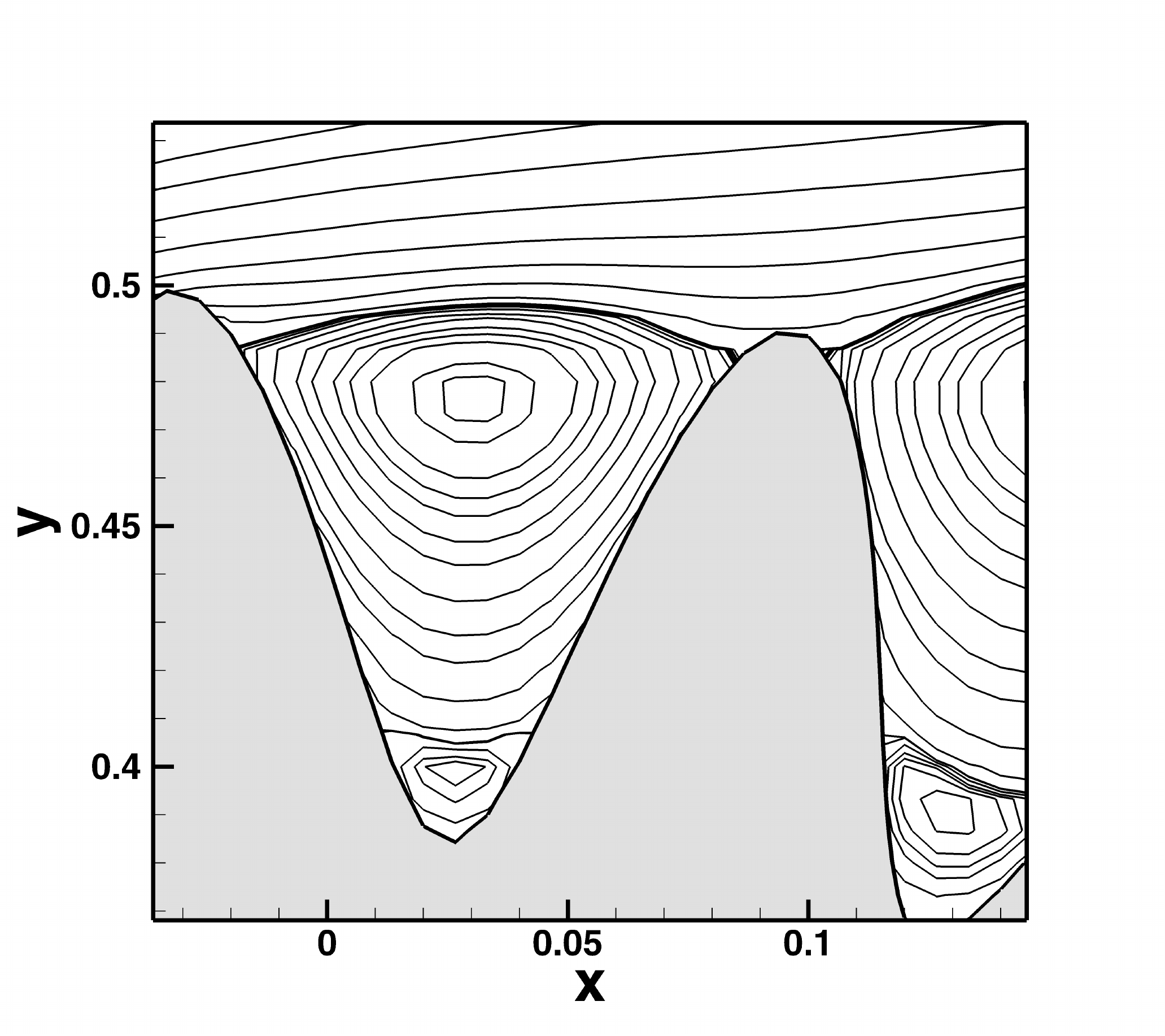} 
\centering (a) 
\end{minipage}           
\begin{minipage}[b]{.45\linewidth}
\includegraphics[scale=0.4]{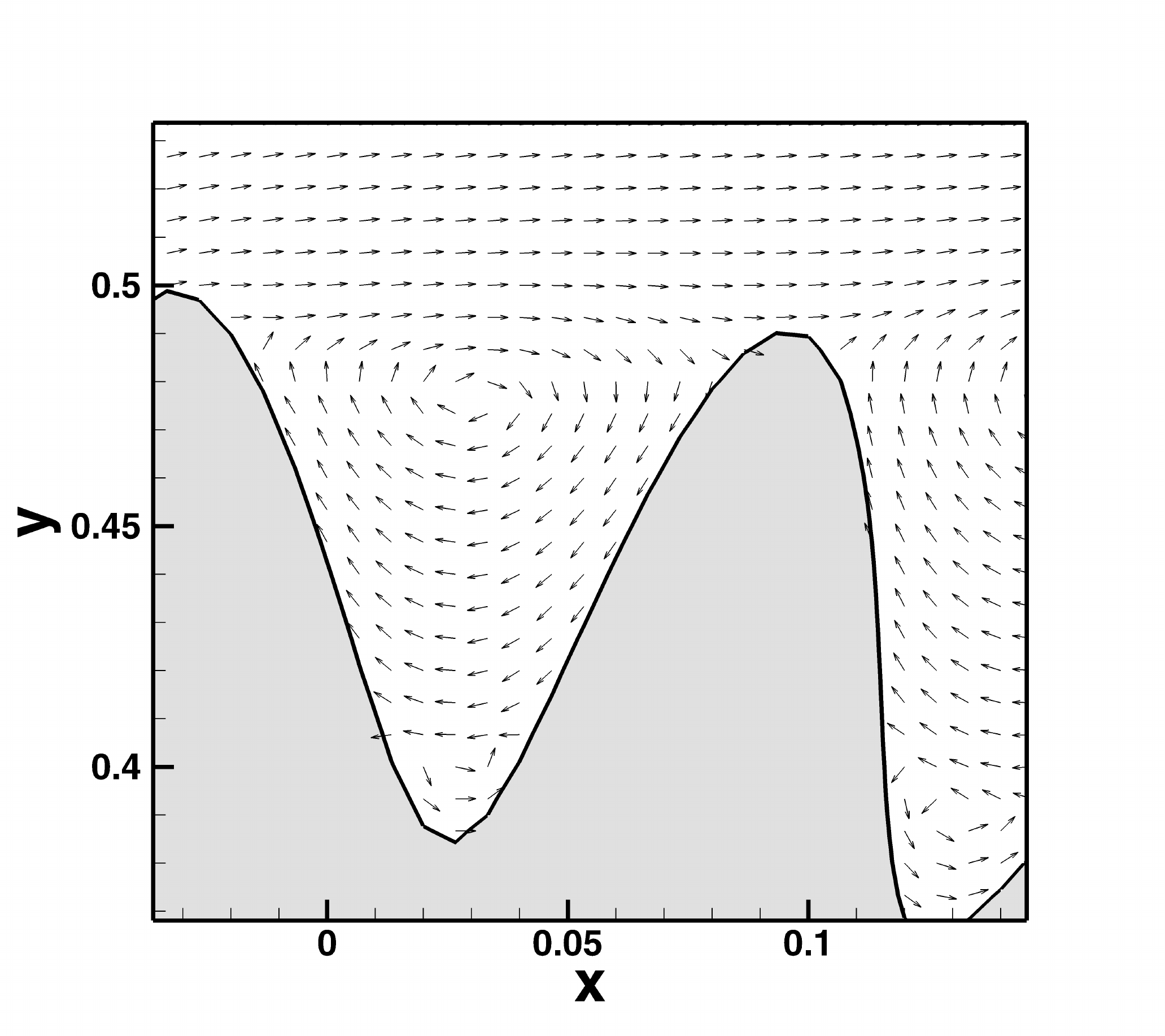}
\centering (b) 
\end{minipage} 
 \caption{{\sl Instantaneous (a) streamlines and (b) velocity vector plots for flow past a twenty four spike cactus cylinder for  $Re=300$.} }
\label{vec_sf_300}
\end{figure}

\begin{figure}[!ht]
\begin{minipage}[b]{.45\linewidth}   \
\hfill \centering\psfig{file=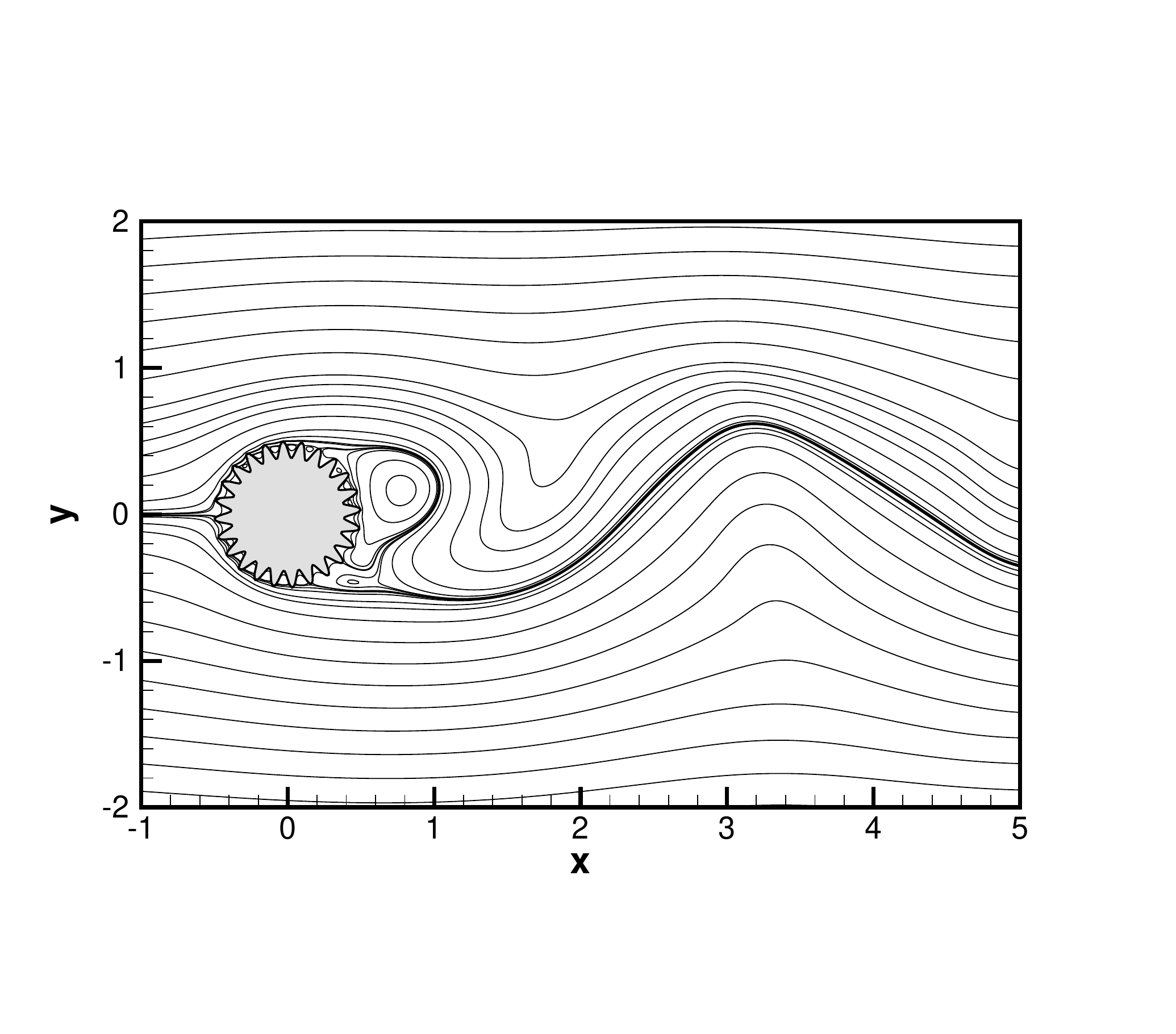,width=.7\linewidth}(a)
\end{minipage}
\begin{minipage}[b]{.45\linewidth}
\centering\psfig{file=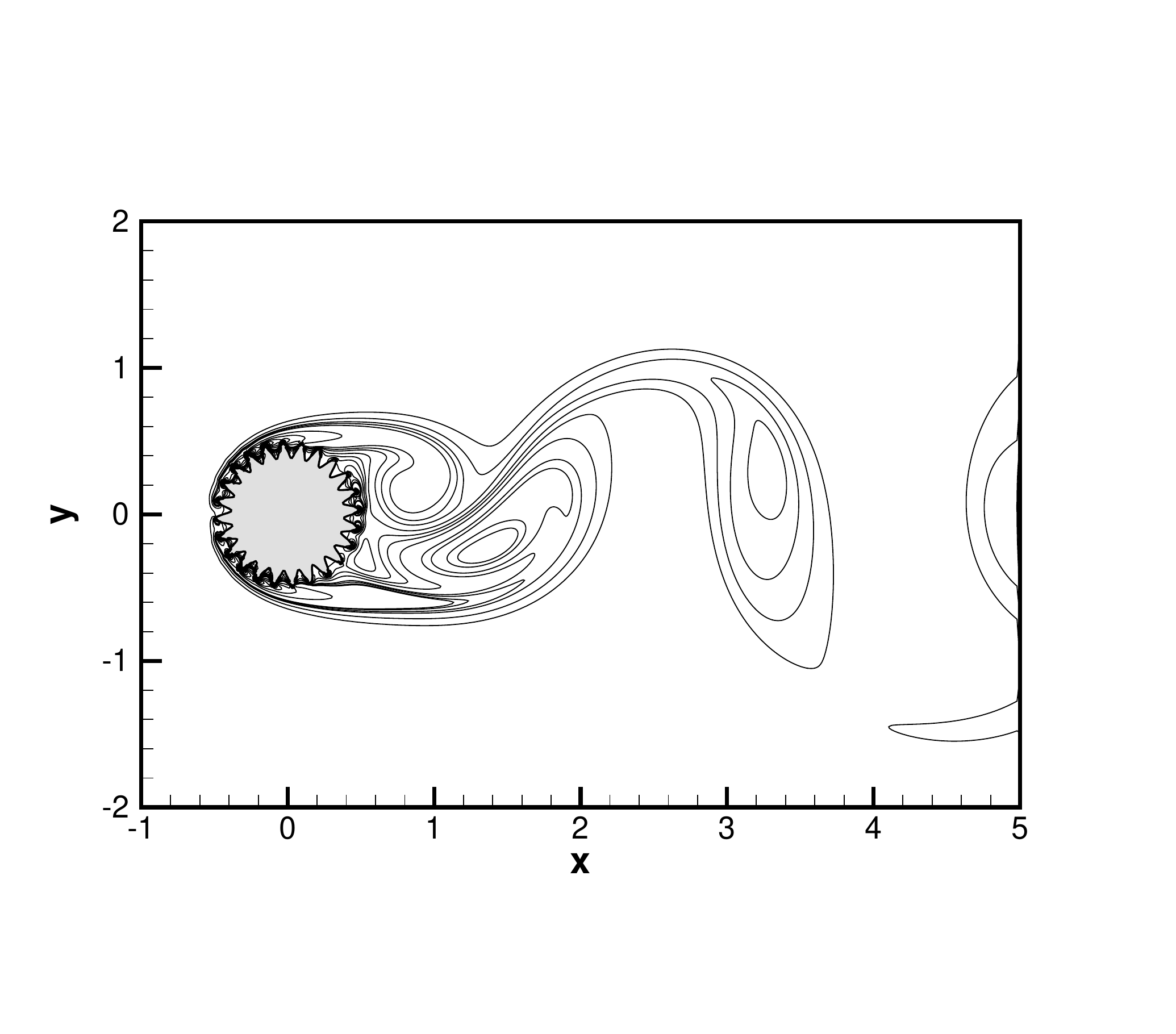,width=.7\linewidth}(a)
\end{minipage}
\begin{minipage}[b]{.45\linewidth}   \
\hfill \centering\psfig{file=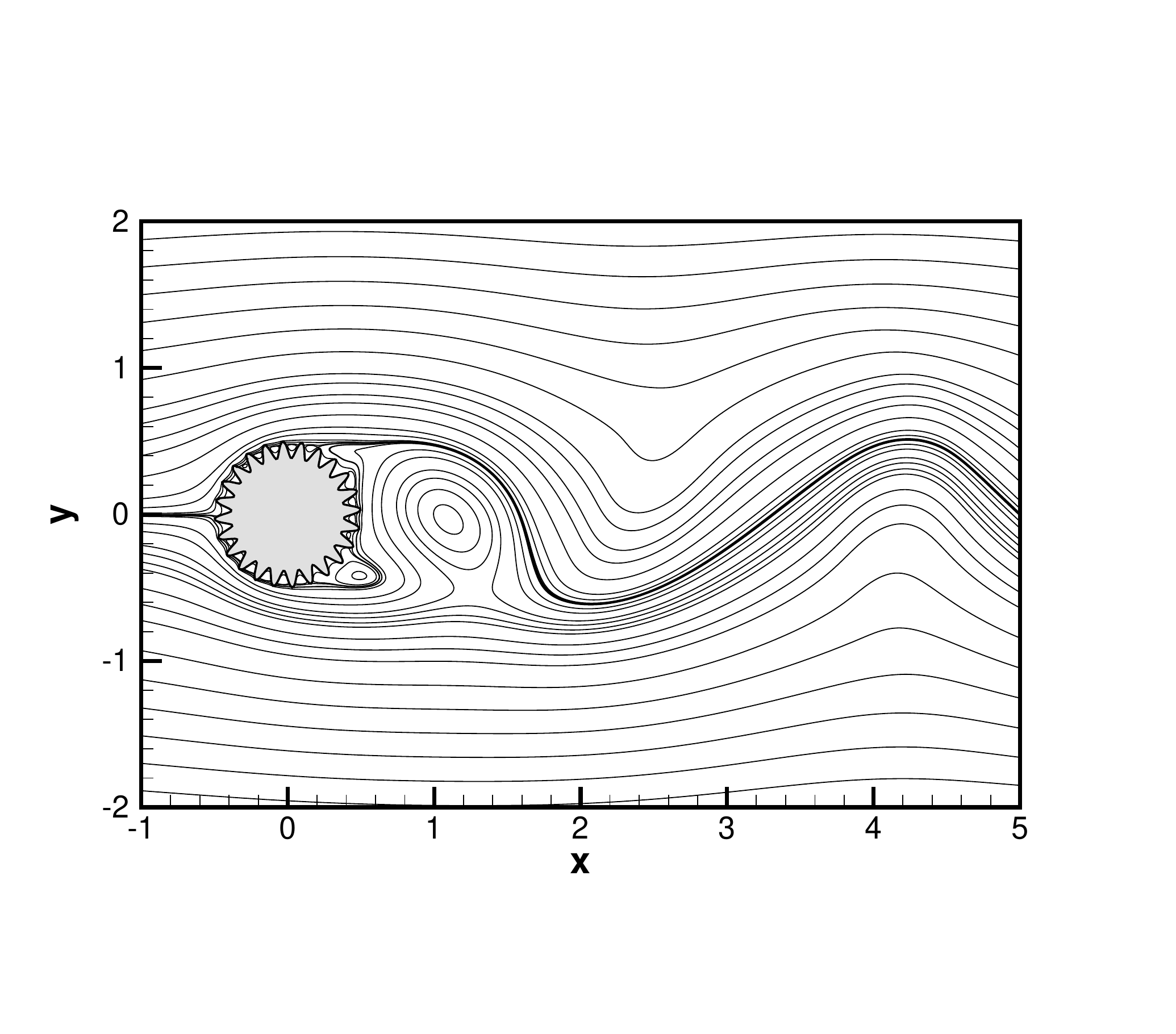,width=.7\linewidth}
 (b)
\end{minipage}
\begin{minipage}[b]{.45\linewidth}
\centering\psfig{file=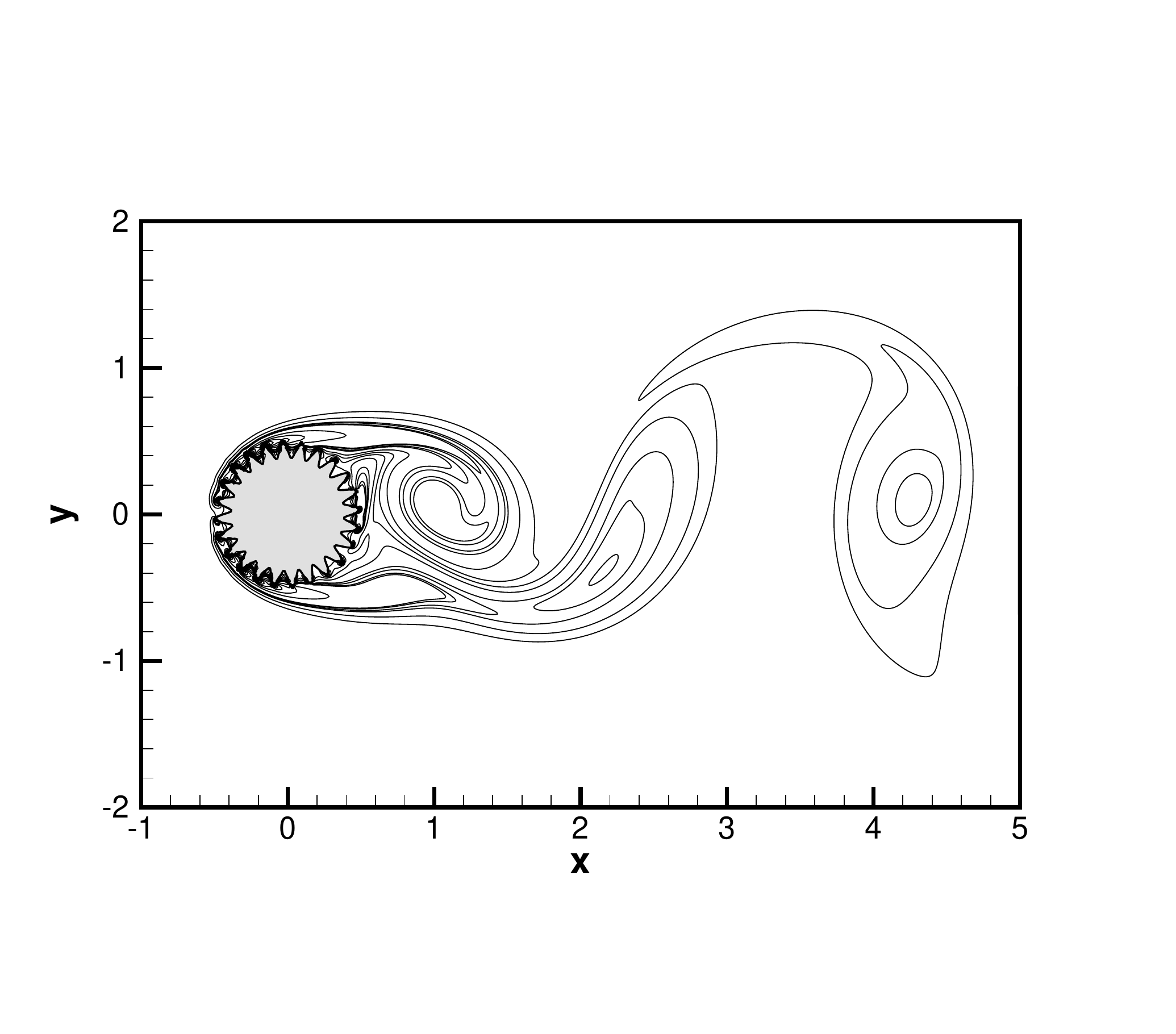,width=.7\linewidth}
 (b)
\end{minipage}
\begin{minipage}[b]{.45\linewidth}   \
\hfill \centering\psfig{file=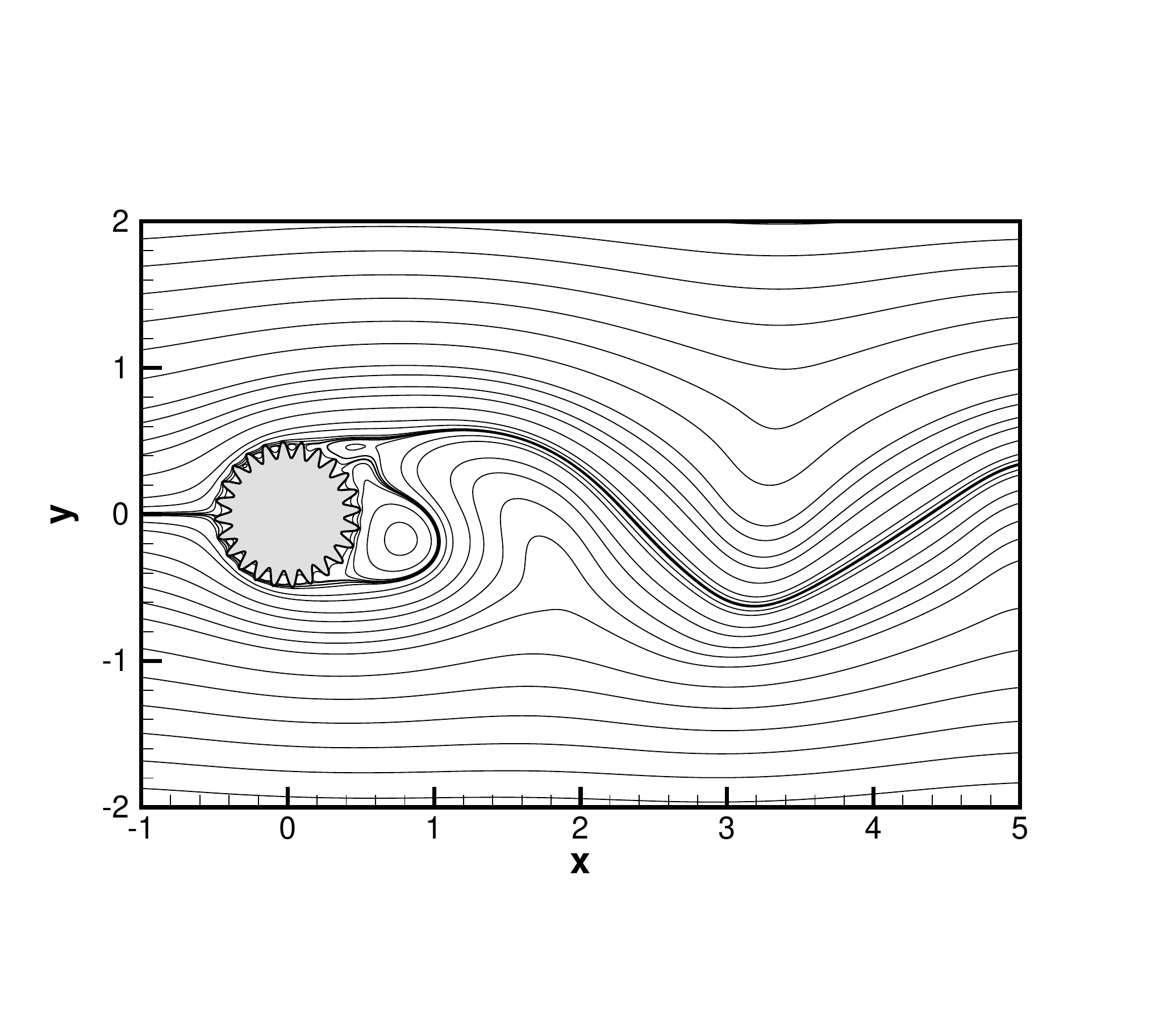,width=.7\linewidth}
 (c)
\end{minipage}
\begin{minipage}[b]{.45\linewidth}
\centering\psfig{file=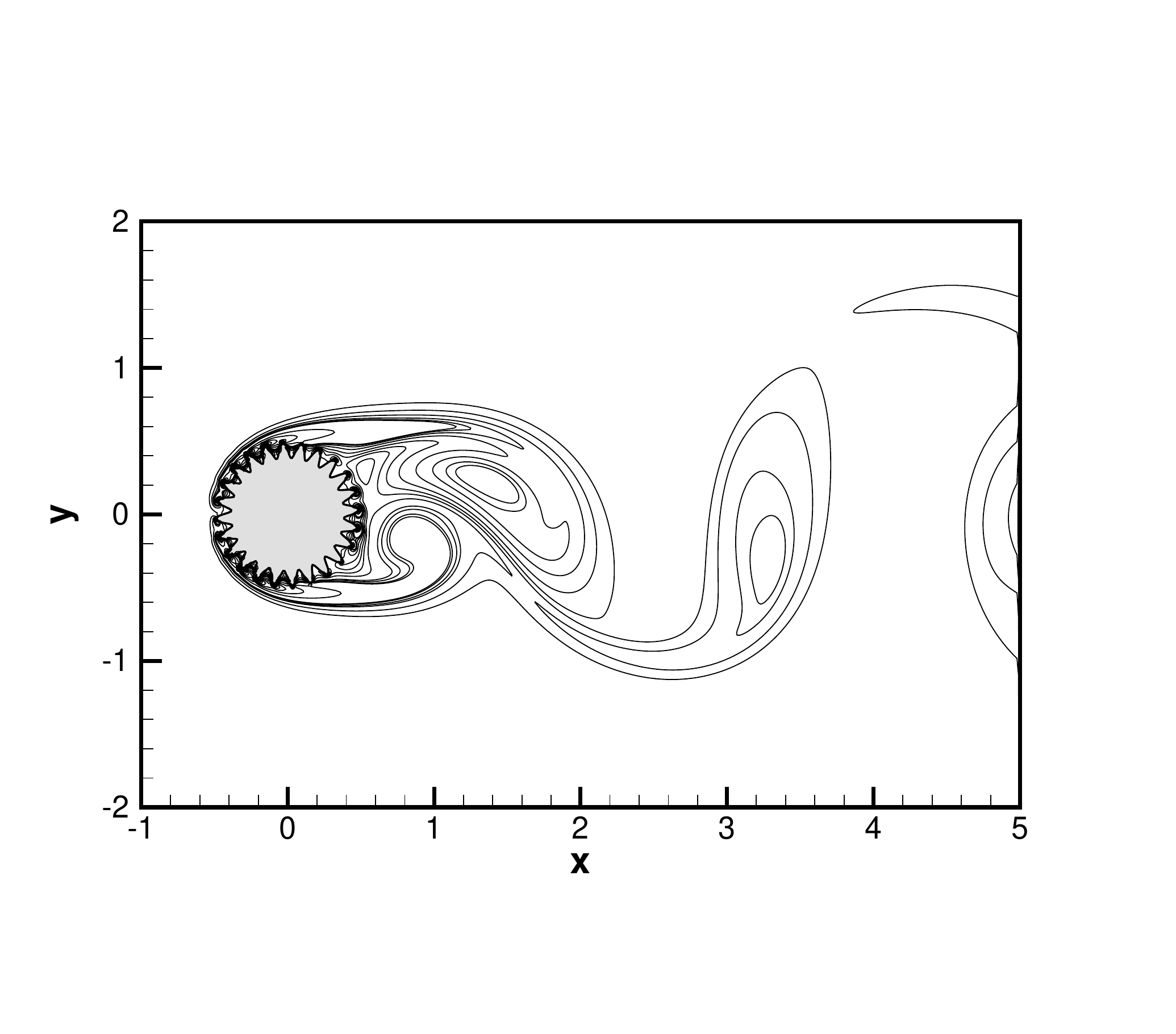,width=.7\linewidth}
 (c)
\end{minipage}
\begin{minipage}[b]{.45\linewidth}   \
\hfill \centering\psfig{file=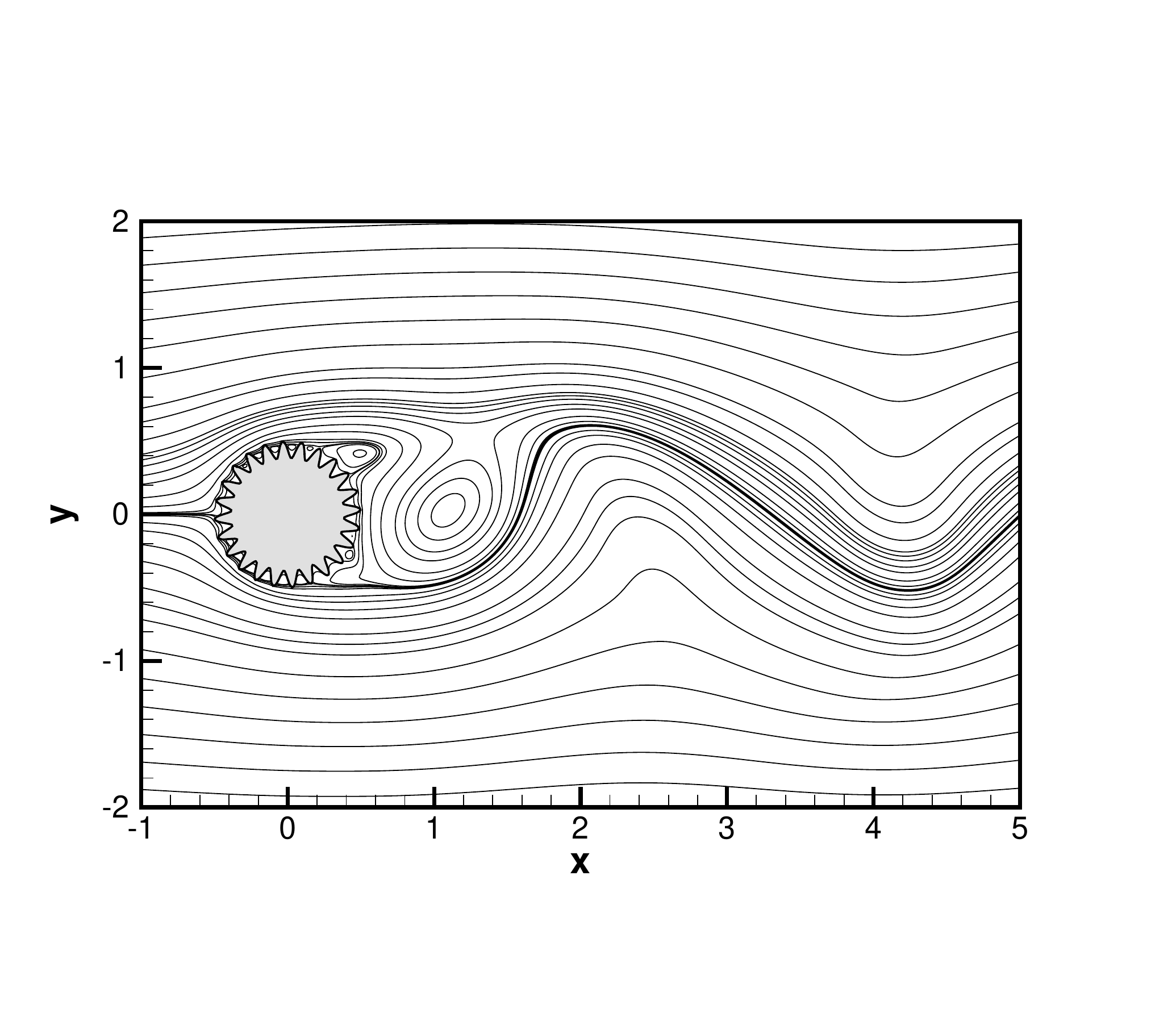,width=.7\linewidth}
 (d)
\end{minipage}
\begin{minipage}[b]{.45\linewidth}   \
\centering\psfig{file=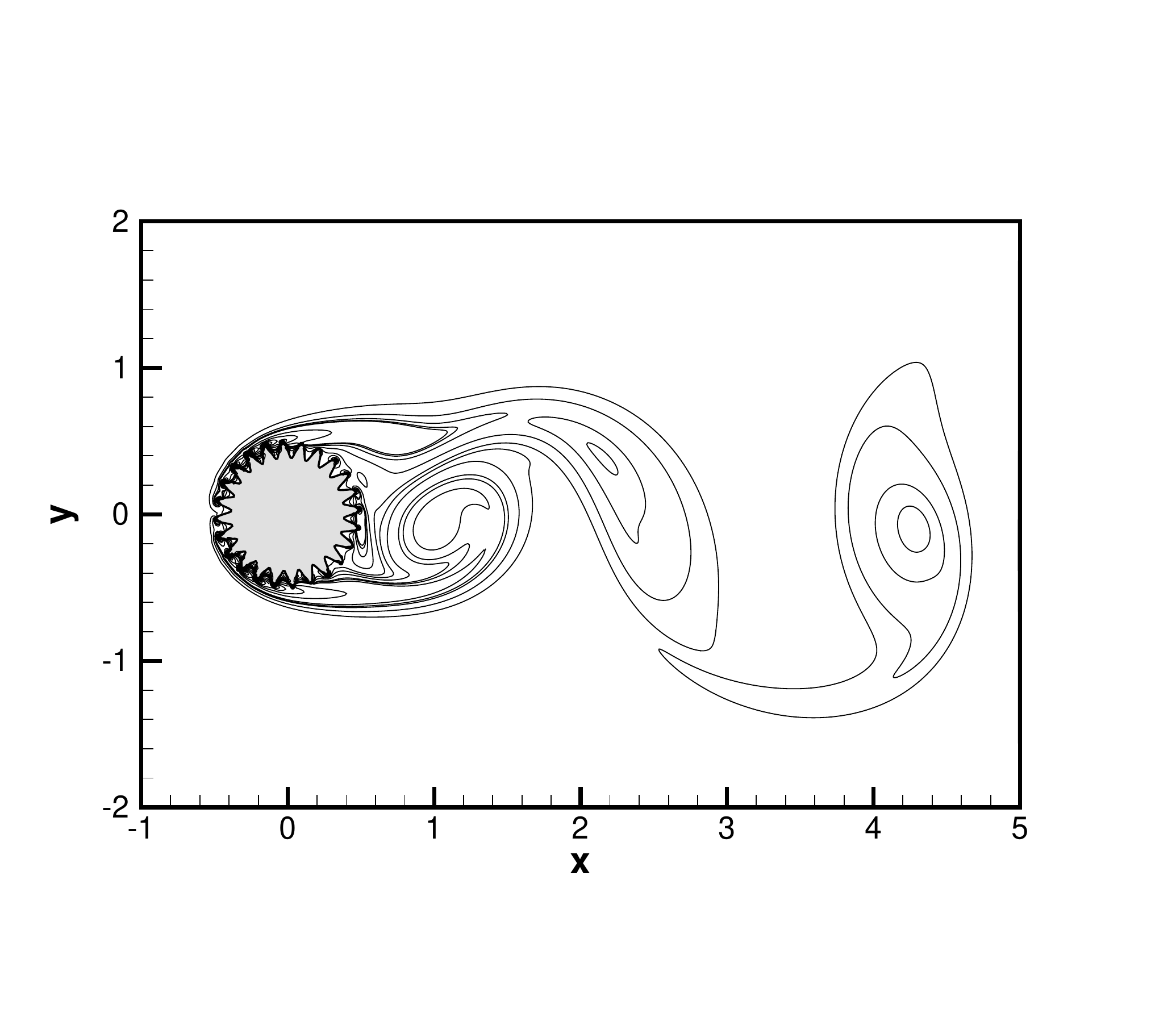,width=.7\linewidth}
 (d)
\end{minipage}
\begin{minipage}[b]{.45\linewidth}
\hfill \centering\psfig{file=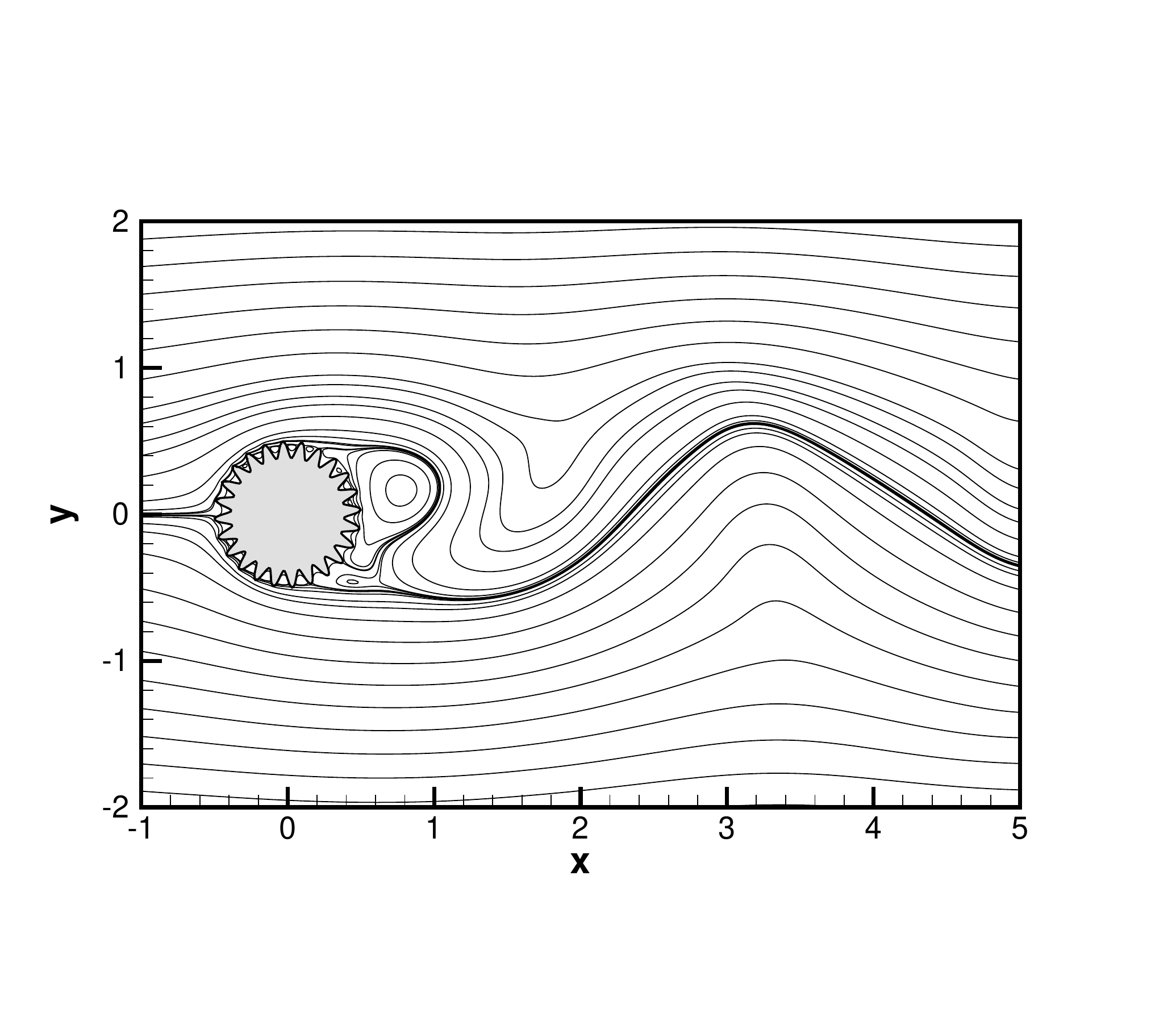,width=.65\linewidth}
 (e)
\end{minipage}
\begin{minipage}[b]{.45\linewidth}   \
\centering\psfig{file=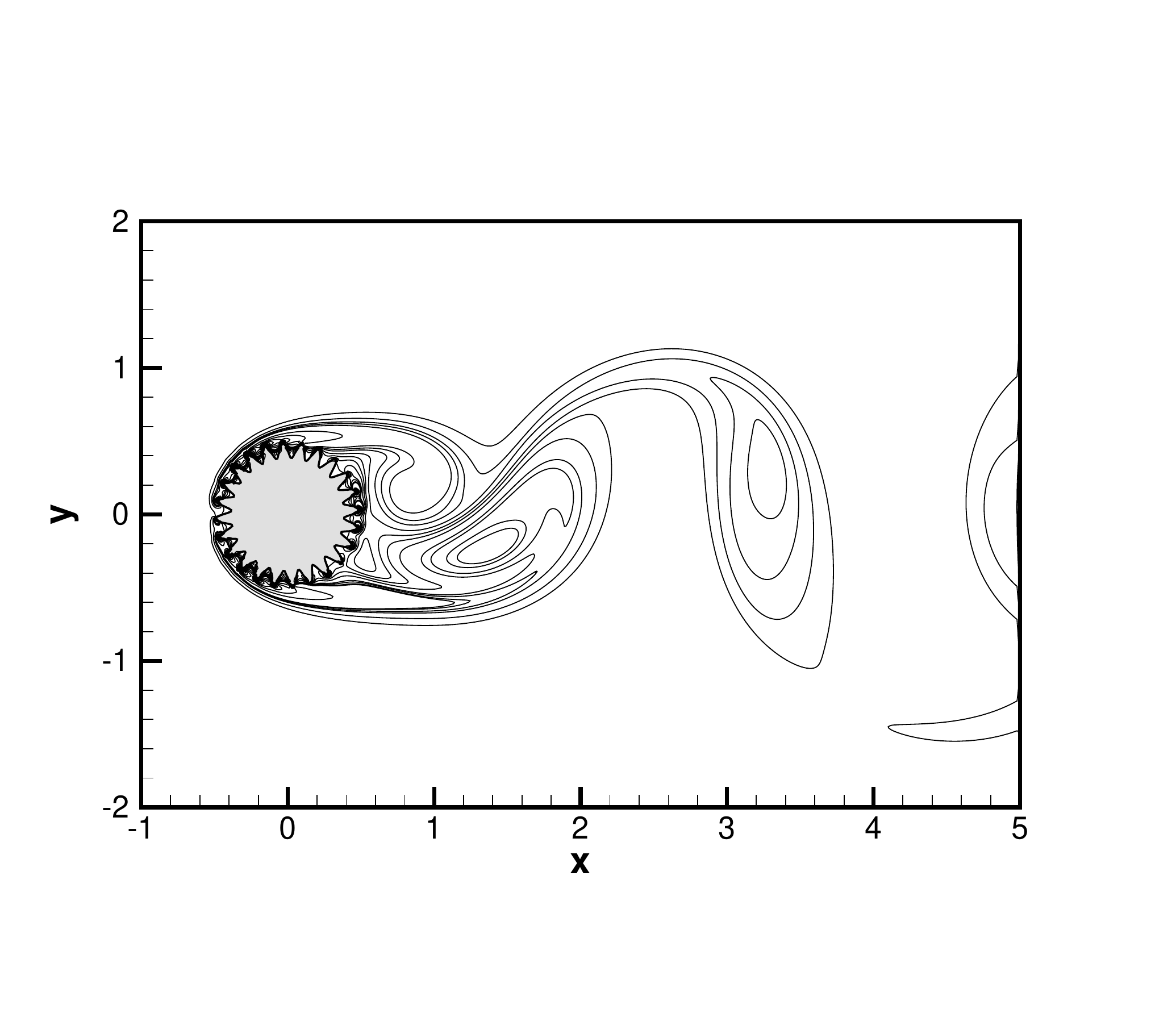,width=.65\linewidth}
 (e)
\end{minipage}
\caption{\sl {Streamfunction (left) and vorticity (right) contours for the flow past a $24$ spike cactus for $Re=300$: : (a) $t = 0$, (b) $t = \frac{\pi}{2}$, (c) $t =\pi$, (d) $t =  \frac{3\pi}{2}$ and (e) $t =2\pi$.}}\label{cac_300_sf_vt}
\end{figure}

\begin{figure}[!h]
\begin{center} 
\includegraphics[scale=0.5]{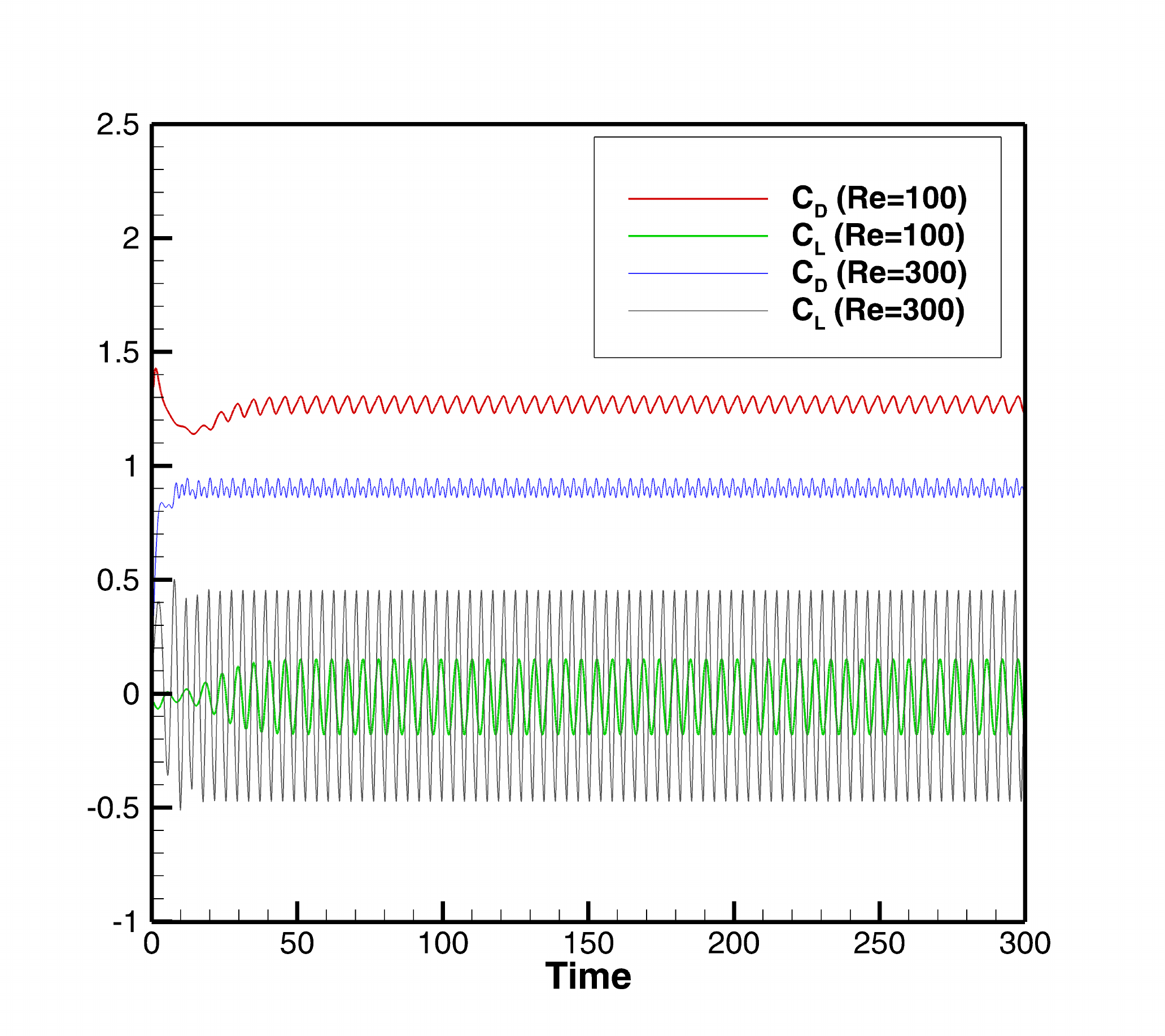}  
\end{center}
\caption{{\sl Time histories of drag and lift coefficients for the cactus shaped cylinder for $Re=100$ and $300$.} }
\label{his_cac}
\end{figure}

Computations were carried out for $Re=100$ and $300$ along with that for the smooth cylinder till periodic vortex shedding stage is reached. Opposed to the smooth cylinder, the flow pattern in the neighbourhood of the surface changes significantly for cactus shaped cylinder. Figure  \ref{vec_sf_300}(a)-(b) shows the instantaneous streamlines and velocity vector plots respectively inside a cactus groove for $Re=300$. These plots clearly indicate the presence of recirculation zones inside the grooves; one can also see the existence of a secondary zone which is consistent with the findings of \cite{babu2008aerodynamic}. In order to gain further insight into the flow field variation around the cactus cylinders, we plot the streamfunction (left column) and vorticity contours (right column) in figures \ref{cac_300_sf_vt}(a)-(e) at five different phases within a shedding cycle. Note that figure \ref{cac_300_sf_vt}(c) is a mirror image of \ref{cac_300_sf_vt}(a) and \ref{cac_300_sf_vt}(e) while figure \ref{cac_300_sf_vt}(d) is a mirror image of figure \ref{cac_300_sf_vt}(b). This is because of the fact the shedding of vortices takes place from the upper and the lower parts of the cylinder alternately during a half-cycle period, thus exemplifying the efficiency of our approach in accurately capturing the phenomenon. The time history of drag and lift coefficients for these two $Re$s, depicting the periodic nature of the flow is shown in figure \ref{his_cac}. 

\begin{table}[]
\centering
\caption{ Comparison of Strouhal number, drag and lift coefficients of the periodic flow for Cactus shaped and Circular Cylinders for $Re =100$ and $300$.}
\begin{tabular}{|c|ccc|ccc|}
\hline
$Re$                               & \multicolumn{3}{c|}{$100$}                & \multicolumn{3}{c|}{$300$}                               \\ \hline
Flow parameters &
$Cactus$ &
\multicolumn{1}{c}{$Circular$} &
$ \% difference$ &
\multicolumn{1}{c}{$Cactus$} &
\multicolumn{1}{c}{$Circular$} &
$\% difference$ \\ \hline
$St$& $0.175$       & $0.180$                 & -           & $0.220$ & $0.251$            & $-$    \\
$C_D$                         & $1.268$ & $1.402$                 & $9.56$           & $0.905$ & $1.205$                 & $24.89$                  \\
$C_L$                           & $\pm0.151$ & $\pm 0.233$            &$35.19$           & $\pm 0.454$ & $\pm 0.596$                 &$23.82$                  \\
 \hline
\end{tabular}
\label{dl_comp}
\end{table}
\begin{figure}[!t]
 \begin{minipage}[b]{.45\linewidth}  
\includegraphics[scale=0.4]{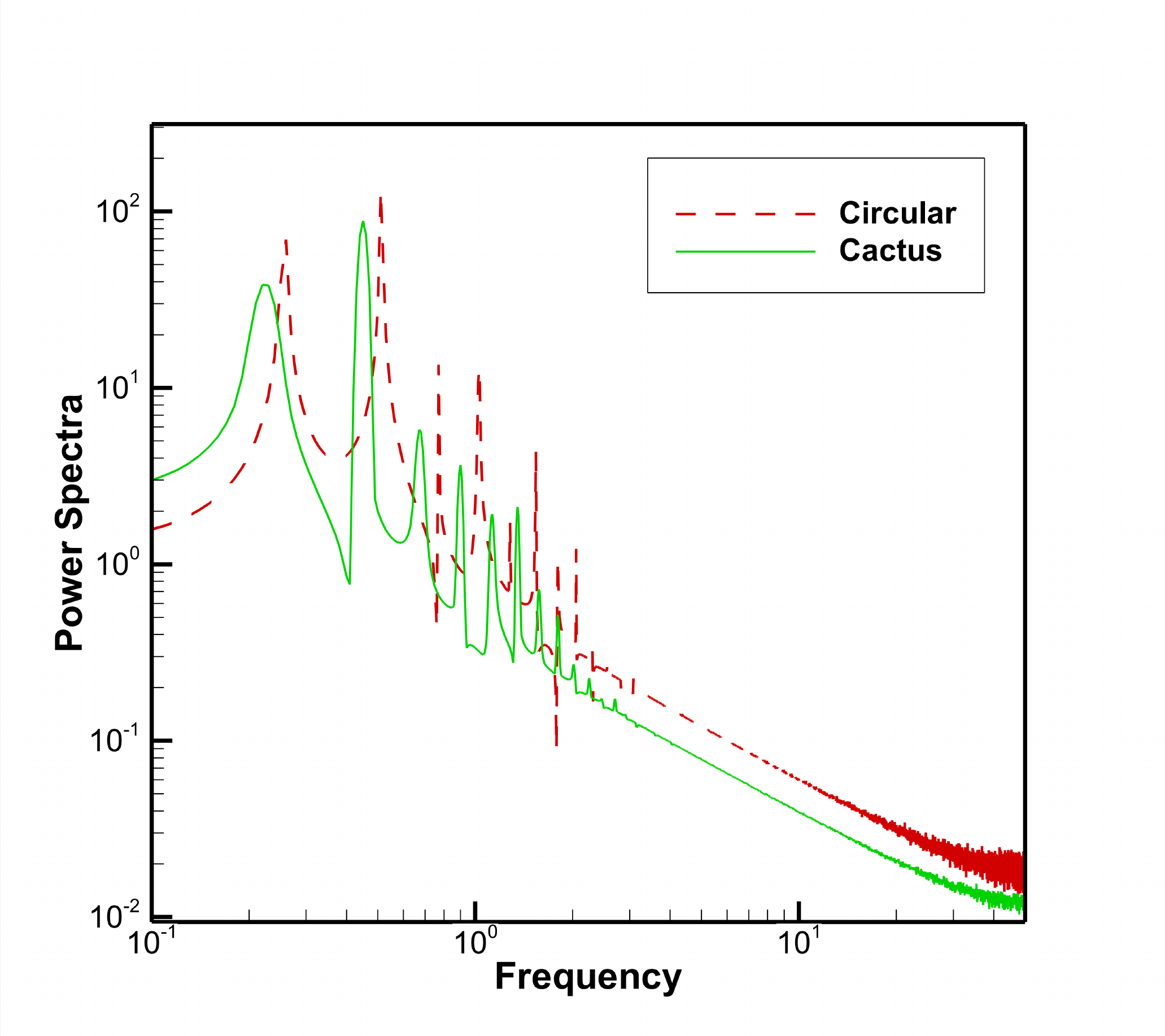} 
\centering (a) 
\end{minipage}           
\begin{minipage}[b]{.45\linewidth}
\includegraphics[scale=0.4]{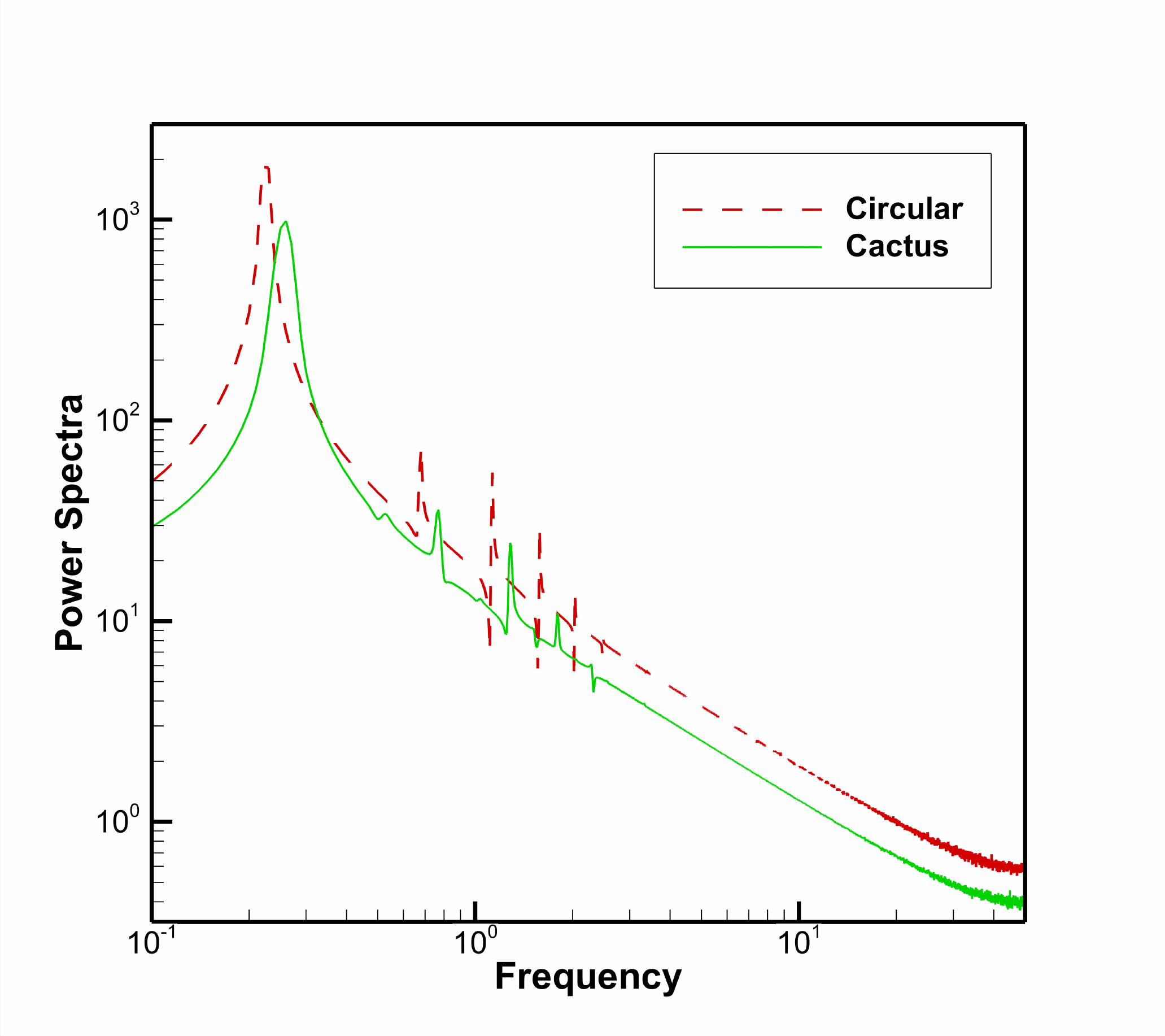}
\centering (b) 
\end{minipage} 
 \caption{{\sl Comparison of power spectra of the cactus shaped  and circular (smooth) cylinders for  $Re=300$ based on the time histories of (a)drag  and (b) lift coefficients.} }
\label{spectra_cactus_dl_300}
\end{figure}
In table \ref{dl_comp}, we compare the force coefficients and Strouhal numbers for the periodic flows for the cactus shaped and circular cylinders for $Re=100$ and $300$. One can clearly see a reduction in the unsteady loads for the cactus shaped cylinder compared to the smooth circular cylinder. With increase in Reynolds number, the drag reduction becomes more prominent. Our observations are consistent with the findings of  Babu and Mahesh \cite{babu2008aerodynamic}. The decrement in percentage difference from reference \cite{babu2008aerodynamic} may be attributed to the fact that while their valleys were extremely sharp,  the ones employed in our computation are smooth (see figure \ref{cac_sch}(a)). The reduction in the unsteady loads is also evident from the comparison of the power spectra of the cactus shaped and circular cylinder shown in figure \ref{spectra_cactus_dl_300} which reflects the decrease in the higher frequency contents of the loads. It is worth mentioning that while \cite{babu2008aerodynamic} had used $436 \times 10^3$ and $7 \times 10^6$ hexahedral elements for flow computations for $Re=100$ and $300$ respectively, we have used only $110 \times 10^3$ nodes for both the cases. Moreover, while the elements of the grid converged solution were $0.0052D$ in the azimuthal direction and $0.0006D$ radially on the surface of the cylinder in their case, we used a step length $0.018D$ in both the vertical and horizontal directions in our computation (see figure \ref{cac_sch}(b)).  
\clearpage

\subsection{Flow Past two tandem circular cylinders}
Flow induced oscillations of multiple circular cylinders is an extremely complex flow problem. It has garnered immense interest over the last few decades owing to its real life applications in the field of off-shore oil drilling rigs and tall chimneys, heat exchanger and riser tubes, cooling of nuclear fuel rods amongst others. In the same vein, the next problems considered here are the flow past two tandem cylinders of equal diameter $D$ separated by a distance $PD$, where the downstream cylinder is stationary and the upstream cylinder is either kept stationary or oscillating transversely. 

\subsubsection{Test Case 4: Stationary upstream cylinder}
\begin{figure}[!h]
\begin{minipage}[b]{.5\linewidth}  
\includegraphics[scale=0.43]{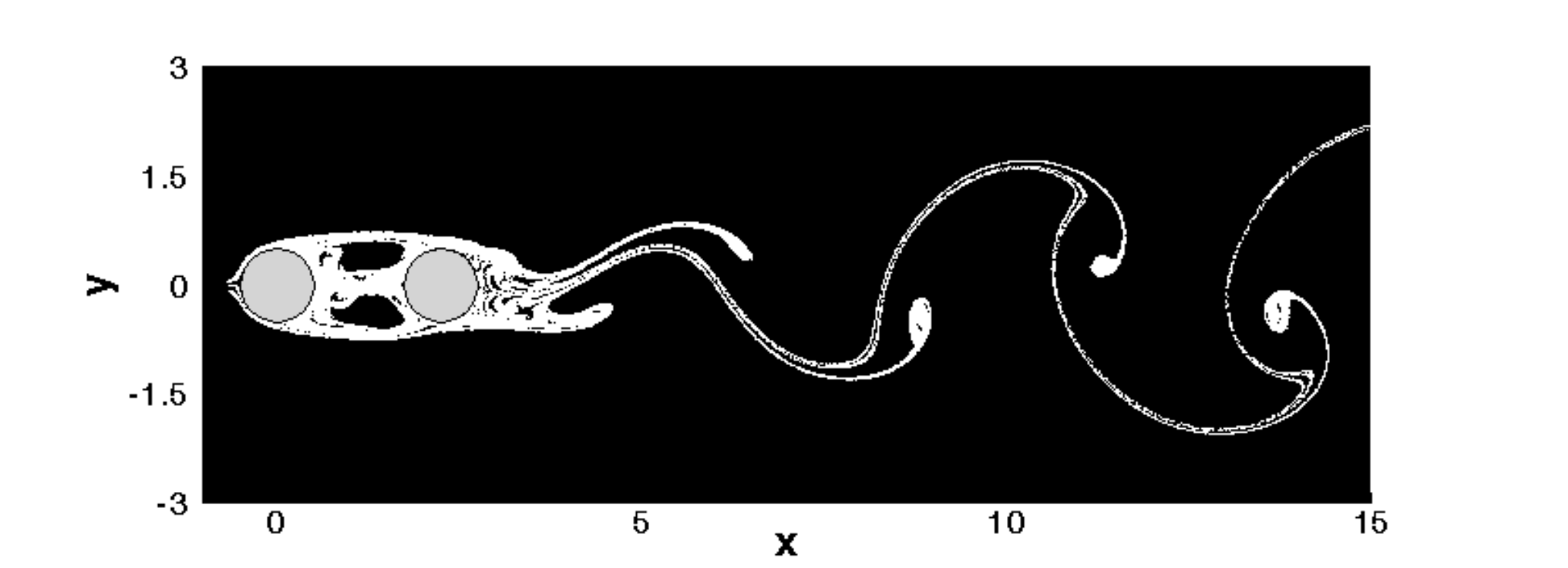} 
 \centering (a)
\end{minipage}            \hspace{-2.mm}
\begin{minipage}[b]{.5\linewidth}
\includegraphics[scale=0.43]{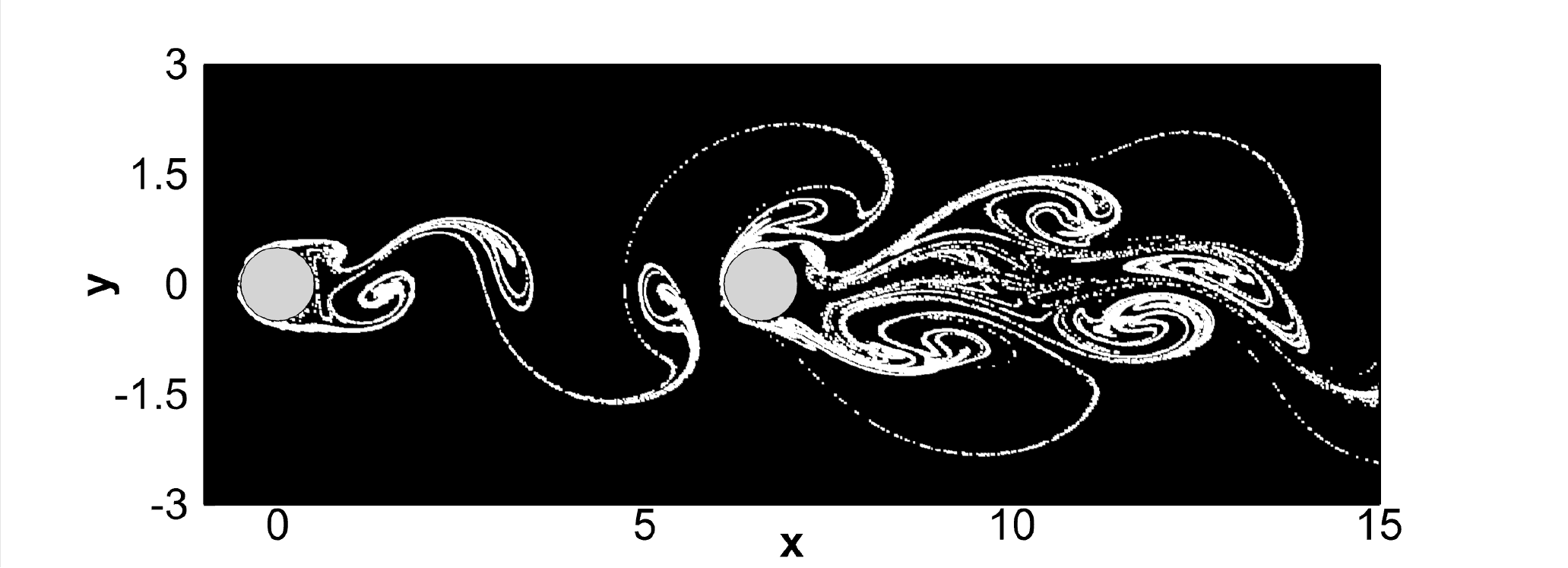} 
 \centering (b)
\end{minipage} 
\begin{minipage}[b]{.5\linewidth}  
\includegraphics[scale=0.43]{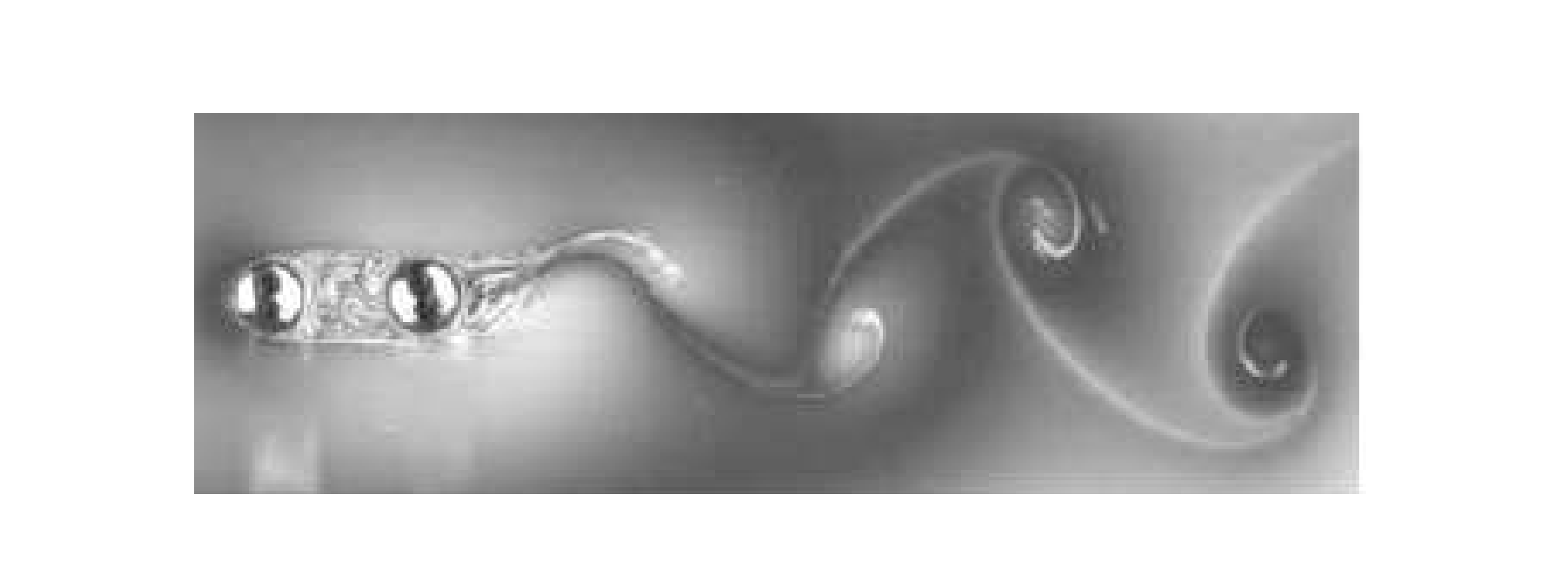} 
 \centering (c)
\end{minipage}            \hspace{-2.mm}
\begin{minipage}[b]{.5\linewidth}
\includegraphics[scale=0.43]{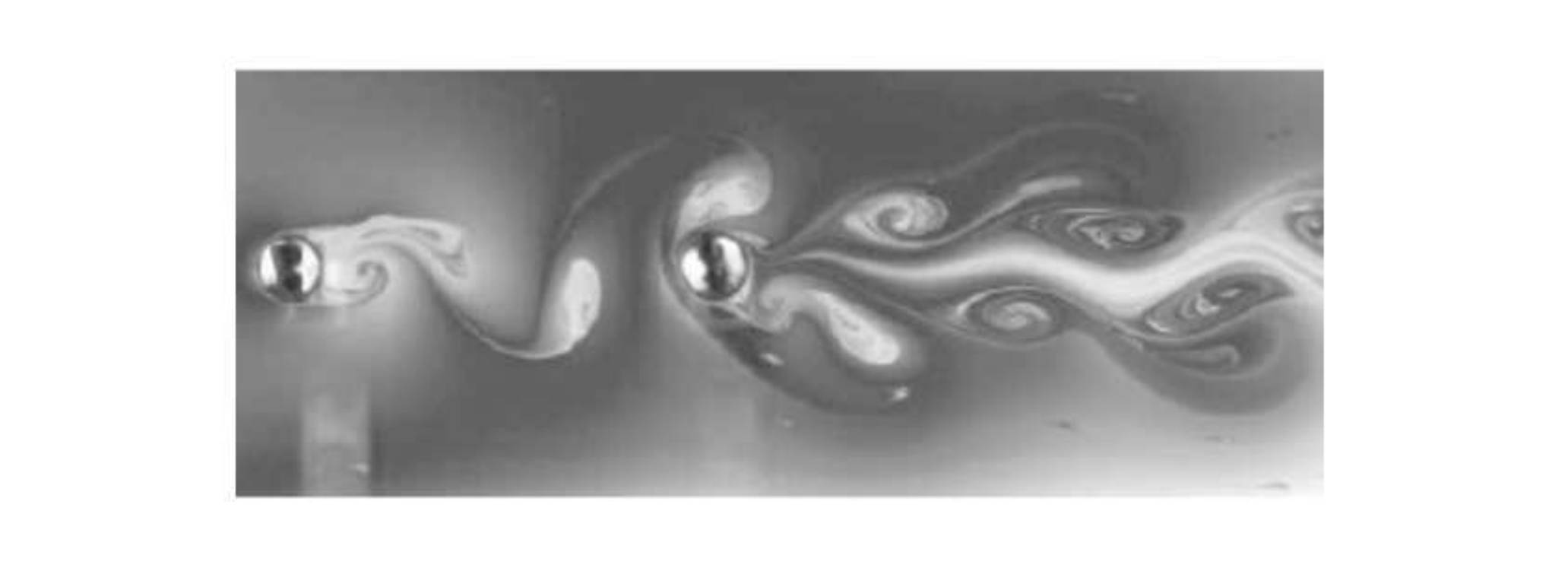} 
 \centering (d)
\end{minipage} 
\caption{{\sl  Comparison of (a)-(b) our computed streaklines for $P=2.25$ and $6.58$ respectively with the (c)-(d) experimental visualization of Yang {\it et al.} \cite{yang2019critical} for Test Case 4. } }
\label{streak}
\end{figure}

Here the flow configuration is similar to the one shown in figure \ref{setup_bluff} except the fact that the single bluff body in the figure is replaced by a stationary circular cylinder at the origin and another cylinder of same dimension is placed to its right at a distance $PD$ apart for $Re=200$. We have considered $P=2.25$ and $6.58$ which characterises the medium pitch and long-pitch regime categorized recently by Hoisseini {\it et al} \cite{hosseini2020flow}. The gap between the cylinders  are chosen so as to compare our simulations with the recent experimental visualizations of Yang {\it et al.} \cite{yang2019critical} and in the process validate our simulations. The streaklines resulting from our computations are presented top and bottom along with the visualizations of Yang {\it et al.} \cite{yang2019critical} in figure \ref{streak}. Our simulations are extremely close to experimental ones exemplifying the efficiency of our immersed interface approach. From figure \ref{streak}(a),(c), one can spot that no shedded vortex is visible in the gap between the cylinders in the medium pitched regime and wake behind the downstream cylinder resembles that of an isolated cylinder. On the other hand, for the long pitch regime, vortex shedding reminiscent of an isolated cylinder is observed in the gap as well as behind the downstream cylinder (see figure \ref{streak}(b),(d)). Our observations are consistent with the experimental results of  \cite{yang2019critical} and the numerical simulations of \cite{hosseini2020flow}, thus establishing the robustness of the current approach.

\subsubsection{Test Case 5: Oscillating upstream cylinder}
\begin{figure}[!h]
 \begin{minipage}[b]{.61\linewidth}  
\includegraphics[scale=0.525]{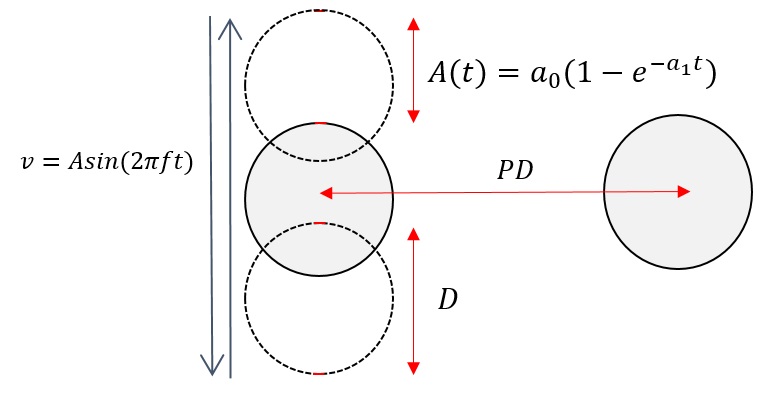} 
\centering (a) 
\end{minipage}           
\begin{minipage}[b]{.37\linewidth}
\includegraphics[scale=0.35]{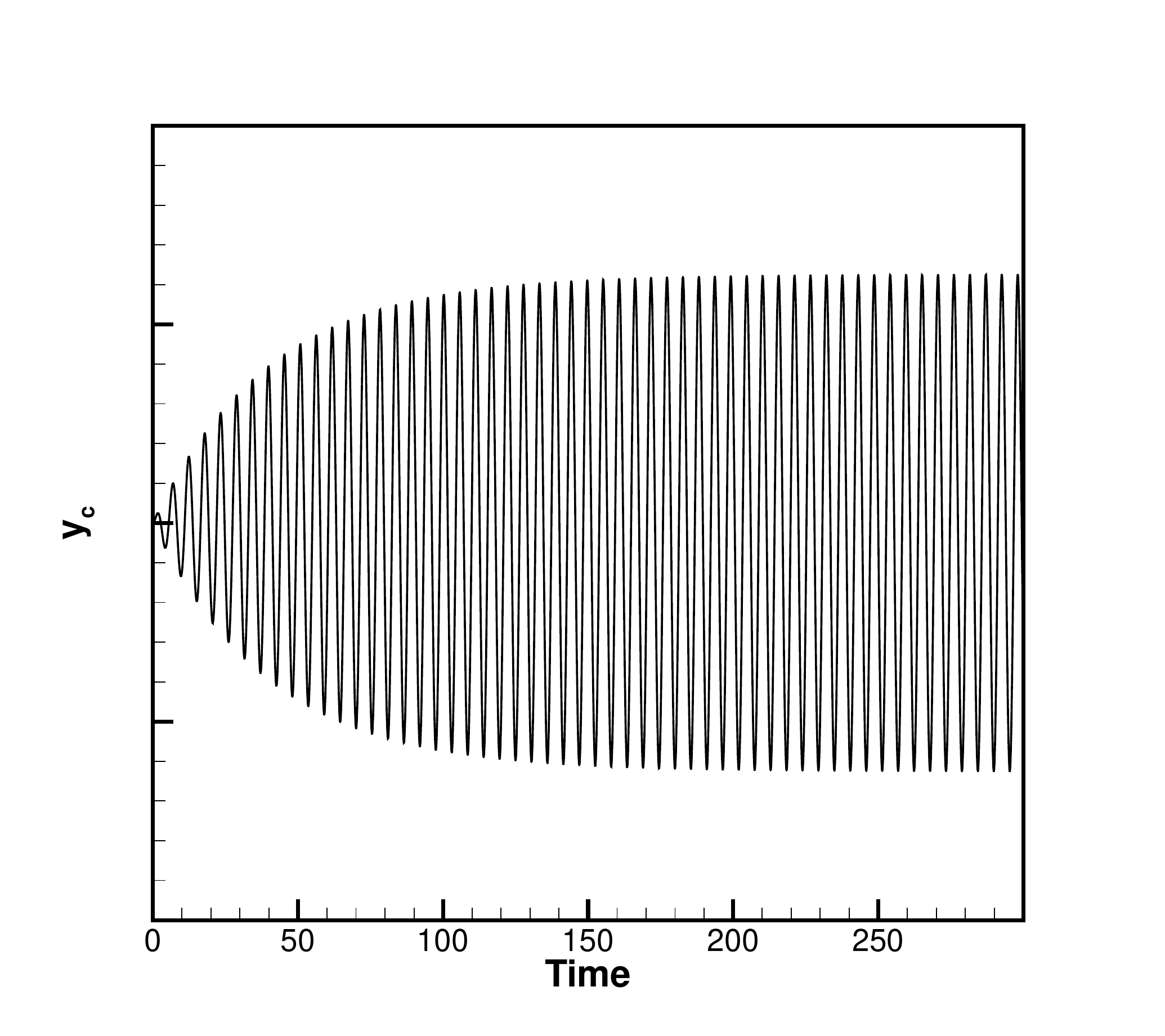}
\centering (b) 
\end{minipage} 
 \caption{{\sl (a) Schematic of the upstream and downstream cylinders and (b) Displacement of the upstream cylinder center for Test Case 5.} }
\label{tandem}
\end{figure}

Next, we consider the case, where the stationary upstream cylinder in the above example is now replaced by a transversely oscillating one. This study is similar to the third experimental arrangement by Kim {\it et al.} \cite{kim2009flow} where a fixed cylinder was placed in the wake of a transversely oscillating cylinder to suppress vortex induced vibration (VIV). However, in our computation, we allow the vibrating amplitude of the upstream cylinder to vary and evolve with time $t$ through the function $\displaystyle A(t)=a_0(1-e^{-a_1t})$; for  $a_0=0.0$, it reduces to a stationary cylinder. Under the imposed oscillation, the displacement of the $y$-coordinate of the center of this cylinder is given by $\displaystyle A(t)\sin(2\pi f t)$. In our computations, we have chosen $a_0=0.626,\; a_1=0.025,\;P=2.95\; {\rm and }\; f=0.182$. The schematic of the problem is similar to figure \ref{setup_bluff} except the fact that the bluff body shown in that figure is now replaced by the upstream oscillating cylinder with velocity $\displaystyle (u,v)=(0,2A(t)\pi f\cos(2\pi f t))$ and another stationary cylinder of the same dimension is placed to its right. The schematic of these two cylinders and the time history of the displacement of the upstream cylinder is shown in figures \ref{tandem}(a)-(b) respectively.

\begin{figure}[!ht]

\begin{minipage}[b]{.45\linewidth}   \
\hfill \centering\psfig{file=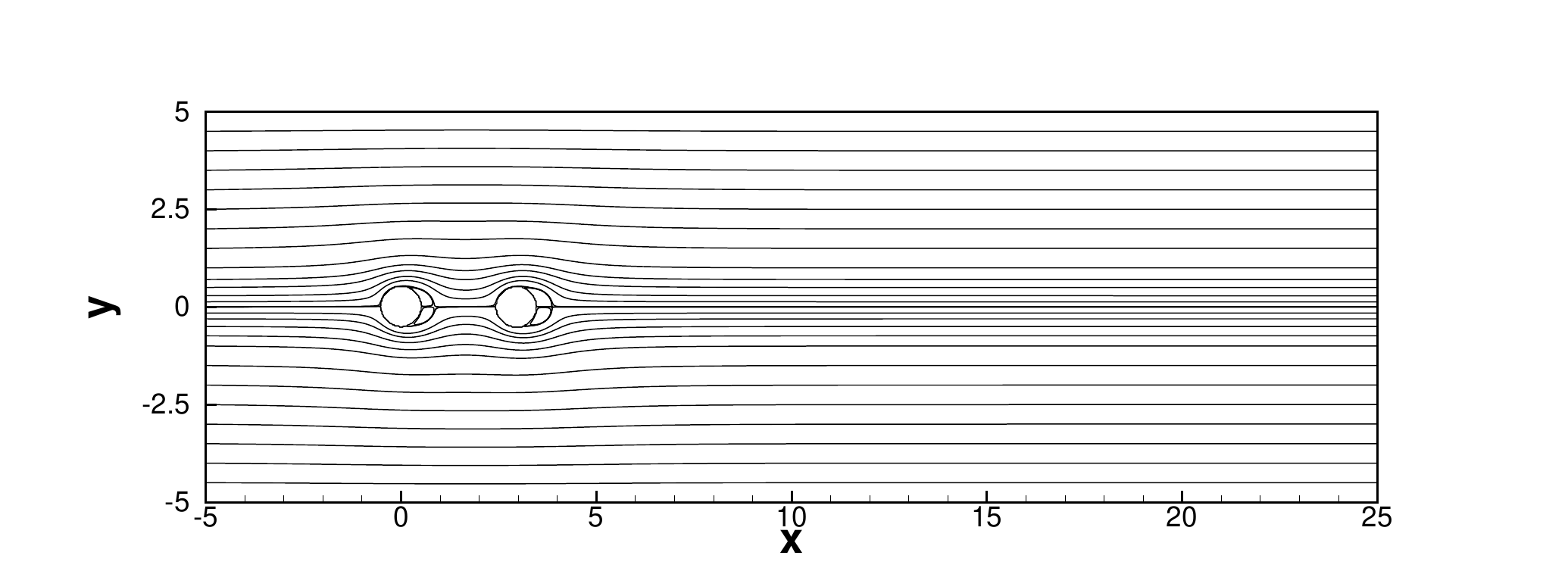,width=1.0\linewidth}
\end{minipage} 
\hspace{-0.75cm}
\begin{minipage}[b]{.45\linewidth}
\centering\psfig{file=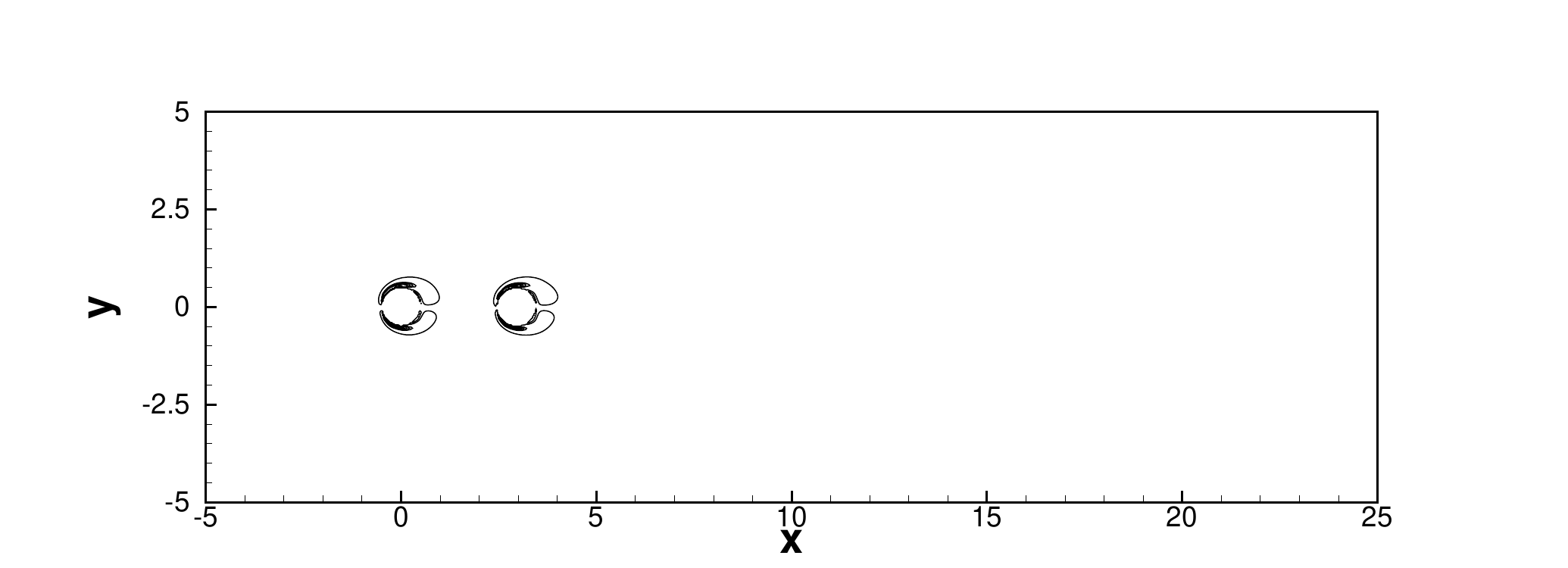,width=1.0\linewidth}
\end{minipage}(a)
 \begin{minipage}[b]{.45\linewidth}   \
\hfill \centering\psfig{file=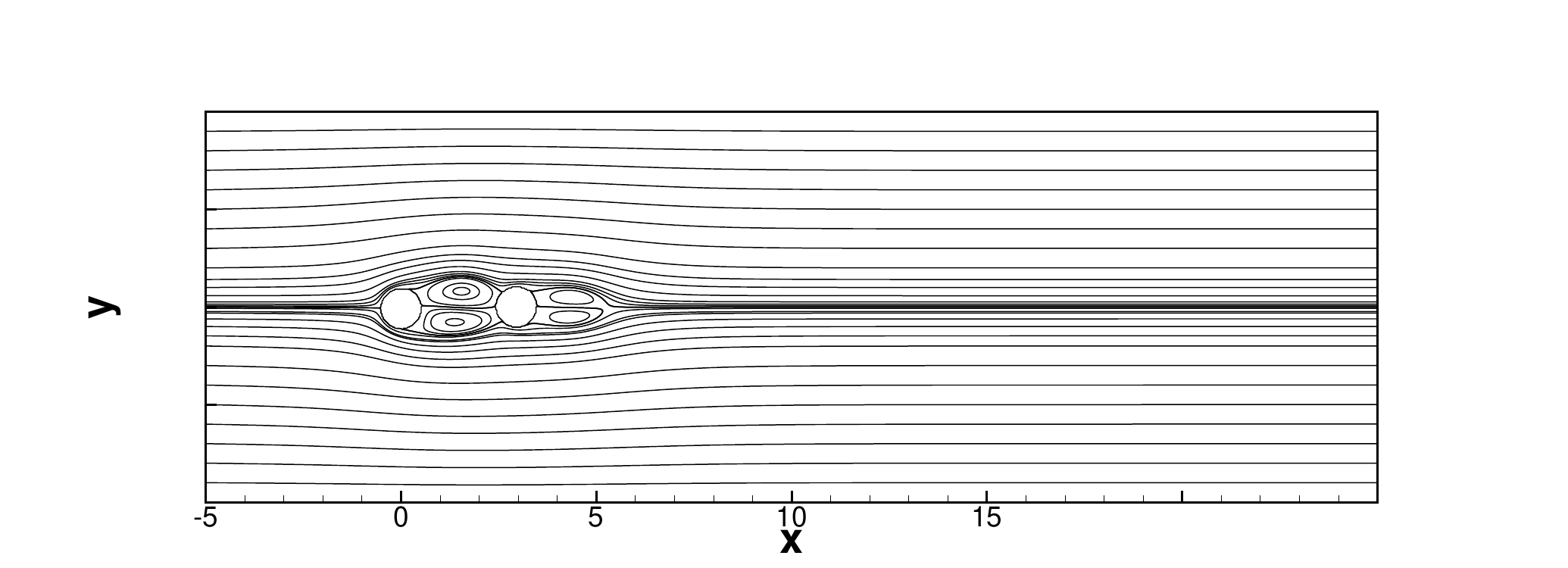,width=1.0\linewidth}
\end{minipage}
\hspace{-0.75cm}
\begin{minipage}[b]{.45\linewidth}
\centering\psfig{file=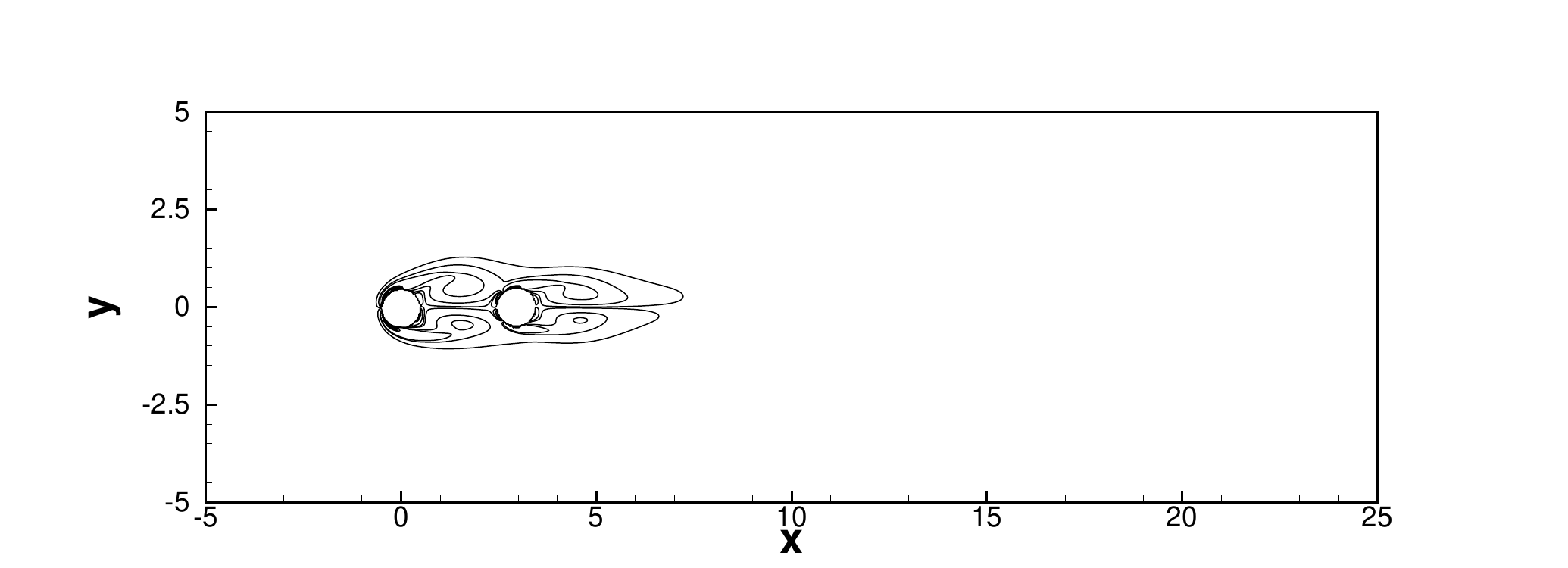,width=1.0\linewidth}
\end{minipage}(b)
  \begin{minipage}[b]{.45\linewidth}   \
\hfill \centering\psfig{file=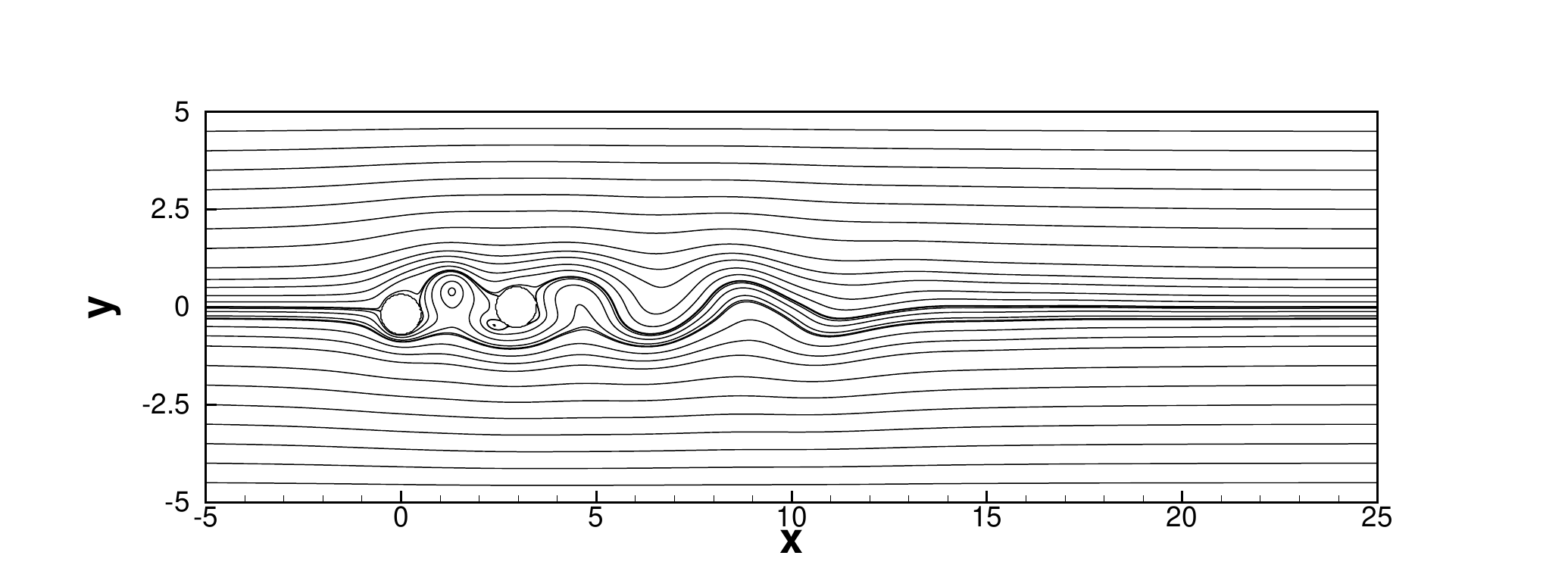,width=1.0\linewidth}
\end{minipage}
\hspace{-0.75cm}
\begin{minipage}[b]{.45\linewidth}   \
\centering\psfig{file=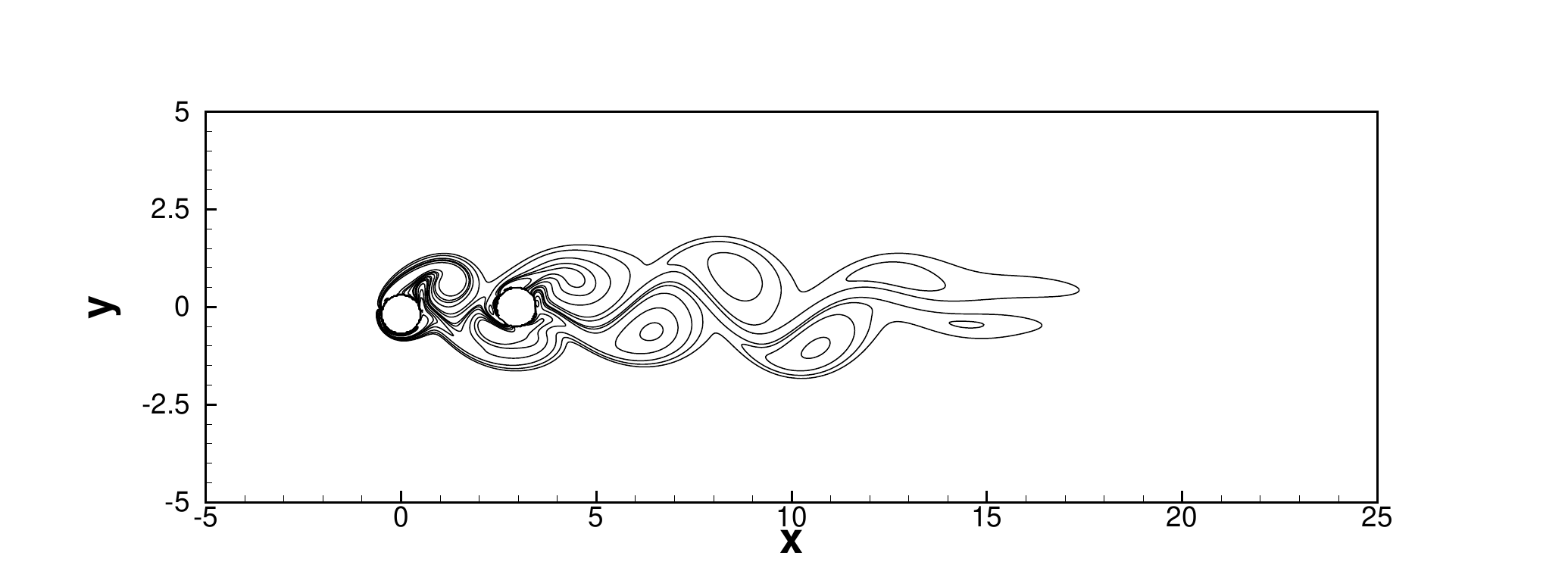,width=1.0\linewidth}
\end{minipage}(b)
\begin{minipage}[b]{.45\linewidth}
\hfill \centering\psfig{file=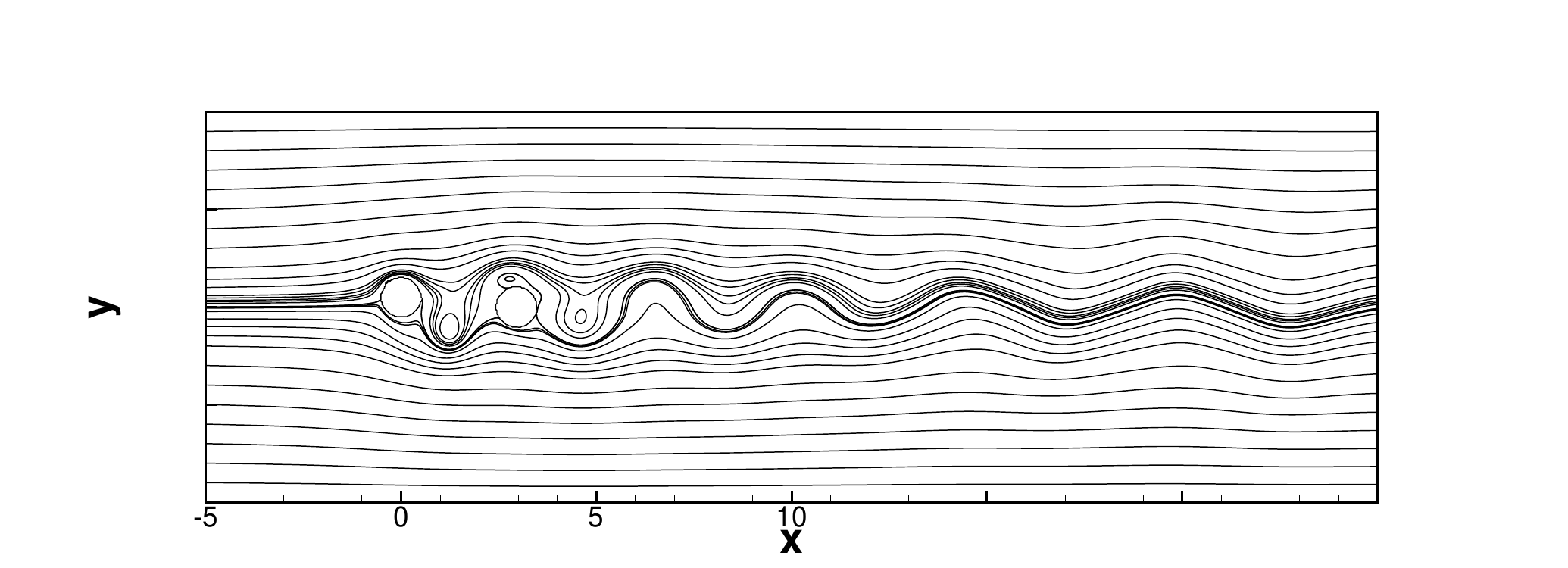,width=1.0\linewidth}
\end{minipage}
\hspace{-0.75cm}
\begin{minipage}[b]{.45\linewidth}   \
\centering\psfig{file=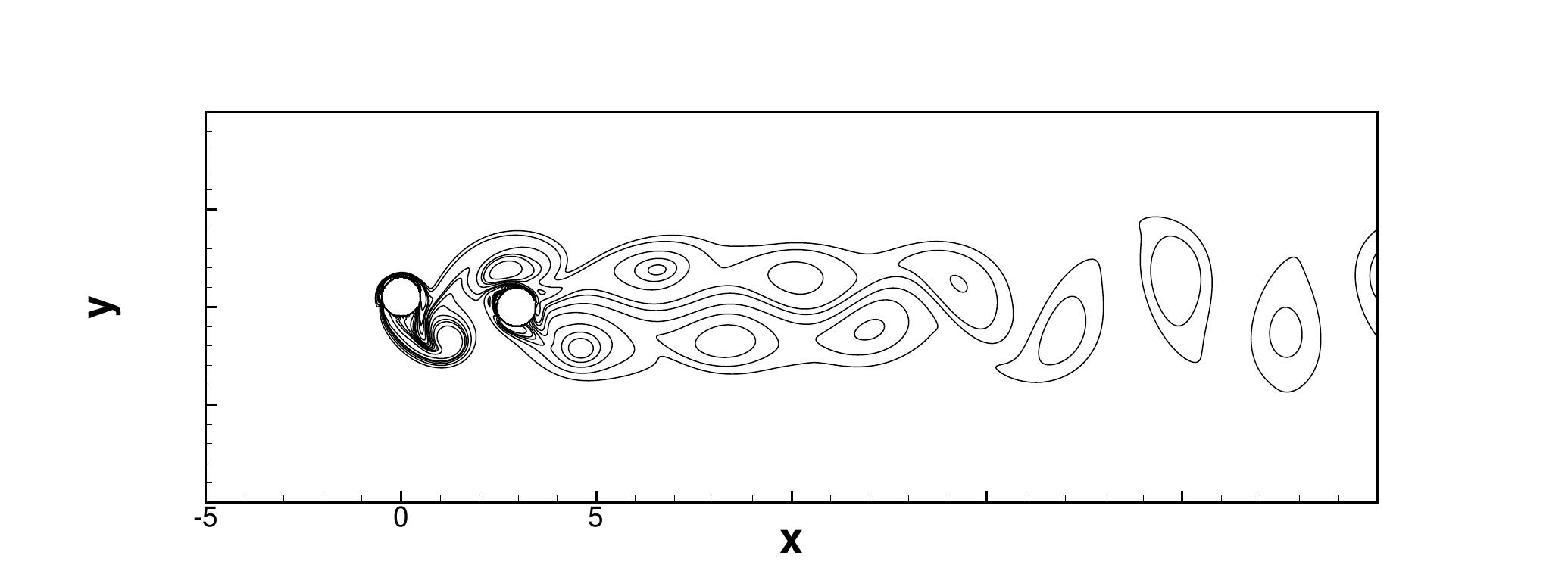,width=1.0\linewidth}
\end{minipage}(d)
\begin{minipage}[b]{.45\linewidth}
\hfill \centering\psfig{file=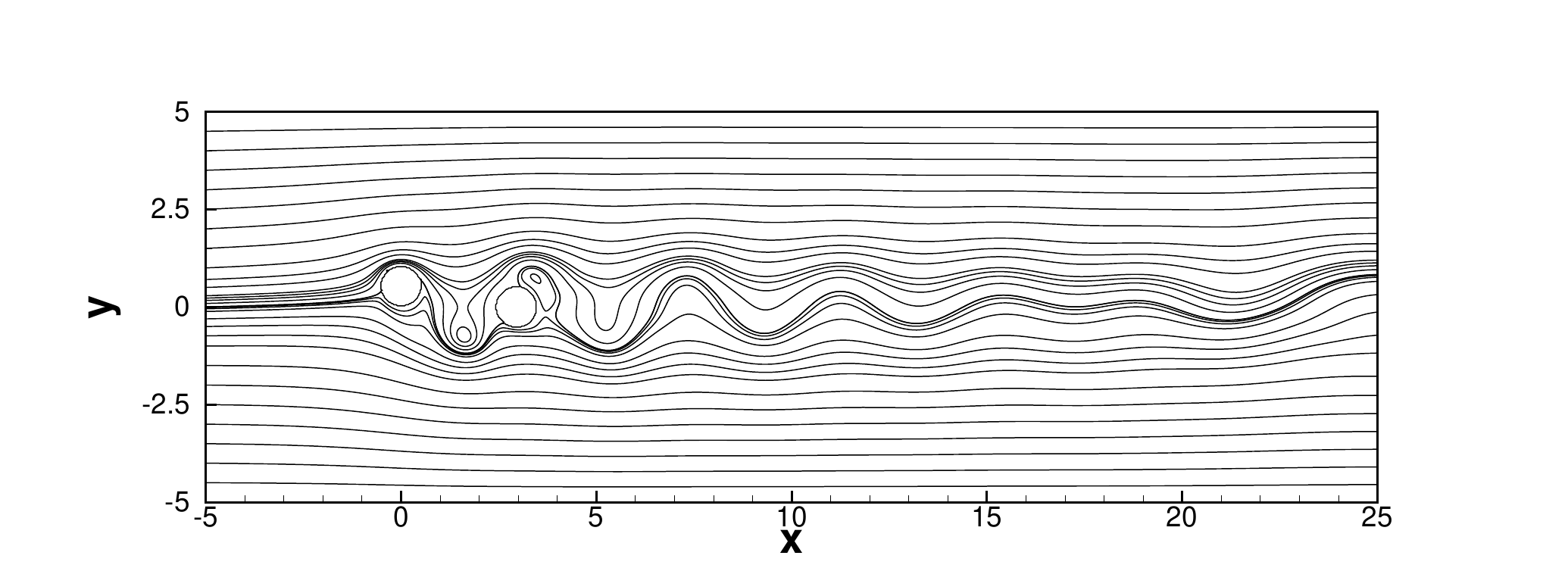,width=1.0\linewidth}
\end{minipage}
\hspace{-0.75cm}
\begin{minipage}[b]{.45\linewidth}   \
\centering\psfig{file=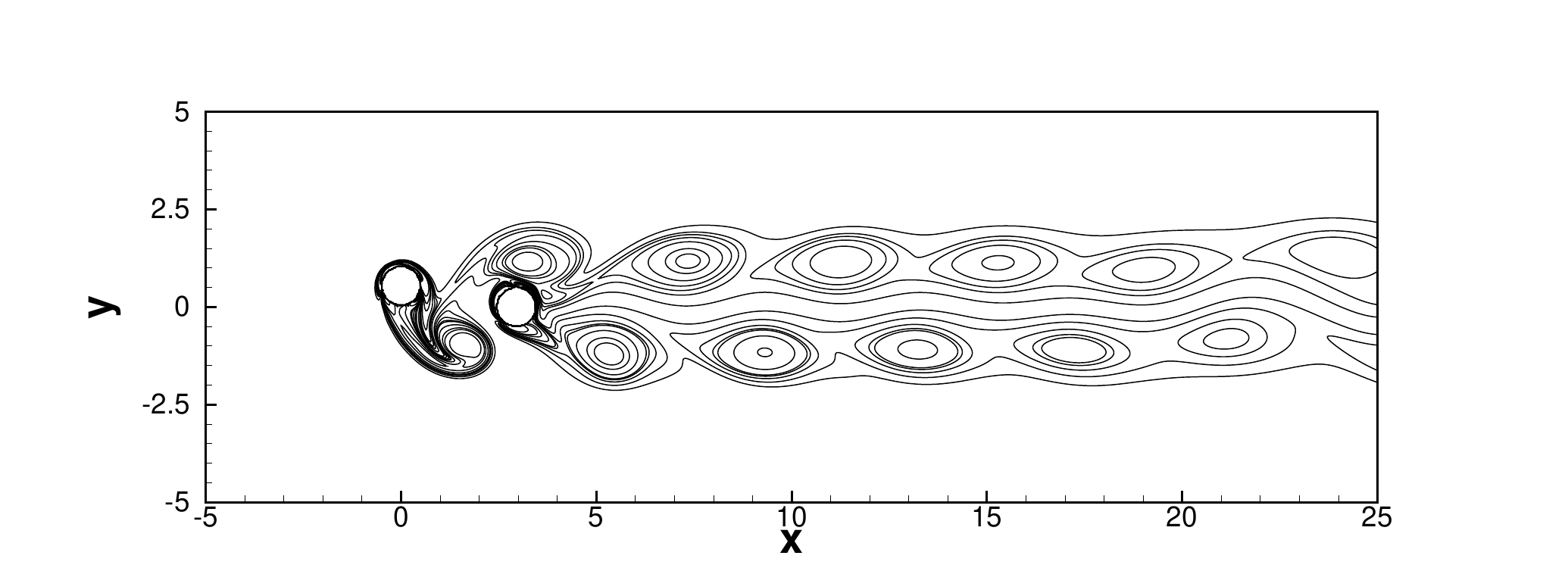,width=1.0\linewidth}
\end{minipage}(e)
\begin{minipage}[b]{.45\linewidth}
\hfill \centering\psfig{file=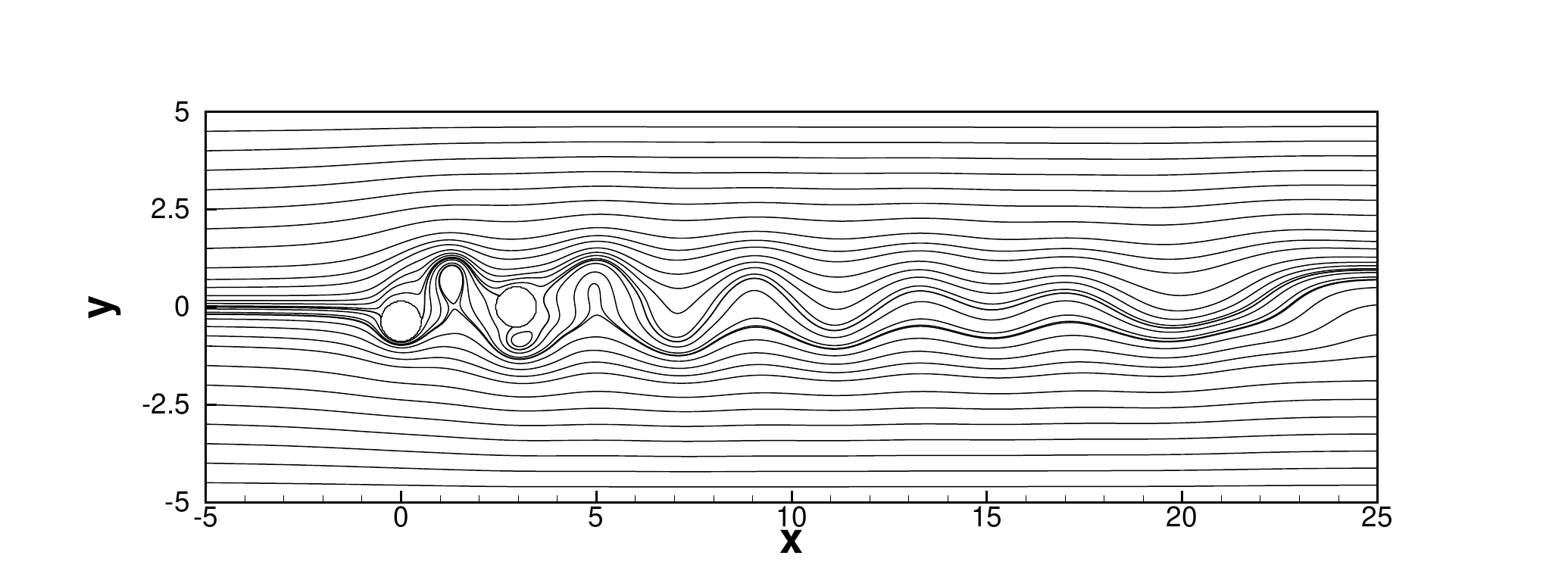,width=1.0\linewidth}
\end{minipage}
\hspace{-0.75cm}
\begin{minipage}[b]{.45\linewidth}
\centering\psfig{file=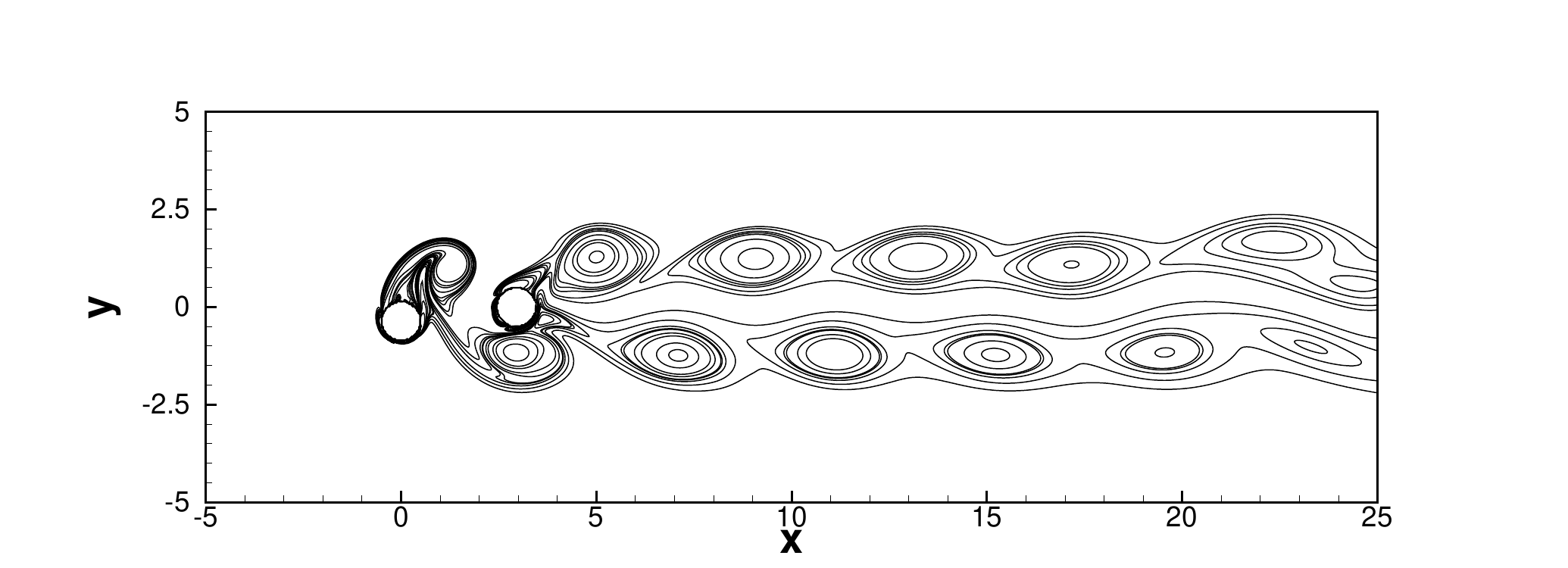,width=1.0\linewidth}
\end{minipage}(f)
\caption{\sl {Streamlines (left) and vorticity contours (right) for flow past two tandem circular cylinders at: (a) $t = 1$, (b) $t = 5$, (c) $t = 20$, (d) $t = 50$ (e) $t = 100$,  and (f) $t = 300$ for Test Case 5.}}\label{tan_sf_vt}
\end{figure}

\begin{figure}[!h]
 \includegraphics[scale=0.8]{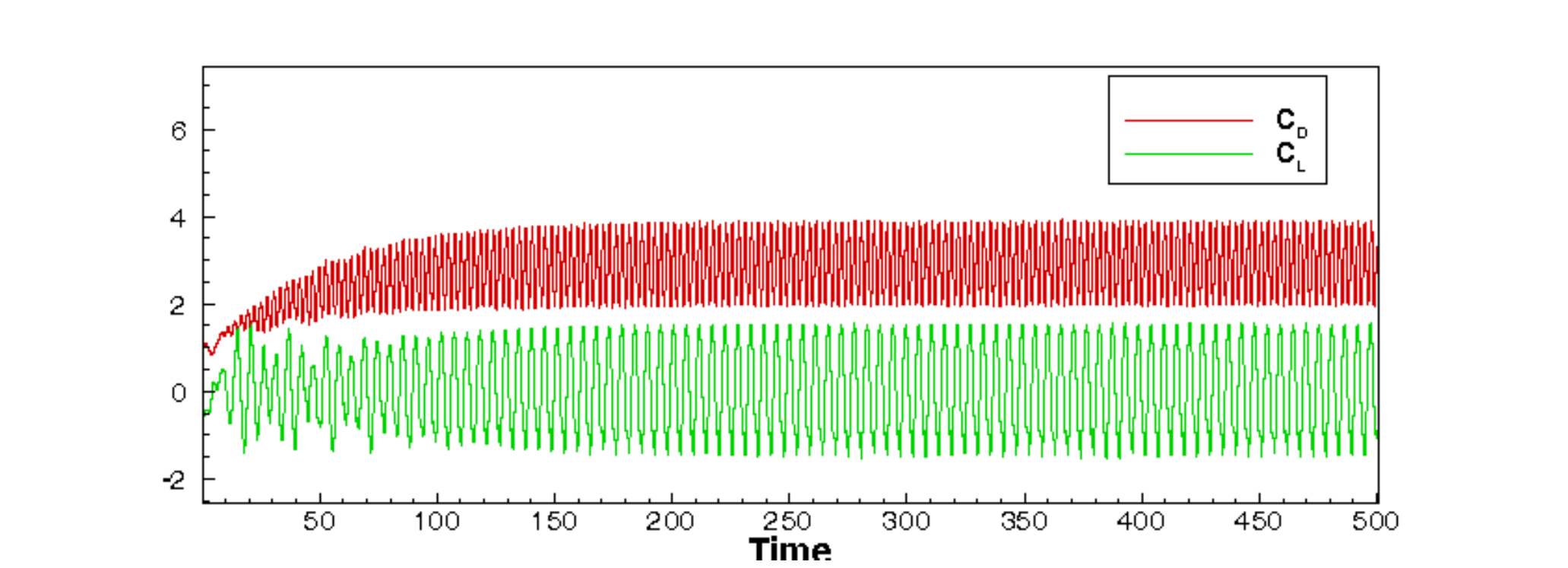} 
         \caption{{\sl History of drag and lift coefficients of the tandem cylinders for case 5.} }
\label{dl_tandem}
\end{figure}
In figure \ref{tan_sf_vt}, we show the evolution of the flow by plotting the streamlines and vorticity contours for this flow for $Re=100$ at time stations $t=1,\;5,\; 20,\;50,\;100\; {\rm and}\;300$. It is evident that the wake behind the cylinders settles into a periodic motion eventually, which is also confirmed by figure \ref{dl_tandem}, where we show the time history of the drag and lift coefficients of the upstream cylinder. A {\bf 2S} shedding mode, typical of this kind of flow \cite{behara2019flow} is obvious here. This can be further observed in the accompanying video "tandem.avi" where the flow evolution during $0 \leq t \leq 200$ is shown. It is heartening to note that without the inclusion of any external circles embedded into the figures, the interfaces have been captured very smoothly by our immersed interface approach on a relatively coarse grid with step-length $h=l=0.0295$. 
\begin{figure}[!h]
 \begin{minipage}[b]{.45\linewidth}  
\includegraphics[scale=0.4]{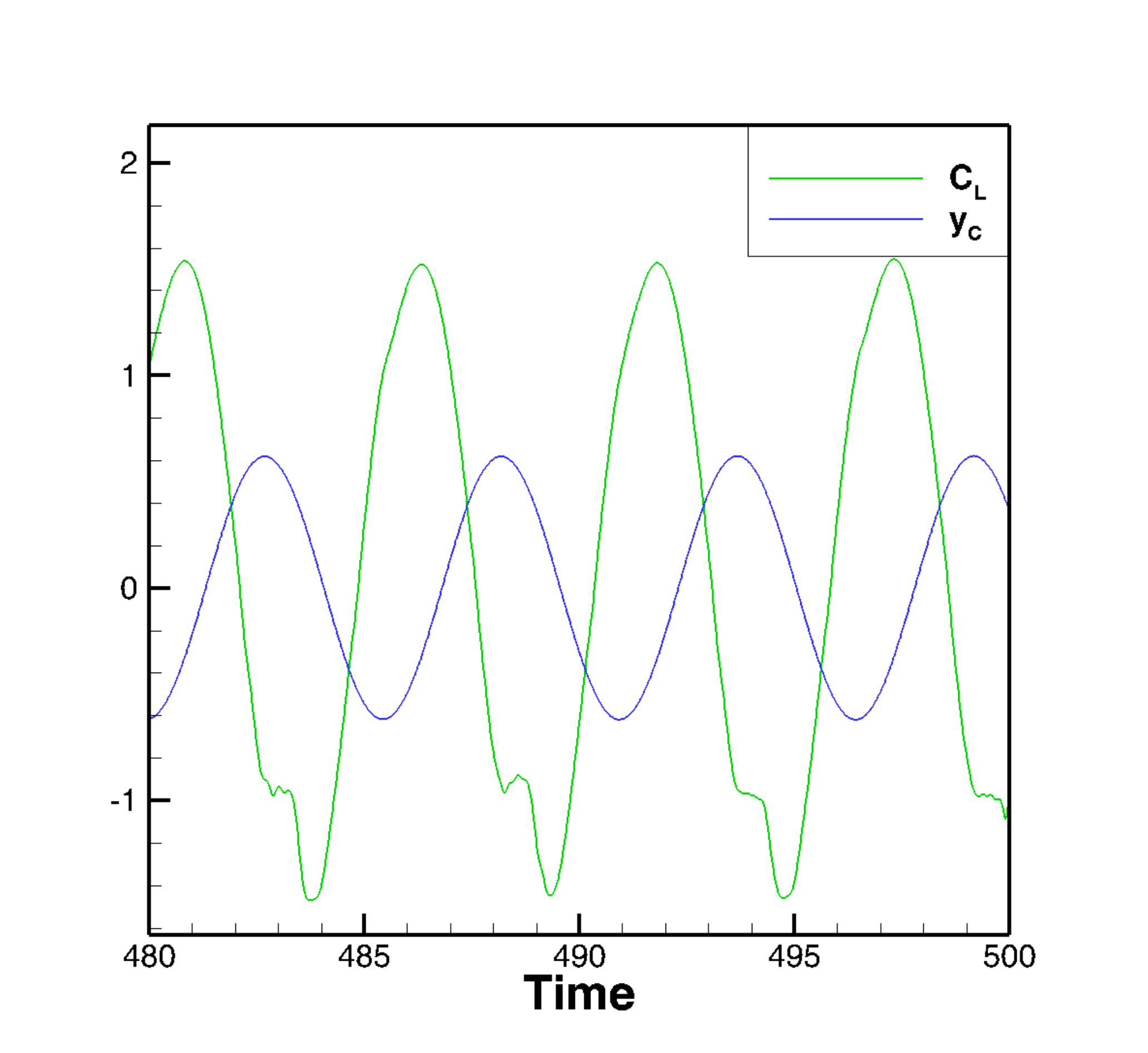} 
\centering (a) 
\end{minipage}           
\begin{minipage}[b]{.45\linewidth}
\includegraphics[scale=0.4]{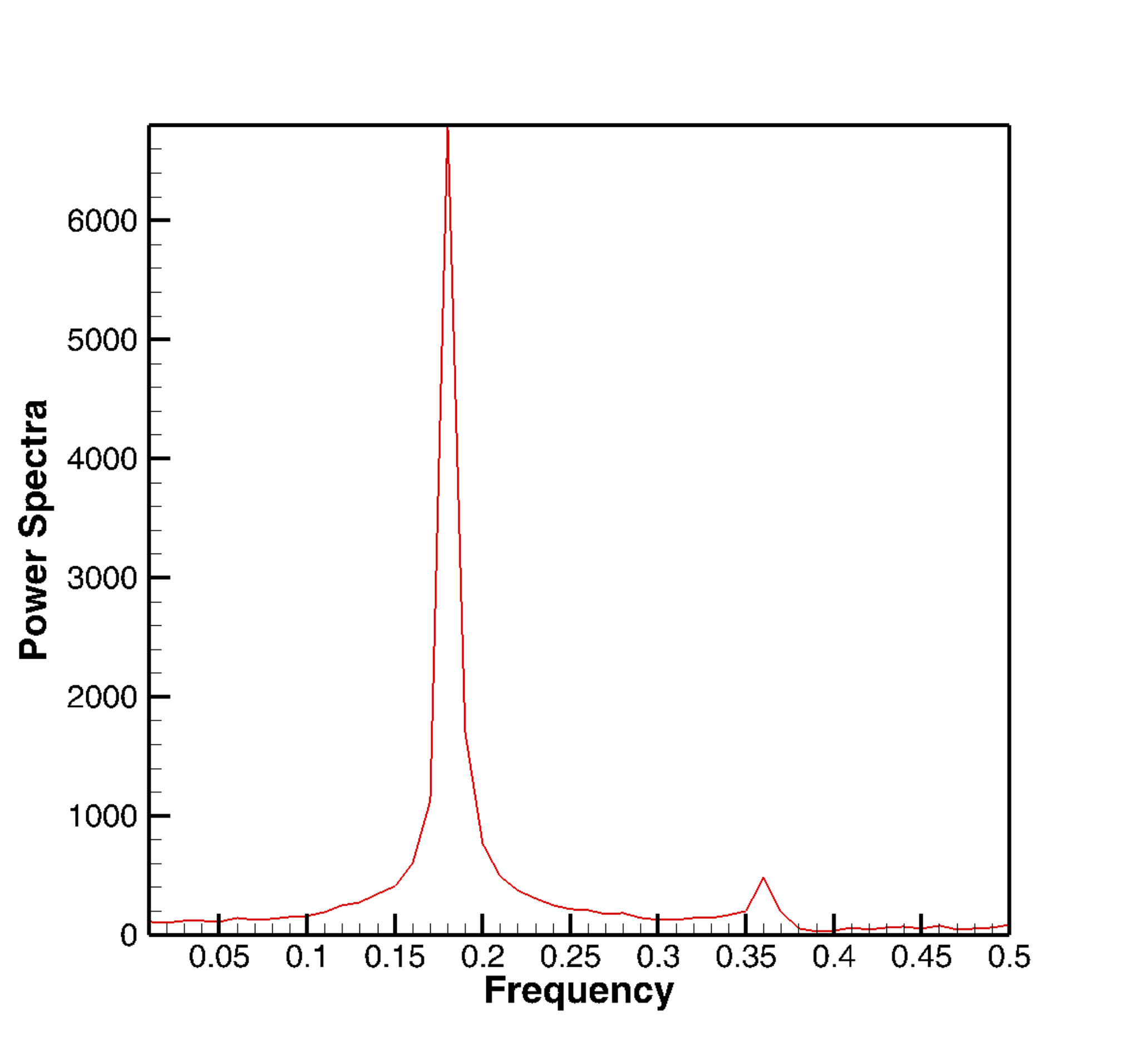}
\centering (b) 
\end{minipage} 
 \caption{{\sl (a) Time history of the  $y$-coordinate of the oscillating cylinder and lift coefficient, and (b) Power spectra of the lift coefficient displaying the Strouhal number for Test Case 5.} }
\label{lock_tandem}
\end{figure}

In figure \ref{lock_tandem}(a), we show time history of the displacement of the oscillating cylinder side by side with the time history of the corresponding lift coefficient for a very short time range $480 \leq t \leq 500$. These figures demonstrate that the primary vortex shedding frequency and the frequency of vibration of the cylinder is almost equal. Figure \ref{lock_tandem}(b) reconfirms the fact as the Strouhal number $St$ (as defined in section \ref{st_cyl}) computation from the Fast Fourier Transform of the lift coefficient history results in a value of $0.1801$, which is extremely close to the frequency $0.182$ of the upstream cylinder. Thus we conclude that the flow situation under consideration here exhibits the characteristics of a  {\bf lock-on} regime. 
\section{Conclusion}
The current work is concerned with the development of a hybrid explicit jump immersed interface approach in conjunction with a higher order compact (HOC) scheme for simulating transient complex flows on Cartesian grids. Originally developed for parabolic equations with discontinuities in the solutions, source terms and the coefficients across the interface, this approach was seen to easily accommodate the N-S equations for simulating flow past bluff bodies immersed in the flow. $\psi$-$\zeta$ formulation of the N-S equations for incompressible viscous flows has been utilized for this purpose. A novel strategy has been adopted for the jump conditions at the irregular points across the interface using Lagrangian interpolation on a Cartesian grid. A compact and concise form of the matrix equations resulting from the discretization of the  $\psi$ and $\zeta$ equations have also been provided.

Firstly a parabolic problem having a known analytical solution is solved in order to establish the spatial rate of convergence of the proposed approach. Our approach was seen to reduce magnitude of the error with a much faster decay rate of $O(h^4)$ in comparison to other existing methods, thus establishing the theoretical rate of convergence in the way. Next, it was employed to simulate several complex fluid flow problems past bluff bodies having real life applications, including flows involving multiple and moving bodies. This includes the flow past a stationary circular and a twenty-four edge cactus cylinder, flows past two tandem cylinders, where in one situation both are fixed and in another, one of them is transversely oscillating in the flow. Except for the stationary circular cylinder cases, opposed to most of the earlier computations which were performed by either finite volume or finite element approach in extremely finer grids , our approach accomplishes the same in relatively coarse grids in  FD set-up. Moreover, to the best of our knowledge, the tandem cylinder case, where one cylinder is oscillating with variable amplitudes along with the stationary cactus, have been tackled for the first time by such an approach employing the $\psi$-$\zeta$ formulation of the N-S equations against the primitive variable formulation in earlier simulations. Furthermore, in the process, we also provide elaborate description of the steps to compute drag and lift coefficients in multiply connected domains. 

In all cases, our computed solutions are extremely close to existing numerical and experimental results. Thus, apart from exemplifying the accuracy and the robustness of the proposed approach, our simulations aptly demonstrates its ability in handling complicated geometries and varied flow situations. Currently we are working on the expansion of the proposed approach to problems involving heat and mass transfer, and multiphase flows.

\section{Appendix}
We have $\vec{x}= x \hat{i} + y \hat{j}$, $\vec{u}= u \hat{i} + v \hat{j}$ and $\vec{u_{s}}= u_{s} \hat{i} + v_{s} \hat{j}$. Note that $\vec{\omega}=  \zeta\hat{k}$ here so that
\begin{equation}
\lbrace \vec{ \omega} (\vec{x} \times \vec{u})\rbrace . \hat{n} ds= (\vec{x} \times \vec{u}) \lbrace \zeta \hat{k} . dy \hat{i} - dx \hat{j}\rbrace = \vec{0}. \label{m4}
\end{equation}
Also, 
\begin{equation}
\frac{d}{dt} \int_{V(t)} \vec{u} dV  = \int_{V(t)} \frac{\partial \vec{u}}{\partial t} dV +\oint_{S(t)} \hat{n}. \vec{u_s} \vec{u} \,ds + \oint_{S_b(t)} \hat{n}. \vec{u_s} \vec{u} \,ds.
 \end{equation}
 If the surface is fixed, the surface integral over $S$ reduces to zero. Likewise, for a fixed volume, viz., when $V(t)$ is independent of time, the surface integral over $S_b(t)$ vanishes. A sufficient (although the necessary) condition for this surface integral to vanish is for $\vec{u_s} . \hat{n}$ to be equal to zero. With all these assumptions, \eqref{m1} reduces to 
\begin{equation}
\vec{F}=-\frac{d}{dt} \int_{V(t)} \vec{u} dV + \oint_{S(t)} \hat{n}. \gamma_{mom} \,ds.  \label{m5}
\end{equation}
for control volume enclosing solid bodies as the last term in \eqref{m1} vanishes for such bodies, under these assumptions.
\begin{equation}
\gamma_{mom}= \frac{\rho}{2} |\vec{u}|^{2} \mathbf{I} + \rho\left [(\vec{u_{s}} - \vec{u}) \vec{u}- \vec{u} (\vec{x} \times \vec{\zeta}) \right] - \rho\left[ \left( \vec{x}. \frac{\partial \vec{u}}{\partial t} \mathbf{I} - \vec{x} \frac{\partial \vec{u}}{\partial t} \right) \right] + \left[ \vec{x}. (\nabla . \mathbf{T}) \mathbf{I} - \vec{x} (\nabla. \mathbf{T}) \right] + \mathbf{T}. \label{m6}
\end{equation} 
We will evaluate the second term on the RHS of \eqref{m6} using each of terms on the RHS of \eqref{m6} (for the cases considered in our study).
\begin{equation}
 \frac{1}{2} |\vec{u}|^{2} \mathbf{I} . \hat{n} ds= \frac{1}{2} (u^2+v^2) dy \hat{i} - \frac{1}{2} (u^2+v^2) dx \hat{j}.\label{m7}
\end{equation}
\begin{equation}
 (\vec{u_{s}} - \vec{u}) \vec{u} . \hat{n} ds= \lbrace (u_{s}-u)(u dy - v dx) \hat{i} + (v_{s}-v) (u dy - v dx) \hat{j}  \rbrace.   \label{m8}
\end{equation}
\begin{equation}
- \lbrace  \vec{u} (\vec{x} \times \vec{\omega}) \rbrace  . \hat{n} ds=  (-u dy + v dx) (y \zeta \hat{i} - x \zeta \hat{j}).\label{m9}
\end{equation}
\begin{eqnarray}
\left[ \left( \vec{x}. \frac{\partial \vec{u}}{\partial t} \mathbf{I} - \vec{x} \frac{\partial \vec{u}}{\partial t} \right) \right] . \hat{n} ds= & (dy \hat{i} - dx \hat{j}) \left( x \frac{\partial u}{\partial t} + y \frac{\partial v}{\partial t}\right) + (x dy - y dx) \left(  \frac{\partial u}{\partial t} \hat{i} + \frac{\partial v}{\partial t} \hat{j}\right), \nonumber \\
 =& y \left( dx \frac{\partial u}{\partial t} + dy \frac{\partial v}{\partial t}\right) \hat{i} - x \left( dx \frac{\partial u}{\partial t} + dy \frac{\partial v}{\partial t}\right) \hat{j}. \label{m10}
\end{eqnarray}
Now $\vec{F} = F_D \hat{i} + F_L \hat{j}$, we get
\begin{equation}
F_D= - \iint_V \left(\frac{\partial u }{\partial t} + u \frac{\partial u}{\partial x} + v \frac{\partial v}{\partial y}\right) \,dx\,dy + \oint_{S(t)}  \left(-v (u_{s} -u)+yv \zeta -y \frac{\partial v}{\partial t}\right) dx +  \left( \frac{1}{2} (u^2+v^2) + u (u_{s} - u) -yu \zeta - y \frac{\partial v}{\partial t} \right) dy \,ds 
\end{equation}
\begin{equation}
F_L= - \iint_V \left(\frac{\partial v }{\partial t} + u \frac{\partial v}{\partial x} + y \frac{\partial v}{\partial y}\right) \,dx\,dy + \oint_{S(t)}  \left( -\frac{1}{2} (u^2+v^2) - v (v_{s} - v) -xv \zeta - x \frac{\partial u}{\partial t} \right) dx +  \left(u (v_{s} -v)+xu \zeta +x \frac{\partial u}{\partial t}\right) dy  \,ds 
\end{equation}
and when the body is stationary then $u_{s}=0$.
\begin{equation}
\oint_{S} \lbrace \vec{x} . (\nabla . \mathbf{T}) \mathbf{I} - \vec{x} (\nabla . \mathbf{T}) + \mathbf{T}\rbrace \hat{n} \,ds
\end{equation}
\begin{equation}
=\mu \oint_{S(t)} (x \hat{i} + y \hat{j}). (\nabla^{2} u \hat{i} + \nabla^2 v \hat{j}) (dy \hat{i} - dx \hat{j}) - (x dx - y dx) (\nabla^2 u \hat{i} + \nabla^2 v \hat{j}) + \left( 2 \frac{\partial u}{\partial x} + \frac{\partial u}{\partial y} + \frac{\partial v }{\partial x}  \right) dy \hat{i} -  \left(  \frac{\partial u}{\partial y} + \frac{\partial v}{\partial x} + 2 \frac{\partial v }{\partial y}  \right) dx \hat{j} \,ds ,
\end{equation}
$\left( \textnormal{where} \quad \nabla. \mathbf{T}= \mu (\nabla^2 u \hat{i} + \nabla^2 v \hat{j}) \quad \textnormal{and} \quad \mathbf{T}=\mu \left( 2 \frac{\partial u}{\partial x} + \frac{\partial u}{\partial y} + \frac{\partial v }{\partial x}  \right) \hat{i} + \mu \left(  \frac{\partial u}{\partial y} + \frac{\partial v}{\partial x} + 2 \frac{\partial v }{\partial y}  \right) \hat{j} \right).$
\begin{eqnarray}
=&\mu \oint_{S(t)}  \lbrace ( x \nabla^{2} u  + y \nabla^2 v ) dy\rbrace \hat{i} -& \lbrace( x \nabla^{2} u  + y \nabla^2 v ) dx \rbrace \hat{j} - ( x \nabla^{2} u  dy - y \nabla^2 u dx)  \hat{i} - ( x \nabla^{2} v  dy - y \nabla^2 v dx)  \hat{j}  \nonumber\\
+& \left( 2 \frac{\partial u}{\partial x} + \frac{\partial u}{\partial y} + \frac{\partial v }{\partial x}  \right) dy \hat{i} -&  \left(  \frac{\partial u}{\partial y} + \frac{\partial v}{\partial x} + 2 \frac{\partial v }{\partial y}  \right) dx \hat{j} \,ds,
\end{eqnarray}
\begin{equation}
=\mu  \oint_{S(t)}  \lbrace y \nabla^{2} u dx + \left( y \nabla^2 v + 2 \frac{\partial u}{\partial x} + \frac{\partial u}{\partial y} + \frac{\partial v }{\partial x}    dy \right)  \rbrace \hat{i} - \lbrace x \nabla^{2} v dy + \left( x \nabla^2 u +  \frac{\partial u}{\partial y} + \frac{\partial v}{\partial x} + 2 \frac{\partial v }{\partial y}    dx \right)  \rbrace \hat{j} \,ds.
\end{equation}

\bibliographystyle{amsplain}	

\end{document}